%
\newif\ifloadreferences\loadreferencestrue
%
%
%
%
%
\let\myfrac=\frac%
\input eplain %
\let\frac=\myfrac%
\let\myfootnote=\footnote%
\input amstex \input epsf %
\let\footnote=\myfootnote%
%
%
\loadeufm\loadmsam\loadmsbm\message{symbol names}\UseAMSsymbols\message{,}%
%
\font\myfontdefault=cmr10%
\newif\ifmakebiblio%
\newif\ifinappendices%
\newif\ifundefinedreferences%
\newif\ifchangedreferences%
\makebibliofalse%
\undefinedreferencesfalse%
\changedreferencesfalse%
%
%
%
%
%
\def\setcatcodes{\catcode`\!=0 \catcode`\\=11}%
{\global\let\noe=\noexpand%
\catcode`\@=11 \catcode`\_=11 \setcatcodes%
!global!def!_@@internal@@makeref#1{%
!global!expandafter!def!csname #1ref!endcsname##1{%
!csname _@#1@##1!endcsname%
!expandafter!ifx!csname _@#1@##1!endcsname!relax%
    !write16{#1 ##1 not defined - run saving references}%
    !undefinedreferencestrue%
!fi}}%
!global!def!_@@internal@@makelabel#1{%
!global!expandafter!def!csname #1label!endcsname##1{%
!edef!temptoken{!csname #1info!endcsname}%
!ifloadreferences%
!expandafter!ifx!csname _@#1@##1!endcsname!relax%
!write16{#1 ##1 not hitherto defined - rerun saving references}%
!changedreferencestrue%
!else%
!expandafter!ifx!csname _@#1@##1!endcsname!temptoken%
!else%
!write16{#1 ##1 reference has changed - rerun saving references}%
!changedreferencestrue%
!fi%
!fi%
!else%
!expandafter!edef!csname _@#1@##1!endcsname{!temptoken}%
!edef!textoutput{!write!references{\global\def\_@#1@##1{!temptoken}}}%
!textoutput%
!fi}}%
!global!def!makecounter#1{!_@@internal@@makelabel{#1}!_@@internal@@makeref{#1}}%
!unsetcatcodes%
}
%
%
%
%
%
\def\turnintolatin#1{\ifcase #1 _\or i\or ii\or iii\or iv\or v\or vi\or vii\or viii\or ix\or x\or xi\or xii\or xiii\or xiv\or xv\or xvi\or xvii\or xviii\or xix\or xx\or xxi\or xxii\or xxiii\or xxiv\or xxv\or xxvi\fi}%
\def\alphanum#1{\ifcase #1 _\or A\or B\or C\or D\or E\or F\or G\or H\or I\or J\or K\or L\or M\or N\or O\or P\or Q\or R\or S\or T\or U\or V\or W\or X\or Y\or Z\fi}%
\newwrite\references%
\ifloadreferences{\catcode`\@=11 \catcode`\_=11 \global\def\_@citation@BarBegZegh{1}
\global\def\_@citation@Bianchi{2}
\global\def\_@citation@BonMonSchI{3}
\global\def\_@citation@BonMonSchII{4}
\global\def\_@citation@Clifford{5}
\global\def\_@citation@CliffordI{6}
\global\def\_@citation@CliffordII{7}
\global\def\_@citation@CruForGad{8}
\global\def\_@citation@DajNom{9}
\global\def\_@citation@EnomotoI{10}
\global\def\_@citation@EnomotoII{11}
\global\def\_@citation@EnoKitWei{12}
\global\def\_@citation@FillSmi{13}
\global\def\_@citation@GalMarMil{14}
\global\def\_@citation@GalMirI{15}
\global\def\_@citation@GalMirII{16}
\global\def\_@citation@Garling{17}
\global\def\_@citation@Goldstein{18}
\global\def\_@citation@HarveyLawson{19}
\global\def\_@citation@HitchinII{20}
\global\def\_@citation@HitchinI{21}
\global\def\_@citation@KitigawaII{22}
\global\def\_@citation@KitigawaI{23}
\global\def\_@citation@KokUmeYam{24}
\global\def\_@citation@KokRossSajUmeYam{25}
\global\def\_@citation@LabIII{26}
\global\def\_@citation@LabI{27}
\global\def\_@citation@LabII{28}
\global\def\_@citation@LabIV{29}
\global\def\_@citation@Mira{30}
\global\def\_@citation@Lounesto{31}
\global\def\_@citation@LychRub{32}
\global\def\_@citation@LychRubChekI{33}
\global\def\_@citation@McDuffSalamon{34}
\global\def\_@citation@Milnor{35}
\global\def\_@citation@Pinkall{36}
\global\def\_@citation@Rubtsov{37}
\global\def\_@citation@SajUmeYam{38}
\global\def\_@citation@Sasaki{39}
\global\def\_@citation@SasakiNomizu{40}
\global\def\_@citation@Sch{41}
\global\def\_@citation@SmiSLC{42}
\global\def\_@citation@SmiEW{43}
\global\def\_@citation@SmiMS{44}
\global\def\_@citation@SmiAS{45}
\global\def\_@citation@SmiQMAK{46}
\global\def\_@citation@SmiQCKS{47}
\global\def\_@citation@Spivak{48}
\global\def\_@citation@Tamburelli{49}
\global\def\_@citation@Volkert{50}
\global\def\_@citation@Weiner{51}
\global\def\_@citation@Yau{52}
\global\def\_@head@Introduction{1}
\global\def\_@subhead@Introduction{1.1}
\global\def\_@subhead@CECSurfaces{1.2}
\global\def\_@subhead@AFinalWord{1.3}
\global\def\_@subhead@Acknowledgements{1.4}
\global\def\_@head@CliffordStructures{2}
\global\def\_@subhead@OverviewCliffordStructures{2.1}
\global\def\_@subhead@PseudoInvolutions{2.2}
\global\def\_@proc@PseudoInvolutionsHaveEigenspaces{2.2.1}
\global\def\_@eqn@JordanDecompositionOfIdentity{\relax \unhbox \voidb@x \hbox {{\relax \tenrm (2.1)}}}
\global\def\_@proc@InvariantPlanes{2.2.2}
\global\def\_@rmk@InvariantPlanes{\relax \unhbox \voidb@x \hbox {2.2.1}}
\global\def\_@subhead@CliffordAlgebras{2.3}
\global\def\_@eqn@FormOfX{\relax \unhbox \voidb@x \hbox {{\relax \tenrm (2.2)}}}
\global\def\_@eqn@EigenspacesOfMx{\relax \unhbox \voidb@x \hbox {{\relax \tenrm (2.3)}}}
\global\def\_@proc@NegativeUnitImaginaryElementsAreBalanced{2.3.1}
\global\def\_@proc@InvariantPlanesAreBiLagrangian{2.3.2}
\global\def\_@subhead@InvariantPlanes{2.4}
\global\def\_@proc@OrthogonalityAndNullityOfEigenspaces{2.4.1}
\global\def\_@proc@ExceptionalPlanesI{2.4.2}
\global\def\_@rmk@ExceptionalPlanesI{\relax \unhbox \voidb@x \hbox {2.4.1}}
\global\def\_@proc@ExceptionalPlanesII{2.4.3}
\global\def\_@rmk@ExceptionalPlanesII{\relax \unhbox \voidb@x \hbox {2.4.2}}
\global\def\_@proc@PrincipalDirectionsI{2.4.4}
\global\def\_@proc@PrincipalDirectionsII{2.4.5}
\global\def\_@rmk@PrincipalDirectionsII{\relax \unhbox \voidb@x \hbox {2.4.3}}
\global\def\_@proc@PrincipalPlanes{2.4.6}
\global\def\_@rmk@PrincipalPlanes{\relax \unhbox \voidb@x \hbox {2.4.4}}
\global\def\_@head@CliffordStructuresAndExtrinsicCurvature{3}
\global\def\_@subhead@OverviewLegendrianCliffordStructures{3.1}
\global\def\_@subhead@AnExplicitCliffordStructure{3.2}
\global\def\_@eqn@DefinitionOfA{\relax \unhbox \voidb@x \hbox {{\relax \tenrm (3.1)}}}
\global\def\_@eqn@DescriptionOfFix{\relax \unhbox \voidb@x \hbox {{\relax \tenrm (3.2)}}}
\global\def\_@proc@DescriptionOfFix{3.2.1}
\global\def\_@proc@ExistenceOfCliffordIsomorphism{3.2.2}
\global\def\_@eqn@ExplicitFormulaOfIJK{\relax \unhbox \voidb@x \hbox {{\relax \tenrm (3.3)}}}
\global\def\_@subhead@TheCliffordBundle{3.3}
\global\def\_@eqn@OmegaKIsFlat{\relax \unhbox \voidb@x \hbox {{\relax \tenrm (3.4)}}}
\global\def\_@proc@OmegaKIsFlat{3.3.1}
\global\def\_@eqn@IJKAreCovariantConstant{\relax \unhbox \voidb@x \hbox {{\relax \tenrm (3.5)}}}
\global\def\_@proc@IJKAreCovariantConstant{3.3.2}
\global\def\_@subhead@IntroducingBilegendrianSurfaces{3.4}
\global\def\_@proc@SymmetryOfCubicForm{3.4.1}
\global\def\_@eqn@CSymmetryOfJ{\relax \unhbox \voidb@x \hbox {{\relax \tenrm (3.6)}}}
\global\def\_@proc@CSymmetryOfJ{3.4.2}
\global\def\_@eqn@VanishingValuesOfCCF{\relax \unhbox \voidb@x \hbox {{\relax \tenrm (3.7)}}}
\global\def\_@proc@VanishingValuesOfCCF{3.4.3}
\global\def\_@proc@FlatMetric{3.4.4}
\global\def\_@subhead@SurfacesOfConstantPositiveCurvature{3.5}
\global\def\_@eqn@DerivativeOfGaussLift{\relax \unhbox \voidb@x \hbox {{\relax \tenrm (3.8)}}}
\global\def\_@proc@DerivativeOfKGaussLift{3.5.1}
\global\def\_@proc@KGaussLiftIsInvariant{3.5.2}
\global\def\_@eqn@FlatMetricOfKSurface{\relax \unhbox \voidb@x \hbox {{\relax \tenrm (3.9)}}}
\global\def\_@proc@FlatMetricOfKSurface{3.5.3}
\global\def\_@rmk@FlatMetricOfKSurface{\relax \unhbox \voidb@x \hbox {3.5.1}}
\global\def\_@subhead@SurfacesOfConstantNegativeCurvature{3.6}
\global\def\_@proc@KGaussLiftIsInvariantII{3.6.1}
\global\def\_@rmk@KGaussLiftIsInvariantII{\relax \unhbox \voidb@x \hbox {3.6.1}}
\global\def\_@eqn@FlatMetricOfKSurfaceB{\relax \unhbox \voidb@x \hbox {{\relax \tenrm (3.10)}}}
\global\def\_@proc@FlatMetricOfKSurfaceB{3.6.2}
\global\def\_@rmk@FlatMetricOfKSurfaceB{\relax \unhbox \voidb@x \hbox {3.6.2}}
\global\def\_@eqn@FFSNegCurv{\relax \unhbox \voidb@x \hbox {{\relax \tenrm (3.11)}}}
\global\def\_@subhead@HilbertsTheorem{3.7}
\global\def\_@eqn@CurvatureEquationNegativeCase{\relax \unhbox \voidb@x \hbox {{\relax \tenrm (3.12)}}}
\global\def\_@proc@CurvatureEquationNegativeCase{3.7.1}
\global\def\_@eqn@HyperbolicSineGordon{\relax \unhbox \voidb@x \hbox {{\relax \tenrm (3.13)}}}
\global\def\_@proc@HyperbolicSineGordon{3.7.2}
\global\def\_@eqn@HazzidakiMetric{\relax \unhbox \voidb@x \hbox {{\relax \tenrm (3.14)}}}
\global\def\_@eqn@HazzidakiInequality{\relax \unhbox \voidb@x \hbox {{\relax \tenrm (3.15)}}}
\global\def\_@proc@HazzidakiInequality{3.7.3}
\global\def\_@proc@HilbertsTheorem{3.7.4}
\global\def\_@proc@HilbertsTheoremB{3.7.5}
\global\def\_@eqn@NoSolutionOfSinGordon{\relax \unhbox \voidb@x \hbox {{\relax \tenrm (3.16)}}}
\global\def\_@proc@NoSolutionOfSinGordon{3.7.6}
\global\def\_@head@BilegendrianSurfacesOverTheSphere{4}
\global\def\_@subhead@OverviewBilegendrianSurfaces{4.1}
\global\def\_@eqn@IntroWaveEquation{\relax \unhbox \voidb@x \hbox {{\relax \tenrm (4.1)}}}
\global\def\_@subhead@Quaternions{4.2}
\global\def\_@eqn@DefinitionOfAQuatCaseI{\relax \unhbox \voidb@x \hbox {{\relax \tenrm (4.2)}}}
\global\def\_@eqn@DefinitionOfAQuatCaseII{\relax \unhbox \voidb@x \hbox {{\relax \tenrm (4.3)}}}
\global\def\_@proc@DefinitionOfAQuatCase{4.2.1}
\global\def\_@eqn@DefinitionOfIJKForFlatSurfaces{\relax \unhbox \voidb@x \hbox {{\relax \tenrm (4.4)}}}
\global\def\_@subhead@BilegendrianSurfaces{4.3}
\global\def\_@eqn@FlatMetricOfBilegImm{\relax \unhbox \voidb@x \hbox {{\relax \tenrm (4.5)}}}
\global\def\_@eqn@BLengthOfIthVector{\relax \unhbox \voidb@x \hbox {{\relax \tenrm (4.6)}}}
\global\def\_@eqn@COverDiagB{\relax \unhbox \voidb@x \hbox {{\relax \tenrm (4.7)}}}
\global\def\_@eqn@COverDiagC{\relax \unhbox \voidb@x \hbox {{\relax \tenrm (4.8)}}}
\global\def\_@eqn@NormalIsLeftAndRightTransported{\relax \unhbox \voidb@x \hbox {{\relax \tenrm (4.9)}}}
\global\def\_@proc@NormalIsLeftAndRightTransported{4.3.1}
\global\def\_@eqn@ProductCriterionIsSatisfied{\relax \unhbox \voidb@x \hbox {{\relax \tenrm (4.10)}}}
\global\def\_@proc@ProductCriterionIsSatisfied{4.3.2}
\global\def\_@subhead@Factorization{4.4}
\global\def\_@eqn@FactorizationOfBilegendrians{\relax \unhbox \voidb@x \hbox {{\relax \tenrm (4.11)}}}
\global\def\_@proc@FactorizationOfBilegendrians{4.4.1}
\global\def\_@eqn@GrassmannianOfOrientedPlanes{\relax \unhbox \voidb@x \hbox {{\relax \tenrm (4.12)}}}
\global\def\_@eqn@FactorisationCriterionI{\relax \unhbox \voidb@x \hbox {{\relax \tenrm (4.13)}}}
\global\def\_@eqn@FactorisationCriterionII{\relax \unhbox \voidb@x \hbox {{\relax \tenrm (4.14)}}}
\global\def\_@proc@FactorisationCriterion{4.4.2}
\global\def\_@eqn@DoublyPeriodicCaseI{\relax \unhbox \voidb@x \hbox {{\relax \tenrm (4.15)}}}
\global\def\_@eqn@DoublyPeriodicCaseII{\relax \unhbox \voidb@x \hbox {{\relax \tenrm (4.16)}}}
\global\def\_@proc@DoublyPeriodicCase{4.4.3}
\global\def\_@subhead@UnitTorsionCurvesAndTheHopfFibration{4.5}
\global\def\_@eqn@DefinitionOfAlphaAndW{\relax \unhbox \voidb@x \hbox {{\relax \tenrm (4.17)}}}
\global\def\_@proc@HopfProjectionOfQuasiperiodicCurve{4.5.1}
\global\def\_@proc@RiemannianSubmersion{4.5.2}
\global\def\_@eqn@DefinitionOfArea{\relax \unhbox \voidb@x \hbox {{\relax \tenrm (4.18)}}}
\global\def\_@proc@DefinitionOfArea{4.5.3}
\global\def\_@eqn@ContactPropertyOfAlpha{\relax \unhbox \voidb@x \hbox {{\relax \tenrm (4.19)}}}
\global\def\_@proc@ContactPropertyOfAlpha{4.5.4}
\global\def\_@proc@HolonomyCondition{4.5.5}
\global\def\_@subhead@Curvature{4.6}
\global\def\_@eqn@ProjectionOfFirstVector{\relax \unhbox \voidb@x \hbox {{\relax \tenrm (4.20)}}}
\global\def\_@eqn@ProjectionOfSecondVector{\relax \unhbox \voidb@x \hbox {{\relax \tenrm (4.21)}}}
\global\def\_@eqn@ExplicitFormulaForC{\relax \unhbox \voidb@x \hbox {{\relax \tenrm (4.22)}}}
\global\def\_@proc@ExplicitFormulaForC{4.6.1}
\global\def\_@eqn@COverDiagA{\relax \unhbox \voidb@x \hbox {{\relax \tenrm (4.23)}}}
\global\def\_@proc@COverDiagA{4.6.2}
\global\def\_@proc@FramedCurve{4.6.3}
\global\def\_@eqn@CurvatureOfAsymptoticCurves{\relax \unhbox \voidb@x \hbox {{\relax \tenrm (4.24)}}}
\global\def\_@proc@CurvatureOfAsymptoticCurves{4.6.4}
\global\def\_@eqn@TorsionOfAsymptoticCurves{\relax \unhbox \voidb@x \hbox {{\relax \tenrm (4.25)}}}
\global\def\_@proc@TorsionOfAsymptoticCurves{4.6.5}
\global\def\_@rmk@TorsionOfAsymptoticCurves{\relax \unhbox \voidb@x \hbox {4.6.1}}
\global\def\_@eqn@WaveOfThetaVanishes{\relax \unhbox \voidb@x \hbox {{\relax \tenrm (4.26)}}}
\global\def\_@eqn@WaveOfThetaVanishesII{\relax \unhbox \voidb@x \hbox {{\relax \tenrm (4.27)}}}
\global\def\_@proc@WaveOfThetaVanishes{4.6.6}
\global\def\_@rmk@WaveOfThetaVanishes{\relax \unhbox \voidb@x \hbox {4.6.2}}
\global\def\_@eqn@ClassificationOfCompleteImmersions{\relax \unhbox \voidb@x \hbox {{\relax \tenrm (4.28)}}}
\global\def\_@proc@ClassificationOfCompleteImmersions{4.6.7}
\global\def\_@proc@AreaIsInteger{4.6.8}
\global\def\_@head@CalculatingTheCurvature{A}
\global\def\_@eqn@CurvatureFormula{\relax \unhbox \voidb@x \hbox {{\relax \tenrm (A.1)}}}
\global\def\_@proc@CurvatureFormula{A.0.9}
\global\def\_@head@Bibliography{B}
 }%
\else{\openout\references=references.tex }%
\fi%
%
%
\newcount\headno%
\global\headno=0%
\def\headinfo{\ifinappendices\alphanum\headno\else\the\headno\fi}%
\def\nextheadno{\global\advance\headno by 1 \global\subheadno=0 \global\eqnno=0 \headinfo}%
\makecounter{head}%
%
%
\newcount\subheadno%
\global\subheadno=0%
\def\subheadinfo{\headinfo.\the\subheadno}%
\def\nextsubheadno{\global\advance\subheadno by 1 \global\procno=0 \global\rmkno=0 \subheadinfo}%
\makecounter{subhead}%
%
%
\newcount\procno%
\global\procno=0%
\def\procinfo{\subheadinfo.\the\procno}%
\def\nextprocno{\global\advance\procno by 1 \procinfo}%
\makecounter{proc}%
%
%
\newcount\figno%
\global\figno=0%
\def\figinfo{\subheadinfo.\the\figno}%
\def\nextfigno{\global\advance\figno by 1 \figinfo}%
\makecounter{fig}%
%
%
\newcount\eqnno%
\global\eqnno=0%
\def\eqninfo{\text{{\rm (\headinfo.\the\eqnno)}}}%
\def\nexteqnno[#1]{\global\advance\eqnno by 1 \eqninfo\hbox{\eqnlabel{#1}}}%
\makecounter{eqn}%
%
%
\newcount\rmkno%
\global\rmkno=0%
\def\rmkinfo{\text{\subheadinfo.\the\rmkno}}%
\def\nextrmkno[#1]{\global\advance\rmkno by 1 \rmkinfo\hbox{\rmklabel{#1}}}%
\makecounter{rmk}%
%
%
%
%
%
\def\gobbleeight#1#2#3#4#5#6#7#8{}%
\newcount\citationno%
\global\citationno=0%
\def\citationinfo{\the\citationno}%
\makecounter{citation}%
\newwrite\biblio%
\def\newref#1#2{%
\def\temptext{#2}%
\edef\bibliotextoutput{\expandafter\gobbleeight\meaning\temptext}%
\global\advance\citationno by 1\citationlabel{#1}%
\ifmakebiblio%
    \edef\fileoutput{\write\biblio{\noindent\hbox to 0pt{\hss$[\the\citationno]$}\hskip 0.2em\bibliotextoutput\medskip}}%
    \fileoutput%
\fi}%
\def\cite#1{%
$[\citationref{#1}]$%
\ifmakebiblio%
    \edef\fileoutput{\write\biblio{#1}}%
    \fileoutput%
\fi%
}%
%
%
%
%
\let\mypar=\par%
\edef\Pagetitle={Blank}\headline={\hfil\Pagetitle\hfil}%
\edef\Pagefooter={Blank}\footline={\hfil\Pagefooter\hfil}%
%
%
\newcount\showpagenumflag%
\global\showpagenumflag=0 %
\def\nextoddpage%
{\newpage\ifodd\pageno%
\else\global\showpagenumflag=0 %
\null\vfil\eject%
\global\showpagenumflag=1 %
\fi}%
%
%
\font\headfont=cmb12%
\def\newhead#1[#2]%
{\ifhmode\mypar\fi%
\ifnum\headno=0 \else\goodbreak\bigskip\fi%
{\headfont\noindent\nextheadno\ - #1.}\headlabel{#2}%
\nobreak\medskip}%
%
%
\def\newsubhead#1[#2]%
{\ifhmode\mypar\fi%
\ifnum\subheadno=0 \else\goodbreak\medskip\fi%
{\bf\noindent\nextsubheadno\ - #1.\ }\subheadlabel{#2}}%
\def\newsubheadspecial#1[#2]%
{\ifhmode\mypar\fi%
\ifnum\subheadno=0 \else\goodbreak\medskip\fi%
{\bf\noindent\nextsubheadno\ - #1\ }\subheadlabel{#2}}%
%
%
\newif\ifinproclaim%
\global\inproclaimfalse%
\def\proclaim#1{%
\goodbreak\medskip
\bgroup\inproclaimtrue%
\noindent{\bf #1}%
\nobreak\medskip\sl}%
\def\noskipproclaim#1{%
\goodbreak\medskip%
\bgroup\inproclaimtrue%
\noindent{\bf #1}\nobreak\sl}%
\def\endproclaim{\mypar\egroup\nobreak\medskip\ignorespaces}%
%
%
%
\newcount\xpos\newcount\ypos
\def\makelabelgrid{%
\xpos=-5 \ypos=-5 %
\loop\ifnum\xpos<6 %
{\loop\ifnum\ypos<6 %
\def\labeltext{x}%
\ifnum\xpos=0\def\labeltext{+}\fi%
\ifnum\ypos=0\def\labeltext{+}\fi%
\placelabel[\xpos][\ypos]{\labeltext}%
\advance\ypos by 1 %
\repeat}%
\advance\xpos by 1 %
\repeat}%
\def\placelabel[#1][#2]#3{{%
\setbox10=\hbox{\raise #2cm \hbox{\hskip #1cm #3}}%
\ht10=0pt \dp10=0pt \wd10=0pt \box10}}%
%
%
%
%
\def\myitem#1{\noindent\hbox to .5cm{\hfill#1\hss}}%
%
%
%
%
%
%
%
%
%
\font\sansseriften=cmss10%
\font\sansserifseven=cmss7%
\font\sansseriffive=cmss5%
\newfam\sansseriffam%
\textfont\sansseriffam=\sansseriften%
\scriptfont\sansseriffam=\sansserifseven%
\scriptscriptfont\sansseriffam=\sansseriffive%
\def\mathsf{\fam\sansseriffam}%
%
%
%
\font\boldten=cmb10%
\font\boldseven=cmb7%
\font\boldfive=cmb5%
\newfam\mathboldfam%
\textfont\mathboldfam=\boldten%
\scriptfont\mathboldfam=\boldseven%
\scriptscriptfont\mathboldfam=\boldfive%
\def\mathbf{\fam\mathboldfam}%
%
%
%
\font\mycmmiten=cmmi10%
\font\mycmmiseven=cmmi7%
\font\mycmmifive=cmmi5%
\newfam\mycmmifam%
\textfont\mycmmifam=\mycmmiten%
\scriptfont\mycmmifam=\mycmmiseven%
\scriptscriptfont\mycmmifam=\mycmmifive%
\def\hexa#1{\ifcase #1 0\or 1\or 2\or 3\or 4\or 5\or 6\or 7\or 8\or 9\or A\or B\or C\or D\or E\or F\fi}%
\mathchardef\mathi="7\hexa\mycmmifam7B%
\mathchardef\mathj="7\hexa\mycmmifam7C%
%
%
\font\mymsbmten=msbm10 at 8pt%
\font\mymsbmseven=msbm7 at 5.6pt
\font\mymsbmfive=msbm5 at 4pt%
\newfam\mymsbmfam%
\textfont\mymsbmfam=\mymsbmten%
\scriptfont\mymsbmfam=\mymsbmseven%
\scriptscriptfont\mymsbmfam=\mymsbmfive%
\mathchardef\mybeth="7\hexa\mymsbmfam69%
\mathchardef\mygimmel="7\hexa\mymsbmfam6A%
\mathchardef\mydaleth="7\hexa\mymsbmfam6B%
%
%
%
%
\def\proof{{\noindent\bf Proof:\ }}%
\def\remark[#1]{{\noindent\bf Remark \nextrmkno[#1].}}%
\def\qed{~$\square$}%
\def\makeop#1{\global\expandafter\def\csname op#1\endcsname{{\text{#1}}}}%
\def\makeopsmall#1{\global\expandafter\def\csname op#1\endcsname{{\text{\lowercase{#1}}}}}%
%
%
%
\def\minter{\mathop{\cap}}%
%
%
\makeop{Ext}%
\makeop{Int}%
\makeop{Dist}%
\makeop{Diam}%
\makeop{Length}%
%
%
%
%
%
\def\mlim{\mathop{{\text{Lim}}}}%
\def\msup{\mathop{{\text{Sup}}}}%
%
%
%
\makeop{Dim}%
\makeop{Ker}%
\makeop{Coker}%
\makeop{Tr}%
\makeop{Adj}%
\makeop{Det}%
\makeop{End}%
\makeop{Lin}%
\makeop{Symm}%
\makeop{Mult}%
%
%
\makeop{dx}%
\makeop{dy}%
\makeop{dz}%
\makeop{dt}%
\makeop{dVol}%
\makeop{dArea}%
\makeop{Supp}%
\makeop{Hess}%
\makeop{Lip}%
%
%
\makeop{Re}%
\makeop{Im}%
\makeop{Arg}%
\makeop{Log}%
\makeop{Exp}%
%
%
\makeopsmall{Cos}%
\makeopsmall{Sin}%
\makeopsmall{Tan}%
\makeopsmall{Sec}%
\makeopsmall{Cosec}%
\makeopsmall{Cot}%
\makeopsmall{ArcCos}%
\makeopsmall{ArcSin}%
\makeopsmall{ArcTan}%
\makeopsmall{ArcSec}%
\makeopsmall{ArcCosec}%
\makeopsmall{ArcCot}%
%
%
\makeopsmall{Cosh}%
\makeopsmall{Sinh}%
\makeopsmall{Tanh}%
\makeopsmall{ArcCosh}%
\makeopsmall{ArcSinh}%
\makeopsmall{ArcTanh}%
%
%
\makeop{Vol}%
\makeop{Area}%
\makeop{Riem}%
\makeop{Ric}%
\makeop{Scal}%
\makeop{Euc}%
\makeop{Imm}%
\makeop{Emb}%
%
%
\makeop{Id}%
\makeop{Ad}%
\makeop{O}%
\makeop{SO}%
\makeop{SL}%
\makeop{GL}%
\makeop{Conf}%
\makeop{Homeo}%
\makeop{Diff}%
\makeop{Isom}%
%
%
\makeop{Ind}%
\makeop{Sig}%
\makeop{Spec}%
%
%
\makeop{Conv}%
\makeop{Max}%
\makeop{Min}%
\makeop{Mod}%
\makeop{Deg}%
\makeop{loc}%
%
%
%
%
%
%
%
%
%
%
%
%
%
 %
\input graphicx.tex %
%
%
%
%
%
\def\Pagetitle{\hfil\ifnum\pageno=1\null\else{\rm Clifford structures, bilegendrian surfaces, and extrinsic curvature}\fi\hfil}
\def\Pagefooter{\hfil{\myfontdefault\folio}\hfil}
\catcode`\@=11
\def\triplealign#1{\null\,\vcenter{\openup1\jot \m@th %
\ialign{\strut\hfil$\displaystyle{##}\quad$&$\displaystyle{{}##}$\hfil&$\displaystyle{{}##}$\hfil\crcr#1\crcr}}\,}
\def\multiline#1{\null\,\vcenter{\openup1\jot \m@th %
\ialign{\strut$\displaystyle{##}$\hfil&$\displaystyle{{}##}$\hfil\crcr#1\crcr}}\,}
\catcode`\@=12
\makeop{I}%
\makeop{II}%
\makeop{III}%
\makeopsmall{Coth}%
\makeopsmall{Cossec}%
\makeop{dA}%
\makeop{PSL}%
\makeop{Hyp}%
\makeop{Ln}
\makeop{A}
\makeop{J}
\makeop{U}%
\makeop{Umb}%
\makeop{T}%
\makeop{AdS}%
\makeop{sgn}%
\makeop{Spin}%
\makeop{Supp}%
\def\qi{{\mathbf{i}}}%
\def\qj{{\mathbf{j}}}%
\def\qk{{\mathbf{k}}}%

\def\opso{{\frak{so}}}
\makeop{Inv}
\makeop{Fix}
\makeop{Stab}
\makeop{Gr}
\makeop{Exc}
\makeop{Len}
\makeop{Vert}
\makeop{Cl}
\makeop{O}

\def\emph{\ifinproclaim\rm\else\sl\fi}
\def\bil{{\mathbf{b}}}
\makeop{ad}
\def\opUS{{\text{\rm U}}\Bbb{S}^3}

\makeop{Len}
\newref{BarBegZegh}{Barbot T., B\'eguin F., Zeghib A., Prescribing Gauss curvature of surfaces in $3$-dimensional spacetimes: application to the Minkowski problem in the Minkowski space, {\sl Ann. Inst. Fourier}, {\bf 61}, no. 2, (2011), 511–-591}
\newref{Bianchi}{Bianchi L., Sulle superficie a curvatura nulla in geometria ellittica, {\sl Ann. Mat. Pura Appl.}, {\bf 24}, 93--129, (1896)}
\newref{BonMonSchI}{Bonsante F., Mondello G., Schlenker J. M., A cyclic extension of the earthquake flow I, {\sl Geom. Topol.}, {\bf 17}, no. 1, (2013), 157--234}
\newref{BonMonSchII}{Bonsante F., Mondello G., Schlenker J. M., A cyclic extension of the earthquake flow II, {\sl Ann. Sci. Ec. Norm. Sup\'er.}, {\bf 48}, no. 4, (2015), 811–-859}
\newref{Clifford}{Clifford W. K., Preliminary sketch of biquaternons, {\sl Proc. London Math. Soc.}, {\bf IV}, 381--395, (1873)}
\newref{CliffordI}{Clifford W. K., Applications of Grassmann's extensive algebra, {\sl Amer. J. Math.}, {\bf 1}, 350--358, (1878)}
\newref{CliffordII}{Clifford W. K., On the classification of geometric algebras, in {\sl Mathematical papers of William Kingon Clifford}, Tucker R. (ed), MacMillan, London, (1882)}
\newref{CruForGad}{Cruceanu V., Fortuny P., Gadea P. M., A Survey on Paracomplex Geometry, {\sl Rocky Mt. J. Math.}, {\bf 26}, no. 1, 83--115, (1996)}
\newref{DajNom}{Dajczer M., Nomizu K., On flat surfaces in $\Bbb{S}^3$ and $\Bbb{H}^3$, in {\sl Manifolds and Lie Groups}, Hano J., Morimoto A., Murakami S., Okamoto K., Ozeki H. (eds.), Progress in Mathematics, {\bf 14}, Birkh\"auser-Verlag, (1981)}
\newref{EnomotoI}{Enomoto K., The Gauss image of flat surfaces in $\Bbb{R}^4$, {\sl Kodai Math. J.}, {\bf 9}, 19--32, (1986)}
\newref{EnomotoII}{Enomoto K., Global properties of the Gauss image of flat surfaces in $\Bbb{R}^4$, {\sl Kodai Math. J.}, {\bf 10}, 272--284, (1987)}
\newref{EnoKitWei}{Enomoto K., Kitigawa Y., Weiner J. L., A rigidity theorem for the Clifford tori in $\Bbb{S}^3$, {\sl Proc. Amer. Math. Soc.}, {\bf 124}, 265--268, (1996)}
\newref{FillSmi}{Fillastre F., Smith G., Group actions and scattering problems in Teichm\"uller theory, in {\sl Handbook of group actions IV}, Advanced Lectures in Mathematics, {\bf 40}, (2018), 359--417}
\newref{GalMarMil}{G\'alvez J. A., M\'artinez A., Mil\'an F., Flat surfaces in the hyperbolic $3$-space}
\newref{GalMirI}{G\'alvez J. A., Mira P., Embedded isolated singularities of flat surfaces in hyperbolic $3$-space, {\sl Cal. Var. PDE.}, {\bf 24}, 239--260, (2005)}
\newref{GalMirII}{G\'alvez J. A., Mira P., Isometric immersions of $\Bbb{R}^2$ into $\Bbb{R}^4$ and perturbation of Hopf tori, {\sl Math. Zeit.}, {bf 266}, 207--227, (2010)}
\newref{Garling}{Garling D. J. H., {\sl Clifford algebras: an introduction}, Cambridge University Press, (2011)}
\newref{Goldstein}{Goldstein H. P., {\sl Classical mechanics}, Pearson, (2001)}
\newref{HarveyLawson}{Harvey F. R., Lawson H. B., Pseudoconvexity for the special lagrangian potential equation, {\sl Calc.Var. PDEs.}, {\bf 60}, (2021)}
\newref{HitchinII}{Hitchin N. J., The moduli space of special lagrangian submanifolds, {\sl, Annali della Scuola Normale Superiore de Pisa, Classe de Scienze 4${}^{\text{e}}$ s\'erie}, {\bf 25}, no. 3--4, 503--515, (1997)}
\newref{HitchinI}{Hitchin N. J., The moduli space of complex lagrangian submanifolds, {\sl Asian J. Math.}, {\bf 3}, no. 1, 77--92, (1999)}
\newref{KitigawaII}{Kitigawa Y., Periodicity of the asymptotic curves in flat tori in $\Bbb{S}^3$, {\sl J. Math. Soc. Japan}, {\bf 40}, 457--476, (1988)}
\newref{KitigawaI}{Kitagawa Y., Isometric deformations of flat tori in the $3$-sphere with nonconstant mean curvature, {\sl Tohoku Math. J.}, {\bf 52}, no. 2, 283--298, (2000)}
\newref{KokUmeYam}{Kokubu M., Umehara M., Yamada K., Flat fronts in hyperbolic $3$-space, {\sl Pac. J. Math.}, {\bf 216}, 149--175, (2004)}
\newref{KokRossSajUmeYam}{Kokubu M., Rossman W., Saji K., Umehara M., Yamada K., Singularities of flat fronts in hyperbolic $3$-space, {\sl Pac. J. Math.}, {\bf 221}, 303--351, (2005)}
\newref{LabIII}{Labourie F., M\'etriques prescrites sur le bord des vari\'et\'es hyperboliques de dimension $3$, {\sl J. Diff. Geom.}, {\bf 35}, no. 3, (1992), 609--626}
\newref{LabI}{Labourie F., Exemples de courbes pseudo-holomorphes en g\'eom\'etrie riemannienne, in {\sl Holomorphic curves in symplectic geometry}, Progress in Mathematics, {\bf 117}, Birkh\"auser, (1994)}
\newref{LabII}{Labourie F., Probl\`emes de Monge-Amp\`ere, courbes pseudo-holomorphes et laminations, {\sl Geom. Funct. Anal.}, {\bf 7}, (1997), 496--534}
\newref{LabIV}{Labourie F., Un Lemme de Morse pour les surfaces convexes, {\sl Inventiones Mathematicae}, {\bf 141}, No 2, (2000), 239--297}
\newref{Mira}{Le\'on-Guzm\'an M. A., Mira P., Pastor J. A., The space of Lorentzian flat tori in anti de-Sitter $3$-space, {\sl Trans. AMS.}, {\bf 363}, no. 12, 6549--6573, (2011)}
\newref{Lounesto}{Lounesto P., {\sl Clifford algebras and spinors}, London Mathematical Society Lecture Note Series, {\bf 286}, (2001)}
\newref{LychRub}{Lychagin V. V., Rubtsov V. N., Local classification of Monge-Amp\`ere differential equations (Russian), {\sl Dokl. Akad. Nauk SSSR}, {\bf 272}, no. 1, (1983), 34--38}
\newref{LychRubChekI}{Lychagin V. V., Rubtsov V. N., Chekalov I. V., A classification of Monge-Amp\`ere equations, {\sl Ann. Sci. ENS.}, {\bf 26}, no. 3, 281--308, (1993)}
\newref{McDuffSalamon}{McDuff D., Salamon D., {\sl J-holomorphic curves and symplectic topology}, Colloquium publications, {\bf 52}, Amer. Math. Soc., (2012)}
\newref{Milnor}{Milnor T. K., Harmonic maps and classical surface theory in minkowski $3$-space, {\sl Trans. AMS.}, {\bf 280}, no. 1, 161--185, (1983)}
\newref{Pinkall}{Pinkall U., Hopf tori in $\Bbb{S}^3$, {\sl Inv. Math.}, {\bf 81}, 379--386, (1985)}
\newref{Rubtsov}{Rubtsov V. N., Geometry of Monge-Amp\`ere structures, in {\sl  Nonlinear PDEs, Their Geometry, and Applications}, Kycia R. A. et al. (eds.), Springer-Verlag, (2019)}
\newref{SajUmeYam}{Saji K., Umehara M., Yamada K., The geometry of fronts, {\sl Ann. Math.}, {\bf 169}, 491--529, (2009)}
\newref{Sasaki}{Sasaki S., On complete surfaces with Gaussian curvature zero in the $3$-sphere, {\sl Colloq. Math.}, {\bf 26}, 165--174, (1972)}
\newref{SasakiNomizu}{Sasaki S., Nomizu K., {\sl Affine differential geometry: geometry of affine immersions}, Cambridge Tracts in Mathematics, {\bf 111}, Cambridge University Press, (1994)}
\newref{Sch}{Schlenker J. M., Hyperbolic manifolds with convex boundary, {\sl Inv. Math.}, {\bf 163}, (2006), 109--169}
\newref{SmiSLC}{Smith G., Special Lagrangian curvature, Math. Annalen, 355, no. 1, (2013), 57-95}
\newref{SmiEW}{Smith G., On an Enneper-Weierstrass-type representation of constant Gaussian curvature surfaces in 3-dimensional hyperbolic space, arXiv:1404.5006}
\newref{SmiMS}{Smith G., M\"obius structures, hyperbolic ends and $k$-surfaces in hyperbolic space, in {\sl In the Tradition of Thurston, Vol. II}, Ohshika K., Papadopoulos A. (eds.), Springer Verlag, (2022)}
\newref{SmiAS}{Smith G., On the asymptotic geometry of finite-type $k$-surfaces in three-dimensional hyperbolic space, to appear in {\sl J. Eur. Math. Soc.}}
\newref{SmiQMAK}{Smith G., Quaternions, Monge-Amp\`ere structures, and $k$-surfaces, arXiv:2210.02664}
\newref{SmiQCKS}{Smith G., On quasicomplete $k$-surfaces in $3$-dimensional space-forms, arXiv:2211.14868}
\newref{Spivak}{Spivak M., {\sl A Comprehensive Introduction to Differential Geometry, Vol. IV}, Publish or Perish, Berkeley, (1977)}
\newref{Tamburelli}{Tamburelli A., Prescribing metrics on the boundary of anti-de Sitter $3$-manifolds, {\sl Int. Math. Res. Not.}, no. 5, (2018), 1281--1313}
\newref{Volkert}{Volkert K., Space forms: a history, {\sl Bulletin of the Manifold Atlas}, (2013)}
\newref{Weiner}{Weiner J. L., Flat tori in $\Bbb{S}^3$ and their Gauss maps, {\sl Proc. London Math. Soc.}, {\bf 62}, 54--76, (1991)}
\newref{Yau}{Yau S. T. Yau, Submanifolds with constant mean curvature II, {\sl Amer. J. Math.}, {\bf 97}, 76--100, (1975)}
\def\defaultalignment{\rightskip=0pt \leftskip=0pt \spaceskip=0pt \xspaceskip=0pt \parfillskip=0pt plus 1fil \parindent=20pt}%
\def\centre{\rightskip=0pt plus 1fil \leftskip=0pt plus 1fil \spaceskip=.3333em \xspaceskip=.5em \parfillskip=0em \parindent=0em}%
\def\textmonth#1{\ifcase#1\or January\or Febuary\or March\or April\or May\or June\or July\or August\or September\or October\or November\or December\fi}
\font\abstracttitlefont=cmr10 at 14pt {\abstracttitlefont\centre Clifford structures, bilegendrian surfaces, and extrinsic curvature\par}
\bigskip
{\centre 11th August 2023\par}
\bigskip
{\centre Graham Smith{\defaultalignment\numberedfootnote{Departamento de Matem\'atica, Pontif\'\i cia Universidade Cat\'olica do Rio de Janeiro (PUC-Rio), Rio de Janeiro, Brazil}}\par}
\bigskip
{\centre{\sl Dedicated to Fran\c{c}ois Labourie on the occasion of his 60th birthday}\par}
\bigskip
\noindent{\bf Abstract:~}We use Clifford algebras to construct a unified formalism for studying constant extrinsic curvature immersed surfaces in riemannian and semi-riemannian $3$-manifolds in terms of immersed bilegendrian surfaces in their unitary bundles. As an application, we provide full classifications of both complete and compact immersed bilegendrian surfaces in the unit tangent bundle $\opUS$ of the $3$-sphere.
\bigskip
\noindent{\bf AMS Classification:~}53A05, 12E15, 35J96
\bigskip
\noindent{\bf Key words~:~}Clifford algebras, extrinsic curvature, holomorphic curves, bilegendrian curves
\bigskip
%
%
\newhead{Introduction}[Introduction]
\newsubhead{Introduction}[Introduction]
The symplectic nature of riemannian geometry is well-known but rarely emphasized. One theory, however, in which the symplectic perspective has proven itself to be particularly fruitful is that of immersed surfaces of constant extrinsic curvature in riemannian $3$-manifolds. In this paper, we describe a framework for studying constant extrinsic curvature surfaces in terms of bilegendrian surfaces in the unit tangent bundle. We study the interplay between these two classes of surfaces and a third, which we will call $J$-{\emph invariant} surfaces. This interplay, mediated by Clifford structures over the contact distribution of the total space of the unit tangent bundle, provides a unified perspective on many known results together with new insights into others. This paper is intended as a survey, although some new results will also be proven.
\newsubhead{Discussion}[CECSurfaces]
We define an {\emph immersed surface} in a $3$-manifold $X$ to be a pair $S:=(S,e)$, where $S$ is a surface and $e:S\rightarrow X$ is a smooth immersion. When $X$ is riemannian, we define the {\emph extrinsic curvature} of $S$ to be equal to the determinant of its shape operator. For $k\in\Bbb{R}\setminus\{0\}$, we define a {\emph CEC-$k$ surface} in $X$ to be an immersed surface in $X$ of constant extrinsic curvature equal to $k$.\numberedfootnote{CEC-$0$ surfaces, which also exhibit fascinating geometric phenomena, will not be addressed in this paper.} When the curvature is not specified, we will simply speak of {\emph CEC surfaces}.
\par
CEC surfaces can be viewed as solutions of Monge-Amp\`ere type partial differential equations, which may be elliptic or hyperbolic, according to whether the curvature is positive or negative. The use of symplectic structures over the total space of the tangent bundle $TX$ to study such partial differential equations was pioneered by Lychagin--Rubtsov in \cite{LychRub} (see also \cite{LychRubChekI} and \cite{Rubtsov}). However, it was in a series of papers - \cite{LabIII}, \cite{LabI}, \cite{LabII}, \cite{LabIV} - that Labourie used {\emph contact} structures over the total space of the {\emph unit} tangent bundle $UX$ of $X$ to apply these ideas to the study of CEC surfaces. More precisely, given an immersed surface $(S,e)$, we denote by $\hat{e}:=\nu_e$ its unit normal vector field, and we call the surface $\hat{S}:=(S,\hat{e})$ its {\emph Gauss lift}. Labourie showed that, for all positive $k$, there exists an almost complex structure $J_k$ over the contact bundle of $UX$ such that $(S,e)$ is a CEC-$k$ surface if and only if its Gauss lift is $J_k$-pseudoholomorphic. He then used Gromov's theory of pseudoholomorphic curves to prove a simple, yet powerful, compactness result for families of such surfaces. Labourie's work has since revolutionized our understanding of positively-curved CEC surfaces, and has found applications across a variety of fields.\numberedfootnote{For applications to hyperbolic geometry and dynamical systems theory, we refer the reader to Labourie's own work \cite{LabIII} and \cite{LabIV}. For applications to general relativity and Teichm\"uller theory, we refer the reader to \cite{BarBegZegh}, \cite{BonMonSchI}, \cite{BonMonSchII}, \cite{Tamburelli}, and our survey \cite{FillSmi}. Recently, in \cite{SmiMS}, \cite{SmiQMAK} and \cite{SmiQCKS}, we showed how Labourie's ideas yield a complete classification of positively-curved CEC surfaces in $3$-dimensional space-forms, subject to a natural completeness condition.}
\par
In \cite{SmiQMAK}, we showed how Labourie's almost complex structures can be seen as components of quaternionic structures over the contact bundle of $UX$. We had already adopted this point of view implicitely in \cite{SmiSLC}, where we used almost Calabi-Yau structures to generalize Labourie's ideas to hypersurfaces in higher-dimensional manifolds (the reader may also consult the work \cite{HarveyLawson} of Harvey--Lawson for an alternative perspective on these ideas). In the present paper, by replacing quaternionic structures with Clifford structures, we extend Labourie's framework to include, not only negatively-curved CEC surfaces, but also mixed-signature CEC surfaces in semi-riemannian $3$-manifolds.
\par
We will study CEC surfaces in terms of what we will call {\emph bilegendrian surfaces}, which build on the concept of bilagrangian surfaces introduced by Hitchin in \cite{HitchinII} and \cite{HitchinI}. We will also divide bilegendrian surfaces into two classes, namely those of {\emph elliptic type}, which generate elliptic partial differential equations, and those of {\emph hyperbolic type}, which generate hyperbolic ones. Hitchin's bilagrangian surfaces, for example, are of elliptic type, as are Labourie's $k$-surfaces. We will encounter noteworthy bilegendrian surfaces of hyperbolic type presently.
\par
It is interesting to see how our framework of bilegendrian surfaces unifies known results for CEC surfaces. For example, when $X$ is a space-form, we construct in Theorem \procref{FlatMetric} a canonical flat metric, which may be riemannian or semi-riemannian, over the complement of a certain exceptional set of every bilegendrian surface in $UX$. When the curvature is positive, this identifies with the canonical flat structure obtained upon integrating the square root of the Hopf differential, whilst, when the curvature is negative, it identifies with the asymptotic Chebyshev net, used in Hilbert's proof of the non-existence of complete negatively-curved CEC surfaces in Euclidean $3$-space $\Bbb{R}^3$ (see Sections \subheadref{SurfacesOfConstantPositiveCurvature}, \subheadref{SurfacesOfConstantNegativeCurvature}, and \subheadref{HilbertsTheorem}).
\par
The theory of flat fronts is also encompassed within our framework. This theory, developed by Kokubu--Umehara--Yamada in \cite{KokUmeYam}, building on ideas elegantly presented by G\'alvez--Mart\'\i nez--Mil\'an in \cite{GalMarMil}, serves to analyse the structure of singularities of CEC-$1$ surfaces in hyperbolic $3$-space $\Bbb{H}^3$, and has since been studied by various authors in a variety of settings (see, for example, \cite{GalMirI}, \cite{KokRossSajUmeYam}, and \cite{SajUmeYam}). More recently, in \cite{GalMirII}, G\'alvez--Mira introduced a second class of flat fronts in order to study singularities of CEC-$(-1)$ surfaces in the $3$-sphere $\Bbb{S}^3$. In fact, flat fronts of the first type form a special case of the surfaces studied by Labourie (more precisely, they are $J_1$-holomorphic curves). Flat fronts of the second type are not covered by Labourie's theory, but instead constitute an interesting example of bilegendrian surfaces of hyperbolic type.
\par
We conclude this paper by classifying complete and compact bilegendrian surfaces in $\opUS$, which correspond to the flat fronts studied in \cite{GalMirII}. This generalizes the classifications of compact flat surfaces in $\Bbb{S}^3$ given by Kitigawa in \cite{KitigawaII} and by Weiner in \cite{Weiner}. We underline, however, that this generalization is not automatic, as the techniques of Kitigawa and Weiner, which rely on the geometry of the projection onto the base space, break down at points where this projection degenerates. Instead, we will show how the bilegendrian condition naturally interacts with the Lie group structure of $\Bbb{S}^3$ to yield canonical factorizations of the immersions and of their period lattices in the compact case. Furthermore, we will see that these factorization results readily follow from a general factorization criterion for Lie group-valued functions over $\Bbb{R}^2$, making them, in a sense, more algebraic in nature than geometric.
\par
Finally, our bilegendrian framework also has interesting applications to the semi-riemannian case. For example, CEC-$1$ surfaces of mixed signature in the anti de-Sitter space $\opAdS^{2,1}$, studied by Dajczer--Nomizu in \cite{DajNom} and Le\'on-Guzm\'an--Mira--Pastor in \cite{Mira}, correspond to bilegendrian surfaces of hyperbolic type. Likewise, negatively-curved CEC surfaces of mixed signature in spacetimes yield bilegendrian surfaces of elliptic type, and should therefore be amenable to study via Labourie's techniques, with potentially interesting applications. In fact, the elliptic nature of such surfaces was already observed by T. K. Milnor in \cite{Milnor}, but we are not aware of this property having been further examined by later authors.
\newsubheadspecial{Happy 60th Birthday Fran\c{c}ois!}[AFinalWord]
It is now more than 20 years since I first attended Fran\c{c}ois' course on Teichm\"uller theory in the B\^atiment 425 of what was then known as l'Universit\'e Paris-Sud (Paris-XI, now Paris-Saclay). His intuitive approach to geometry, which I later learned had been influenced by Gromov, was for me both novel and inspiring. Over the three years I spent as his PhD student, I learned many important things, including, amongst others, how to correctly serve a glass of wine\numberedfootnote{Never more than half-full.}. Fran\c{c}ois knew how to provide the right encouragement and the right advice at those key moments during my doctoral research which to me seemed the darkest. Above all, however, from that time on, and throughout the rest of my career, the distance of his mathematical vision, and the profundity of his insights, have been an inspiration, impelling me to strive, within the limits of my ability, to achieve similar excellence in my own work. The subject of this paper, as I hope Fran\c{c}ois will recognize, is a tribute to that experience: a new perspective on the original problem that he proposed to me when we first met. With this paper, I wish Fran\c{c}ois a very happy birthday, and many more fruitful years of mathematical creation to come!
\newsubhead{Acknowledgements}[Acknowledgements]
I am grateful to Pablo Mira, Carlos Tomei and Pierre Pansu for helpful comments made to earlier drafts of this paper. I am also grateful to Richard Wentworth and Jean-Marc Schlenker for having invited me to contribute to the volume in which this paper should appear, and for assuring me that the proposed topic would genuinely be of interest. Finally, of course, I am grateful to Fran\c{c}ois Labourie for having introduced me, all those years ago, to CEC surfaces and the marvels and wonders they conceal.
\newhead{Clifford algebras}[CliffordStructures]
\newsubhead{Overview}[OverviewCliffordStructures]
Clifford algebras were introduced by William Clifford in \cite{CliffordI} and \cite{CliffordII} with the aim of unifying Hamilton's theory of quaternions with that of Grassman algebras. In the present paper, they will be of interest to us as a source of endomorphisms $J$ such that $J^2 = -\epsilon\opId$, for some $\epsilon\in\{\pm 1\}$. We will call any such $J$ a {\emph pseudo-involution}, we will call $\epsilon$ its {\emph sign}, and we will say, furthermore, that $J$ is {\emph balanced} whenever its trace vanishes. Pseudo-involutions of positive sign are just complex structures. Balanced pseudo-involutions of negative sign are paracomplex structures.
\par
Given a pseudo-involution $J$ of a $4$-dimensional real vector space $E$, we will be interested in the planes in $E$ that $J$ preserves. By mild abuse of terminology, we will say that a plane $P$ is $J$-{\emph invariant} whenever it is preserved by $J$, but is not a real $J$-eigenspace. The interest in working with Clifford structures is now precisely that they yield pseudo-involutions whose invariant planes are bilegendrian and vice-versa. More precisely, given two symplectic forms $\omega_i$ and $\omega_k$ over $E$, we say that a plane $P$ in $E$ is $(\omega_i,\omega_k)$-{\emph bilagrangian} whenever it is lagrangian with respect to both $\omega_i$ and $\omega_k$. In Theorem \procref{InvariantPlanesAreBiLagrangian}, we will show that, for all the pseudo-involutions $J$ of interest to us, there exists a pair $(\omega_i,\omega_k)$ of symplectic forms such that any plane $P$ in $E$ is $J$-invariant if and only if it is $(\omega_i,\omega_k)$-bilagrangian.
\par
We conclude this section by studying the geometries of $J$-invariant planes. In particular, we identify special vectors which exist for all but a certain set of exceptional planes. These vectors correspond to eigenvectors of symmetric matrices, and will be of use in studying bilegendrian surfaces in the next section.
\newsubhead{Pseudo-isometries and pseudo-involutions}[PseudoInvolutions]
Before introducing Clifford structures, it is useful to study the properties of the objects of which they are comprised, namely pseudo-involutions. Let $E$ be a $2d$-dimensional real vector space. We will say that an endomorphism $M$ of $E$ is a {\emph pseudo-involution} whenever
$$
M^2 = -\epsilon\opId\ ,
$$
for some $\epsilon\in\{\pm 1\}$, which we call its {\emph sign}. We will say that a pseudo-involution $M$ is {\sl balanced} whenever
$$
\opTr(M) = 0\ .
$$
Pseudo-involutions of positive sign, which are always balanced, are just complex structures. They are associated with elliptic partial differential equations, and thus have widespread and well-known geometric applications (see, for example, \cite{McDuffSalamon}). Balanced pseudo-involutions of negative sign are also known as para-complex structures. They are associated with hyperbolic partial differential equations, and have also attracted the interest of many authors (see, for example, \cite{Rubtsov}, and the survey \cite{CruForGad}).
\par
The following construction of pseudo-involutions will be particularly useful. Let $\bil$ be a non-degenerate, symmetric, bilinear form over $E$. We will say that $\bil$ has {\emph positive} sign whenever it has an even number of positive directions, and we will say that it has {\emph negative} sign otherwise. We say that an endomorphism $M$ of $E$ is a {\emph pseudo-isometry} of $\bil$ whenever
$$
\bil(M\cdot,M\cdot) = \pm \bil(\cdot,\cdot)\ .
$$
We denote the Lie group of pseudo-isometries of $\bil$ by $\hat{\opO}(\bil)$, and we denote its Lie algebra by $\opso(\bil)$. Note that every element $M$ of $\opso(\bil)\minter\hat{\opO}(\bil)$ satisfies
$$
M^2 = -M^*M = \pm\opId\ ,
$$
and is thus a pseudo-involution.
\proclaim{Lemma \& Definition \nextprocno}
\noindent If $M$ is a pseudo-involution of negative sign, then $M$ is diagonalizable with eigenvalues $\{\pm 1\}$. In particular, $M$ is balanced if and only if its two eigenspaces have the same dimension.
\endproclaim
\proclabel{PseudoInvolutionsHaveEigenspaces}
\proof Indeed, if $S+N:=M$ denotes the Jordan decomposition of $M$, then
$$
\opId = \pm M^2 = \pm S^2 \pm (2SN+N^2)\ .\eqnum{\nexteqnno[JordanDecompositionOfIdentity]}
$$
Trivially $S^2$ is diagonal, $(2SN+N^2)$ is nilpotent, and the two commute with one-another, so that \eqnref{JordanDecompositionOfIdentity} is the Jordan decomposition of the identity. Since the Jordan decomposition is unique, it follows that
$$
S^2 = \opId\ \text{and}\ (2SN+N^2) = 0\ .
$$
It follows from the first relation that the eigenvalues of $S$ are $\pm 1$. Since $SN$ is nilpotent, the second relation yields
$$
N = \frac{1}{2}S\cdot\bigg(2SN+N^2\bigg)\cdot\sum_{k=0}^\infty\frac{(-1)^k}{2^k}(SN)^k = 0\ ,
$$
so that $M$ is diagonalizable, as desired.\qed
\medskip
For any pseudo-involution $M$ of $E$, let $\opInv(M)$ denote the space of $M$-invariant $d$-dimensional subspaces of $E$, furnished with the topology that it inherits as a subspace of the grassmannian. When $M$ has negative sign and is balanced, this set also contains its two real eigenspaces, which we will find preferable to exclude in the sequel. For this reason, we denote by $\opInv^*(M)$ the set of $M$-invariant $d$-dimensional subspaces of $E$ which are not $M$-eigenspaces. Note that, when $M$ has positive sign, it has no real eigenspaces, and $\opInv^*(M)$ coincides with $\opInv(M)$. By abuse of terminology, in later sections, we will say that a plane is $M$-{\emph invariant} whenever it is an element of $\opInv^*(M)$.
\par
In this paper, we will be concerned with the case where $E$ is $4$-dimensional.
\proclaim{Lemma \nextprocno}
\noindent Let $M$ be a pseudo-involution of some $4$-dimensional real vector space $E$.
\medskip
\myitem{(1)} If $M$ has positive sign, then $\opInv^*(M)$ is homeomorphic to $\Bbb{S}^2$, and
\medskip
\myitem{(2)} if $M$ has negative sign and is balanced, then $\opInv^*(M)$ is homeomorphic to $\Bbb{S}^1\times\Bbb{S}^1$.
\endproclaim
\proclabel{InvariantPlanes}
\remark[InvariantPlanes] In particular, by compactness of the torus, whenever $M$ has negative sign and is balanced, its two eigenspaces are isolated in $\opInv(M)$.
\medskip
\proof Suppose first that $M$ has positive sign. Then $M$ is a complex structure, and $\opInv^*(M)$ identifies with $1$-dimensional complex projective space, which is homeomorphic to $\Bbb{S}^2$, as desired. Suppose now that $M$ is negative and balanced, and let $E_-$ and $E_+$ denote respectively its negative and positive eigenspaces. With respect to the decomposition $E=E_-\oplus E_+$, $M$ has the form
$$
M = \pmatrix -\opId\hfill& 0\hfill\cr 0\hfill& \opId\hfill\cr\endpmatrix\ .
$$
Now let $P$ be an $M$-invariant subspace of $E$ which is not an eigenspace. In particular, $P$ contains a vector $x:=(u,v)$ both of whose components are non-zero. Since $P$ is $M$-invariant, the vector $y:=(-u,v)$ is also an element of $P$, so that $P$ is generated by the pair $(u,0)$ and $(0,v)$. Conversely, any pair of non-vanishing vectors of this form trivially generates an $M$-invariant subspace which is not an eigenspace. It follows that $\opInv^*(M)$ is homeomorphic to the cartesian product of the projective space of $E_-$ with that of $E_+$. Since each of these projective spaces is homeomorphic to $\Bbb{S}^1$, the result follows.\qed
\newsubhead{Clifford algebras}[CliffordAlgebras]
We now present the main properties of Clifford algebras that will be required in the sequel. For a comprehensive introduction to this fascinating theory, we refer the reader to the excellent works \cite{Garling} of Garling and \cite{Lounesto} of Lounesto. Throughout what follows, we will only be concerned with the relatively trivial case of Clifford algebras defined over $2$-dimensional real vector spaces.
\par
Let $\bil$ be a symmetric, bilinear form defined over a $2$-dimensional real vector space $E$. We define the Clifford algebra $\opCl(\bil)$ to be the quotient of the tensor algebra $T^*E$ of $E$ by the two-sided ideal
$$
\Cal{I} := \langle u\otimes v + v\otimes u + 2\bil(u,v)\ |\ u,v\in E\rangle\ .
$$
When $\bil$ is trivial, $\opCl(\bil)$ is none other than the exterior algebra of $E$, also known as its Grassmann algebra, and when $\bil$ is positive-definite, it is straightforward to show that $\opCl(\bil)$ is isomorphic to the quaternion algebra.
\par
We will henceforth only be concerned with the case where $\bil$ is non-degenerate. Let $(e_1,e_2)$ be a $\bil$-orthonormal basis of $E$. The algebra $\opCl(\bil)$ is generated as a vector space by the identity together with the elements
$$
\qi := e_1\ ,\ \qj := e_2\ ,\ \text{and}\ \qk := e_1\cdot e_2\ .
$$
Note that each of these elements has unit (positive or negative) square, and that any two of them anticommute with one-another. The Clifford algebra carries a natural $\Bbb{Z}_2$-grading
$$
\opCl(\bil) = \opCl^+(\bil)\oplus\opCl^-(\bil)\ ,
$$
where $\opCl^+(\bil)$ and $\opCl^-(\bil)$ denote the vector subspaces generated respectively by products of even and odd numbers of elements of $E$. In the present case,
$$
\opCl^+(\bil) = \langle1,\qk\rangle\ \text{and}\ \opCl^-(\bil) = \langle\qi,\qj\rangle\ .
$$
We denote by $\pi^+:\opCl(\bil)\rightarrow\opCl^+(\bil)$ and $\pi^-:\opCl(\bil)\rightarrow\opCl^-(\bil)$ the canonical projections onto each of the factors. This decomposition will play a significant role in the sequel.
\par
The Clifford algebra carries one natural involution and two natural anti-involutions which we now describe. The involution, called the {\emph grade involution}, is defined by
$$
\hat{x} := \pi^+(x) - \pi^-(x)\ .
$$
That is, when
$$
x := a + b\qi + c\qj + d\qk\ ,\eqnum{\nexteqnno[FormOfX]}
$$
its image under the grade involution is
$$
\hat{x} = a - b\qi - c\qj + d\qk\ .
$$
The first anti-involution, called {\emph reversion}, is first defined on products of elements of $E$ by reversing the order of their factors, and is then extended to $\opCl(\bil)$ by linearity. With $x$ as in \eqnref{FormOfX}, its image under reversion is thus
$$
x^\dagger = a + b\qi + c\qj - d\qk\ .
$$
Finally, the second anti-involution, called {\emph Clifford conjugation}, is simply the composition of grade involution with reversion, so that, with $x$ again as in \eqnref{FormOfX}, its Clifford conjugate is
$$
\overline{x} = a - b\qi - c\qj - d\qk\ .
$$
These three operations trivially commute with one another.
\par
In the present, $2$-dimensional case, it makes sense to deepen the analogy with quaternions. We thus denote the $(+1)$-eigenspace of Clifford conjugation by $\Cal{R}$, and we call it the {\emph real subspace}. We likewise denote the $(-1)$-eigenspace of Clifford conjugation by $\Cal{I}$, and we call it the {\emph imaginary subspace}. Trivially,
$$
\Cal{R} = \langle 1\rangle\ \text{and}\ \Cal{I} = \langle \qi, \qj, \qk \rangle\ .
$$
In particular, $\Cal{R}$ is the centre of $\opCl(\bil)$, and $\Cal{I}$ is a Lie sub-algebra. We denote also by $\Cal{R}$ and $\Cal{I}$ the respective projections onto the real and imaginary subspaces. Note that $\Cal{R}$ can be viewed as a trace in the sense that, for all $x,y\in\opCl(\bil)$,
$$
\Cal{R}(x\cdot y) = \Cal{R}(y\cdot x)\ .
$$
\par
Since Clifford conjugation is an anti-involution, for all $x\in\opCl(\bil)$,
$$
\overline{x\cdot\overline{x}} = \overline{\overline{x}}\cdot\overline{x} = x\cdot\overline{x}\ ,
$$
so that $x\cdot\overline{x}$ is always real. We thus define the {\emph inner-product} $g$ over $\opCl(\bil)$ by
$$
g(x,y) := \Cal{R}(x\cdot\overline{y})\ .
$$
Note that the anticommutator of any two imaginary elements $x,y$ satisfies
$$
\{x,y\} = 2g(x,y)\ ,
$$
so that two imaginary elements anticommute if and only if they are orthogonal to one another.
\par
We denote the norm-squared of $g$ by $\|\cdot\|_g^2$. With $x$ as in \eqnref{FormOfX},
$$
\|x\|_g^2 = a^2 - b^2\qi^2 - c^2\qj^2 - d^2\qk^2\ ,
$$
so that the inner-product $g$ is non-degenerate with signature determined by the signs of the squares of $\qi$, $\qj$ and $\qk$. We will say that a unit-length imaginary element of $\opCl(\bil)$ is {\emph positive} or {\emph negative} according to the sign of its norm-squared. Note that this terminology is consistent with that introduced for pseudo-involutions in the preceding section.
\par
Given $x\in\opCl(\bil)$, we denote by $m_x$ the operation of multiplication on the left by this element. For all $x,y,z\in\opCl(\bil)$,
$$
g(m_x(y),m_x(z)) = \|x\|_g^2 g(y,z)\ ,
$$
so that $m_x$ is a pseudo-isometry whenever $x$ has unit length. Furthermore, when, in addition, $x$ is imaginary, for all $y,z\in\opCl(\bil)$,
$$
g(y,m_x(z)) = -g(m_x(y),z)\ ,
$$
so that $m_x$ is $g$-antisymmetric. In particular, for all unit-length, imaginary $x$, we define the non-degenerate alternating form $\omega_x$ by
$$
\omega_x(y,z) := g(y, m_x(z))\ .
$$
\par
Finally, we define the {\emph complementary inner-product} $\hat{g}$ over $\opCl(\bil)$ by
$$
\hat{g}(y,z) := g(\hat{y},z) = \Cal{R}(\hat{y}\cdot\overline{z})\ .
$$
Note that $\opCl^+(\bil)$ and $\opCl^-(\bil)$ are both $g$- and $\hat{g}$-orthogonal to one another. Furthermore
$$
\hat{g}|_{\opCl^+(\bil)} = g|_{\opCl^+(\bil)}\ \text{and}\
\hat{g}|_{\opCl^-(\bil)} = -g|_{\opCl^-(\bil)}\ .
$$
As before, for all $x\in\opCl^-(\bil)$, and for all $y,z\in\opCl(\bil)$,
$$
\hat{g}(m_x(y),m_x(z)) = -\|x\|_g^2\hat{g}(y,z)\ ,
$$
so that $m_x$ is a $\hat{g}$-pseudo-isometry whenever it has unit length. Likewise for all such $x$, $y$ and $z$,
$$
\hat{g}(y,m_x(z)) = g(m_x(y),z)\ ,
$$
so that $m_x$ is also $\hat{g}$-symmetric.
\par
We now characterize planes preserved by the action of multiplication on the left by unit-length imaginary elements of $\opCl(\bil)$.
\proclaim{Lemma \nextprocno}
\noindent If $x$ is a unit-length element of $\Cal{I}$ of negative sign, then $m_x$ is balanced. Furthermore, the eigenspaces of $m_x$ are given by
$$
E_+ = \langle 1+x,y + x\cdot y\rangle\ \text{and}\
E_- = \langle 1-x,y - x\cdot y\rangle\ ,\eqnum{\nexteqnno[EigenspacesOfMx]}
$$
where $y$ is any other unit-length element of $\Cal{I}$ orthogonal to $x$.
\endproclaim
\proclabel{NegativeUnitImaginaryElementsAreBalanced}
\proof Indeed, let $y$ be another unit-length element of $\Cal{I}$ orthogonal to $x$. Since $y$ anticommutes with $x$, $m_y$ interchanges the eigenspaces of $m_x$, so that $m_x$ is balanced, as desired. We verify \eqnref{EigenspacesOfMx} by inspection, and this completes the proof.\qed
\medskip
\noindent We now relate $J$-invariant planes to bilagrangian planes. The following simple observation will play a fundamental role throughout the rest of this paper.
\proclaim{Theorem \nextprocno, {\bf $J$-invariance of bilagrangian planes}}
\noindent Let $x$ be a unit-length element of $\Cal{I}$. A real plane $P\subseteq\opCl(\bil)$ is $m_x$-invariant if and only if it is $\omega_y$-lagrangian for all imaginary $y$ orthogonal to $x$.
\endproclaim
\proclabel{InvariantPlanesAreBiLagrangian}
\proof Note first that if $y$ is a unit-length element of $\Cal{I}$ orthogonal to $x$ then, since $y$ anticommutes with $x$, $(y\cdot x)$ is also a unit-length element of $\Cal{I}$ orthogonal to $x$. Suppose now that $P$ is an element of $\opInv^*(m_x)$. Let $z\in P$ be a non-trivial element such that $w:=m_x(z)\in P$ is not colinear with $z$. With $y$ as above,
$$
\omega_y(z,w) = g(z,m_y(m_x(z))) = g(z,m_{(y\cdot x)}(z)) = \omega_{(y\cdot x)}(z,z) = 0\ ,
$$
so that $\omega_y$ vanishes over $P$, as desired.
\par
Conversely, suppose that $P$ is $\omega_y$-lagrangian for all imaginary $y$ orthogonal to $x$. With $y$ again as above, for all $z,w\in P$,
$$
\omega_y(z,m_x(w)) = g(z,y\cdot x\cdot w) = \omega_{(y\cdot x)}(z,w) = 0.
$$
Since $P$ is $\omega_y$-lagrangian, in particular, it is its own $\omega_y$-annihilator. It follows that $m_x(w)\in P$ for all $w\in P$, so that $P$ is $m_x$-invariant, as desired.
\par
Finally, we verify by inspection of \eqnref{EigenspacesOfMx} that when $x$ has negative sign, and when $y$ is a unit-length element of $\Cal{I}$ orthogonal to $x$, no eigenspace of $m_x$ is $\omega_y$-lagrangian, and this completes the proof.\qed
\medskip
We conclude this section by comparing the present framework with that developed by Lychagin--Rubtsov in \cite{LychRub} (see also \cite{LychRubChekI} and \cite{Rubtsov}). Since the following constructions will not be used in the sequel, we will only sketch the main ideas, leaving the reader to complete the details. Viewing $(\opCl(\bil),g)$ as an inner-product space in its own right, we define $\alpha:\Cal{I}\times\Cal{I}\rightarrow\opso(g)$ by
$$
\alpha(x,y)\cdot z := x\cdot z - z\cdot y\ .
$$
The reader may verify that, for all $(x,y),(x',y')\in\Cal{I}\times\Cal{I}$,
$$
g(x,x') + g(y,y') = -\frac{1}{4}\opTr(\alpha(x,y)\cdot\alpha(x',y'))\ ,
$$
so that, in particular, $\alpha$ is, up to a scalar factor, a linear isomorphism. Since $g$ is non-degenerate, we may also define the linear isomorphism $\beta:\opso(g)\rightarrow\Lambda^2\opCl(\bil)$ by
$$
\beta(M) := g(\cdot,M\cdot)\ .
$$
The metric $g$ also induces a non-degenerate inner-product over $\Lambda^2\opCl(\bil)$, and the reader may verify that, for all $M,M'\in\opso(\opCl(\bil))$,
$$
g(\beta(M),\beta(M')) = -\frac{1}{2}\opTr(M\cdot M')\ .
$$
Note now that the Hodge $*$ operator define a negative pseudo-involution of $\Lambda^2\opCl(\bil)$, and the reader may also verify that its $(+1)$ and $(-1)$ eigenspaces are respectively the images under $(\beta\circ\alpha)$ of $\Cal{I}\times\{0\}$ and $\{0\}\times\Cal{I}$. Combining these facts, it follows that, for all $x,x'\in\Cal{I}$,
$$
\omega_x\wedge\omega_{x'} = g(\omega_x,*\omega_{x'})\opdVol
=g((\beta\circ\alpha)(x,0),*(\beta\circ\alpha)(x',0))\opdVol
=2g(x,x')\opdVol\ ,
$$
so that, in the terminology of Lychagin--Rubtsov, the pair $(\omega_x,\omega_{x'})$ of symplectic forms is effective if and only if the elements $x$ and $x'$ are $g$-orthogonal to one another. Note now that, in Theorem \procref{InvariantPlanesAreBiLagrangian}, the space of imaginary elements of $\opCl(\bil)$ orthogonal to $x$ is $2$-dimensional, and therefore has an orthonormal basis $(y_1,y_2)$ of cardinality $2$. It follows by the preceding discussion that the pair $(\omega_{y_1},\omega_{y_2})$ is effective, whilst, by Theorem \procref{InvariantPlanesAreBiLagrangian}, a real plane is $m_x$-invariant if and only it is $(\omega_{y_1},\omega_{y_2})$-bilegendrian.
\newsubhead{Invariant planes}[InvariantPlanes]
We now study the restrictions of $g$ and $\hat{g}$ to invariant planes. Let $x$ be a unit-length element of $\opCl^-(\bil)$. We first study the eigenspaces of $m_x$ in the case where $x$ has negative sign.
\proclaim{Lemma \nextprocno}
\noindent If $x$ is a unit-length element of $\opCl^-(\bil)$ of negative sign, then
\medskip
\myitem{(1)} the eigenspaces of $m_x$ are $g$-null,
\medskip
\myitem{(2)} the eigenspaces of $m_x$ are $\hat{g}$-orthogonal to one-another, and
\medskip
\myitem{(3)} $\pi_+$ restricts to isometries from $(E_\pm,\hat{g})$ to $(\opCl^+(\bil),2g)$.
\endproclaim
\proclabel{OrthogonalityAndNullityOfEigenspaces}
\proof The first and second assertions follow respectively from the $g$-antisymmetry and $\hat{g}$-symmetry of $m_x$. The third assertion follows by inspection of \eqnref{EigenspacesOfMx} with $y=\qk$, and this completes the proof.\qed
\medskip
\noindent We now study elements of $\opInv^*(m_x)$. We will say that an element $P$ of $\opInv^*(m_x)$ is {\emph exceptional} whenever the restriction of $\hat{g}$ to this plane degenerates. We first classify exceptional invariant planes.
\proclaim{Lemma \nextprocno}
\noindent Let $x$ be a unit-length element of $\opCl^-(\bil)$ of positive sign, and let $P$ be an element of $\opInv^*(m_x)$. The restriction of $\hat{g}$ to $P$ degenerates if and only if it vanishes identically. Furthermore, when this holds, either
\medskip
\myitem{(1)} $P$ is the direct sum of a $g$-null line of $\opCl^+(\bil)$ with a $g$-null line of $\opCl^-(\bil)$, or
\medskip
\myitem{(2)} $P$ is the graph of $\pm m_{(x\cdot\qk)}$ over $\opCl^+(\bil)$.
\endproclaim
\proclabel{ExceptionalPlanesI}
\remark[ExceptionalPlanesI] Note that $\opCl^+(\bil)$ and $\opCl^-(\bil)$ only have $g$-null lines when $\bil$ has negative sign. Furthermore, when this holds, $(1)$ accounts for precisely $2$ planes in $\opCl(\bil)$ since, although $4$ different planes may be obtained by the different direct sums of these null lines, only two of them will be elements of $\opInv^*(m_x)$.
\medskip
\proof Suppose that $\hat{g}$ degenerates over $P$. There exists a non-trivial element $z\in P$ such that, for all other $w\in P$,
$$
\hat{g}(z,w) = 0\ .
$$
Since $m_x$ is $\hat{g}$-symmetric, for all such $w$,
$$
\hat{g}(m_x(z),w) = \hat{g}(z,m_x(w)) = 0\ .
$$
Since $x$ has positive sign, $m_x$ defines a complex structure over $P$. In particular, $z$ and $m_x(z)$ are linearly independent, so that $\hat{g}$ vanishes identically, as desired.
\par
Suppose now that $\hat{g}$ vanishes identically over $P$. Let $\alpha:=z+x\cdot w$ be a non-trivial element of $P$, where $z,w\in\opCl^+(\bil)$. Straightforward calculations yield
$$
\|z+x\cdot w\|_{\hat{g}}^2 = \|z\|_g^2 - \|w\|_g^2\ ,
$$
and
$$
\hat{g}(x\cdot(z+x\cdot w),(z+x\cdot w)) =-2g(z,w)\ .
$$
Since $\hat{g}$ vanishes identically over $P$, it follows that
$$
\|z\|_g^2 - \|w\|_g^2 = g(z,w) = 0\ .
$$
There are now two possibilities. Either $z$ and $w$ are both $g$-null, in which case $P$ is the direct sum of a $g$-null line of $\opCl^+(\bil)$ with a $g$-null line of $\opCl^-(\bil)$. Otherwise, $z$ and $w$ have the same, non-zero $g$-norm-squared and are $g$-orthogonal to one another, so that $w=\pm\qk\cdot z$, and $P$ is the graph of $\pm m_{(x\cdot\qk)}$ over $\opCl^+(\bil)$, as desired.\qed
\proclaim{Lemma \nextprocno}
\noindent Let $x$ be a unit-length element of $\opCl^-(\bil)$ of negative sign, and let $P$ be an element of $\opInv^*(m_x)$. The restriction of $\hat{g}$ to $P$ degenerates if and only if $P$ contains a $\hat{g}$-null real eigenvector of $m_x$.
\endproclaim
\proclabel{ExceptionalPlanesII}
\remark[ExceptionalPlanesII] Note that it is only possible for $m_x$ to have $\hat{g}$-null real eigenvectors when $\bil$ has negative sign. In particular, when $\bil$ has positive sign, $\hat{g}$ never degenerates over $P$.
\medskip
\proof Indeed, by Lemma \procref{InvariantPlanes}, $P=\langle z,w\rangle$, where $z$ and $w$ are eigenvectors of $m_x$ with respective eigenvalues $(-1)$ and $(+1)$. By Item $(2)$ of Lemma \procref{OrthogonalityAndNullityOfEigenspaces}, $z$ and $w$ are $\hat{g}$-orthogonal to one another, so that $\hat{g}$ degenerates over $P$ if and only if one of $z$ or $w$ is $\hat{g}$-null, as desired.\qed
\medskip
Having classified exceptional invariant planes, we now study the properties of non-exceptional invariant planes.
\proclaim{Lemma \& Definition \nextprocno}
\noindent Let $x$ be a unit-length element of $\opCl^-(\bil)$ of positive sign, and let $P$ be a non-exceptional element of $\opInv^*(m_x)$. There exist precisely $4$ vectors $v\in P$ such that $(v,m_x(v))$ is $g$-orthogonal and $\hat{g}$-orthonormal. We call any such vector a {\emph principal vector} of $E$.
\endproclaim
\proclabel{PrincipalDirectionsI}
\proof Note first that, since $m_x$ is $g$-antisymmetric, $g(\cdot,m_x\cdot)$ vanishes identically over $P$. Consider now the bilinear form $h:=\hat{g}(\cdot,m_x\cdot)$. Since $m_x$ is $\hat{g}$-symmetric, $h$ is also symmetric. Since $m_x$ is a $\hat{g}$-anti-isometry, it is also an $h$-anti-isometry. Since $P$ is non-exceptional, by Lemma \procref{ExceptionalPlanesI}, $\hat{g}$, and therefore also $h$, is non-degenerate. Since $h$ admits a non-trivial anti-isometry, it has mixed signature and thus has $2$ null directions. Furthermore, neither of these null directions is $\hat{g}$-null, for otherwise $\hat{g}$ would vanish identically over $P$. We readily verify that any $h$-null vector $v$ with unit $\hat{g}$-norm-squared has the desired properties, and the result follows.\qed
\proclaim{Lemma \& Definition \nextprocno}
\noindent Let $x$ be a unit-length element of $\opCl^-(\bil)$ of negative sign, and let $P$ be a non-exceptional element of $\opInv^*(m_x)$. There exist precisely $8$ vectors $v\in P\otimes\Bbb{C}$ such that $(v,m_x(v))$ is $g$-orthogonal and $\hat{g}$-orthonormal. Furthermore, $4$ of these vectors are real unless $\bil$ has negative sign and $P$ contains eigenvectors of $m_x$ of opposing sign. We call any such vector a {\emph principal direction} of $E$.
\endproclaim
\proclabel{PrincipalDirectionsII}
\remark[PrincipalDirectionsII] We will not be explicitly concerned with the case of imaginary principal directions in the sequel.
\medskip
\proof Let $z,w\in\opCl^+(\bil)$ be such that $\hat{z}:=(z,m_x(z))$ and $\hat{w}:=(w,m_x(w))$ are the $m_x$-eigenvectors which generate $P$. Since $m_x$ is $\hat{g}$-symmetric, these two eigenvectors are $\hat{g}$-orthogonal to one another. Furthermore, since $P$ is non-exceptional, neither vector is $\hat{g}$-null, and we may therefore suppose that they each have unit norm-squared. If their norms-squared have the same sign, then the vectors $(\hat{z}\pm\hat{w})/\sqrt{2}$ have the desired properties. Otherwise, their norms-squared have opposing signs, and the vectors $(\hat{z}\pm i\hat{w})/\sqrt{2}$ have the desired properties. The remaining vectors are obtained by multiplying these vectors by $(\pm1)$ and $(\pm i)$, and the result follows.\qed
\proclaim{Lemma \& Definition \nextprocno}
\noindent Let $x$ be a unit-length element of $\opCl^-(\bil)$ of negative sign, let $P$ be a non-exceptional element of $\opInv^*(m_x)$, and let $v$ be a real principal direction of $P$. The plane $E_v:=\langle v,\qk\cdot x\cdot v\rangle$ has $1$-dimensional intersections with each of $\opCl^+(\bil)$ and $\opCl^-(\bil)$. We call the intersection $E_v\cap\opCl^+(\bil)$ the {\emph principal line} of $P$ associated to $v$.
\endproclaim
\proclabel{PrincipalPlanes}
\remark[PrincipalPlanes] An analogous result also holds for imaginary principal directions, although we will have no use for this in the present paper.
\medskip
\proof Indeed, let $z,w\in\opCl^+(\bil)$ be such that $v=z+x\cdot w$. A straightforward calculation yields
$$
\hat{g}(v,m_x(v)) = (1-\|x\|_g^2)g(z,w) = 2g(z,w)\ .
$$
Since $v$ and $m_x(v)$ are $\hat{g}$-orthogonal to one another, it follows that $z=\pi_+(v)$ and $w=x\cdot\pi_-(v)$ are $g$-orthogonal to one-another. In particular, $z$ and $(\qk\cdot w)$ are linearly dependent, and the result follows.\qed
\newhead{Legendrian Clifford structures}[CliffordStructuresAndExtrinsicCurvature]
\newsubhead{Overview}[OverviewLegendrianCliffordStructures]
We now apply the algebraic theory developed in Section \headref{CliffordStructures} to the setting of contact geometry. Recall that the unit sphere bundle $UX$ of any riemannian $3$-manifold $X$ carries a natural contact structure $\alpha$, with contact bundle $W$, and symplectic form $\omega_i:=d\alpha$. We will construct natural families of Clifford structures over $W$ such that, for all $k\in\Bbb{R}\setminus\{0\}$, there exists a pseudo-involution $J_k$, associated with a fibrewise symplectic form $\omega_k$, having the property that $(S,e)$ is a CEC-$k$ surface if and only if its Gauss lift $(S,\hat{e})$ is $(\omega_i,\omega_k)$-bilegendrian, and thus also $J_k$-invariant.
\par
When $k>0$, the pseudo-involution $J_k$ is nothing other than the Monge-Amp\`ere structure introduced by Labourie in \cite{LabII}, showing how his framework embeds within our own. When the ambient space is either $\Bbb{H}^3$ or $\Bbb{S}^3$, and when $k$ is respectively equal to $(+1)$ or $(-1)$, $(\omega_i,\omega_k)$-bilegendrian surfaces are respectively flat fronts of the first and second types discussed in the introduction.
\par
The main result of this section is Theorem \procref{FlatMetric}, which proves that, whenever $X$ is a space form, every immersed bilegendrian surface in $UX$ carries a canonical flat metric over the complement of a certain set of exceptional points. We underline, above all, how this follows in an elementary manner from the many symmetries of the Clifford structure. Finally, in Sections \subheadref{SurfacesOfConstantPositiveCurvature} and \subheadref{SurfacesOfConstantNegativeCurvature} we show how this metric is related to well-known structures over CEC surfaces, and in Section \subheadref{HilbertsTheorem}, we use this metric to prove Hilbert's theorem concerning the non-existence of complete, negatively-curved CEC surfaces in euclidean $3$-space.
\newsubhead{The Clifford structure}[AnExplicitCliffordStructure]
We now describe the explicit Clifford structure over $\Bbb{R}^2\times\Bbb{R}^2$ which we will use to study surfaces of constant extrinsic curvature. Let $\bil$ be a non-degenerate, symmetric bilinear form over $\Bbb{R}^2$, and let $A$ be an element of $\opso(\bil)\cap\hat{\opO}(\bil)$. Note that
$$
\opVol := -\bil(\cdot,A\cdot)\eqnum{\nexteqnno[DefinitionOfA]}
$$
is a unit-length alternating form, so that, in particular, the choice of $A$ also defines a preferred orientation of $\Bbb{R}^2$.
\par
Given $\eta\in\{\pm 1\}$, we define the {\emph symplectic form} $\omega_i$, the {\emph inner-product} $g^\eta$, and the {\emph grade involution} $\alpha$ over $\Bbb{R}^2\times\Bbb{R}^2$ by
$$\eqalign{
\omega_i((\xi,\mu),(\xi',\mu')) &:= \bil(\xi,\mu') - \bil(\mu,\xi')\ ,\cr
g^\eta((\xi,\mu),(\xi',\mu')) &:= \bil(\xi,\xi') + \eta\bil(\mu,\mu')\ ,\ \text{and}\cr
\alpha(\xi,\mu) &:= (\xi,-\mu)\ .\cr}
$$
We also define the {\emph complementary inner product} $\hat{g}^\eta$ by $g^\eta(\cdot,\alpha\cdot)$, so that
$$
\hat{g}^\eta((\xi,\mu),(\xi',\mu')) := \bil(\xi,\xi') - \eta\bil(\mu,\mu')\ .
$$
\par
Recall now that geometric structures can also be defined in terms of their stabilizer subgroup. Let $\opStab(\omega_i,g^\eta,\alpha)$ denote the stabilizer subgroup of the subgroup of the triple $(\omega_i,g^\eta,\alpha)$ in $\opEnd(\Bbb{R}^2\times\Bbb{R}^2)$.
\proclaim{Lemma \nextprocno}
\noindent $\opStab(\omega_i,g^\eta,\alpha)$ consists of all matrices $M$ of the form
$$
M = \pmatrix N\hfill& 0\hfill\cr 0\hfill& N\hfill\cr\endpmatrix\ ,\eqnum{\nexteqnno[DescriptionOfFix]}
$$
where $N\in\opO(\bil)$.
\endproclaim
\proclabel{DescriptionOfFix}
\proof Indeed, choose $M\in\opStab(\omega_i,\hat{g}^\eta,\alpha)$. Since $M$ preserves $\alpha$,
$$
M = \pmatrix N\hfill& 0\hfill\cr 0\hfill& N'\hfill\cr\endpmatrix\ ,
$$
for some $N,N'\in\opEnd(\Bbb{R}^2)$. Since $M$ preserves $g^\eta$,
$$
N,N'\in\opO(\bil)\ .
$$
Finally, since $M$ preserves $\omega_i$, denoting by $N^*$ the $\bil$-dual of $N$,
$$
N' = (N^*)^{-1} = N\ ,
$$
so that $M$ has the desired form. Conversely, we verify that every linear map of this form preserves $(\omega_i,g^\eta,\alpha)$, and this completes the proof.\qed
\medskip
\noindent We now show how the triple $(\omega_i,g^\eta,\alpha)$ defines a Clifford structure over $\Bbb{R}^2\times\Bbb{R}^2$.
\proclaim{Lemma \nextprocno}
\noindent Suppose that $\bil$ has a positive direction $v$, and let $\qi$ denote the corresponding element of $\opCl(\eta\bil)$. There exists a linear isomorphism $\phi:\Bbb{R}^2\times\Bbb{R}^2\rightarrow\opCl(\eta\bil)$ sending $\omega_i$, $g^\eta$ and $\alpha$ respectively to $\omega_i$, $g$ and the grade involution. Furthermore, $\phi$ is unique up to precomposition by an element of $\opStab(\omega_i,g^\eta,\alpha)$.
\endproclaim
\proclabel{ExistenceOfCliffordIsomorphism}
\proof Since $\qk$ has the same sign as $\bil$, there exists an isometry $\phi_0:(\Bbb{R}^2,\bil)\rightarrow(\opCl^+(\eta\bil),g)$. We define $\phi:\Bbb{R}^2\times\Bbb{R}^2\rightarrow\opCl(\eta\bil)$ by
$$
\phi(\xi,\mu) := \phi_0(\xi) + \qi^{-1}\cdot\phi_0(\mu)\ .
$$
This isomorphism trivially sends $\alpha$ to the grade involution. Furthermore, for all $(\xi,\mu),(\xi',\mu')\in\Bbb{R}^2\times\Bbb{R}^2$,
$$
\omega_i(\phi(\xi,\mu),\phi(\xi',\mu'))=g(\phi_0(\xi) + \qi^{-1}\cdot\phi_0(\mu),\qi\cdot\phi_0(\xi') + \phi_0(\mu'))\ .
$$
Since $\opCl^\pm(\eta\bil)$ are $g$-orthogonal to one another, and since $m_{\qi}$ is $g$-antisymmetric, it follows that
$$
\omega_i(\phi(\xi,\mu),\phi(\xi',\mu'))=g(\phi_0(\xi),\phi_0(\mu')) - g(\phi_0(\mu),\phi_0(\xi'))
=\bil(\xi,\mu') - \bil(\mu,\xi')=\omega_i((\xi,\mu),(\xi',\mu'))\ ,
$$
as desired. Similarly, for all $(\xi,\mu),(\xi',\mu')\in\Bbb{R}^2\times\Bbb{R}^2$,
$$
g(\phi(\xi,\mu),\phi(\xi',\mu')) = g(\phi_0(\xi)+\qi^{-1}\cdot\phi_0(\mu),\phi_0(\xi') + \qi^{-1}\cdot\phi_0(\mu))=g(\phi_0(\xi),\phi_0(\xi')) + \|\qi\|_g^2 g(\phi_0(\mu),\phi_0(\mu'))\ .
$$
Since
$$
\|\qi\|_g^2 = -\qi^2 = \eta\bil(e_1,e_1) = \eta\ ,
$$
it follows that
$$
g(\phi(\xi,\mu),\phi(\xi',\mu')) = g^\eta((\xi,\mu),(\xi',\mu'))\ ,
$$
as desired. This proves existence, uniqueness follows by definition of $\opStab(\omega_i,g^\eta,\alpha)$, and this completes the proof.\qed
\medskip
We now pull back through $\phi$ the objects constructed in Section \headref{CliffordStructures} to obtain what we call a {\sl Clifford structure} over $\Bbb{R}^2\times\Bbb{R}^2$. Since $\opO(\bil)$ has more than one connected component\numberedfootnote{$\opO(\bil)$ has $2$ connected components when $\bil$ has positive sign, and $4$ when $\bil$ has negative sign.}, there is, in fact, a sign ambiguity, which we resolve upon setting
$$
\omega_k^\eta((\xi,\mu),(\xi',\mu')) := \opVol(\xi,\mu) - \eta\opVol(\xi',\mu')\ .
$$
Explicitly, the Clifford structure is given by the matrices $I^\eta$, $J^\eta$, and $K$, defined by
$$
\omega_i =: g^\eta(\cdot,I^\eta\cdot)\ ,\
\omega_k^\eta =: \omega_i(\cdot,J^\eta\cdot)\ ,\ \text{and}\
\omega_k^\eta =: g^\eta(\cdot,K\cdot)\ ,
$$
and we also denote
$$
\hat{I}^\eta := \alpha\circ I^\eta\ ,\ \hat{J}^\eta := \alpha\circ J^\eta\ ,\ \text{and}\ \hat{K}:=\alpha\circ K\ .
$$
We readily verify that
$$
I^\eta = \pmatrix 0\hfill& \opId\hfill\cr -\eta\opId\hfill& 0\hfill\cr\endpmatrix\ ,\
J^\eta = \pmatrix 0\hfill& \eta A\hfill\cr A\hfill& 0\hfill\cr\endpmatrix\ ,\ \text{and}\
K = \pmatrix A\hfill& 0\hfill\cr 0\hfill& -A\hfill\cr\endpmatrix\ ,\eqnum{\nexteqnno[ExplicitFormulaOfIJK]}
$$
and we verify by inspection that these matrices satisfy the algebraic relations of a Clifford structure, namely
$$
I^\eta J^\eta = K\ ,
$$
their squares satisfy
$$
(I^\eta)^2 = -\eta\opId\ ,\ (J^\eta)^2 = -\eta\epsilon\opId\ ,\ \text{and}\ K^2 = -\epsilon\opId\ ,
$$
and their anti-commutators satisfy
$$
\{I^\eta,J^\eta\} = \{I^\eta,K\} = \{J^\eta,K\} = 0\ .
$$
\newsubhead{The Clifford bundle}[TheCliffordBundle]
We now extend the construction of the preceding section to define Clifford structures over bundles. This formalism can, in fact, be fruitfully applied to bundles over general riemannian manifolds (see, for example, \cite{LabII} and \cite{SmiQMAK}). However, we will limit ourselves here to a simple case where all objects in question will be covariant constant. This will nonetheless be sufficient to derive interesting properties of CEC surfaces in space-forms of riemannian and lorentzian signature.
\par
Let $\bil$ be a non-degenerate, symmetric bilinear form over $\Bbb{R}^4$. For $\eta\in\{\pm 1\}$, we define the {\emph symplectic form} $\omega_i$, the {\emph inner product} $g^\eta$, and the {\emph grade involution} $\alpha$ over $\Bbb{R}^4\times\Bbb{R}^4$ by
$$\eqalign{
\omega_{i,(x,y)}((\xi,\mu),(\xi',\mu')) &= \bil(\xi,\mu') - \bil(\mu,\xi')\ ,\cr
g^\eta((\xi,\mu),(\xi',\mu')) &:= \bil(\xi,\xi') + \eta \bil(\mu,\mu')\ ,\ \text{and}\cr
\alpha(\xi,\mu) &:= (\xi,-\mu)\ .\cr}
$$
We likewise define the {\emph complementary inner product} $\hat{g}^\eta$ by $\hat{g}^\eta:=g(\cdot,\alpha\cdot)$, so that
$$
\hat{g}^\eta((\xi,\mu),(\xi',\mu')) := \bil(\xi,\xi') - \eta \bil(\mu,\mu')\ .
$$
\par
We define the open, codimension $1$ submanifold $M$ of $\Bbb{R}^4\times\Bbb{R}^4$ by
$$
M:=\{(x,y)\in\Bbb{R}^4\times\Bbb{R}^4\ |\ \bil(x,x),\bil(y,y)\neq 0\ ,\ \bil(x,y)=0\}\ .
$$
We define the $4$-dimensional subbundle $W$ of $TM$ such that, at every point $(x,y)$ of $M$,
$$
W_{(x,y)} := \{ (\xi,\mu)\in\Bbb{R}^4\times\Bbb{R}^4\ |\ \bil(\xi,x)=\bil(\xi,y)=\bil(\mu,x)=\bil(\mu,y)=0\}\ .
$$
At each point $(x,y)$ of $M$, the triple $(\omega_i,g^\eta,\alpha)$ yields a Clifford structure over $W_{(x,y)}$ as in Section \subheadref{AnExplicitCliffordStructure}, and we resolve the sign ambiguity by setting,
$$
\omega_k^\eta((\xi,\mu),(\xi',\mu')) := \opVol(x,\xi,\mu,y) - \eta\opVol(x,\xi,\mu,y)\ ,
$$
where $\opVol$ here denotes the canonical volume form of $\Bbb{R}^4$. Let $I^\eta$, $J^\eta$ and $K$ denote the pseudo-involutions of $W$ obtained as in \eqnref{ExplicitFormulaOfIJK}, and let $\hat{I}^\eta$, $\hat{J}^\eta$ and $\hat{K}^\eta$ denote their respective compositions with the grade-involution $\alpha$.
\par
We now show that this Clifford structure is covariant constant. To this end, we define the $4$-dimensional bundle $N$ over $M$ such that, at every point $(x,y)$ of $M$,
$$
N_{(x,y)} := \langle (x,0),(y,0),(0,x),(0,y)\rangle\ .
$$
Since $N$ is transverse to $W$ in $\opT(\Bbb{R}^4\times\Bbb{R}^4)$, composing the canonical derivative of $\Bbb{R}^4\times\Bbb{R}^4$ with projection along $N$ yields a covariant derivative over $W$ which we denote by $\nabla$.\numberedfootnote{Although it is usual in riemannian geometry to restrict covariant derivatives to subbundles by composing with orthogonal projection, it is in fact sufficient to project along any transverse subbundle, as is standard practice in affine geometry (see, for example, \cite{SasakiNomizu})} Observe that, at every point $(x,y)$ of $M$, the two bundles $W$ and $N$ are $\omega_i$-, $g^\eta$- and $\hat{g}^\eta$-orthogonal to one another, so that, for all sections $\sigma$ and $\tau$ of $W$ over $M$, and for every tangent vector field $\xi$ over $M$,
$$
(\nabla_\xi\omega_i)(\sigma,\tau) = (\nabla_\xi g^\eta)(\sigma,\tau) = (\nabla_\xi\hat{g}^\eta)(\sigma,\tau) = 0\ .
$$
The analogous result also holds for $\omega_k^\eta$, although it is less immediate.
\proclaim{Lemma \nextprocno}
\noindent For all sections $\sigma$ and $\tau$ of $W$ over $M$, for every tangent vector field $\xi$ over $M$ taking values in $W$, and for each $\eta\in\{\pm 1\}$,
$$
(\nabla_\xi\omega_k^\eta)(\sigma,\tau) = 0\ .\eqnum{\nexteqnno[OmegaKIsFlat]}
$$
\endproclaim
\proclabel{OmegaKIsFlat}
\proclaim{Corollary \nextprocno}
\noindent For every section $\sigma$ of $W$ over $M$, and for every tangent vector field $\xi$ over $M$ taking values in $W$,
$$
(\nabla_\xi I^\eta)\sigma = (\nabla_\xi J^\eta)\sigma = (\nabla_\xi K)\sigma = (\nabla_\xi\alpha)\sigma = 0\ .\eqnum{\nexteqnno[IJKAreCovariantConstant]}
$$
\endproclaim
\proclabel{IJKAreCovariantConstant}
{\bf\noindent Proof of Lemma \procref{OmegaKIsFlat}:\ }It suffices to prove the result for $\omega_k':=\omega_k^{+1}+\omega_k^{-1}$, as the proof for $\omega_k'':=\omega_k^{+1}-\omega_k^{-1}$ is identical. Note first that, since $\xi$ takes values in $W$,
$$
D_\xi x, D_\xi y \in \langle x,y\rangle^\perp\ .
$$
Since $\pi_1(\sigma)$ and $\pi_1(\tau)$ also take values in $\langle x,y\rangle^\perp$, and since this space is $2$-dimensional, it follows that
$$
\opVol(\pi_1(\sigma),\pi_1(\tau),D_\xi x,y) = \opVol(\pi_1(\sigma),\pi_1(\tau),x,D_\xi y) = 0\ .
$$
Thus,
$$
D_\xi\big(\opVol(\pi_1(\sigma),\pi_1(\tau),x,y)\big)
=\opVol(\pi_1(D_\xi\sigma),\pi_1(\tau),x,y) + \opVol(\pi_1(\sigma),\pi_1(D_\xi\tau),x,y)\ .
$$
Let $p$ denote the projection onto $N$ along $W$, and note that
$$
\pi_1(D_\xi\sigma) - \pi_1(\nabla_\xi\sigma) = \pi_1(p(D_\xi\sigma)) \in \langle x,y\rangle\ ,
$$
with a similar identity holding for $\tau$. It follows that
$$\eqalign{
D_\xi\omega_k'(\sigma,\tau)
&=2D_\xi\big(\opVol(\pi_1(\sigma),\pi_1(\tau),x,y)\big)\cr
&=2\opVol(\pi_1(\nabla_\xi\sigma),\pi_1(\tau),x,y) + 2\opVol(\pi_1(\sigma),\pi_1(\nabla_\xi\tau),x,y)\cr
&=\omega_k'(\nabla_\xi\sigma,\tau) + \omega_k'(\sigma,\nabla_\xi\tau)\ ,\cr}
$$
so that
$$
\big(\nabla_\xi\omega_k')(\sigma,\tau) = 0\ ,
$$
as desired.\qed
\newsubhead{Bilegendrian surfaces}[IntroducingBilegendrianSurfaces]
We now introduce the class of surfaces of interest to us. We define an {\emph immersed surface} in $M$ to be a pair $(S,\phi)$, where $S$ is a smooth surface, and $\phi:S\rightarrow M$ is a smooth immersion. We define an {\emph integral surface} of $W$ to be an immersed surface $(S,\phi)$ whose every tangent plane is contained in $W$, that is, such that $\opIm(D\phi(x))\subseteq W_{\phi(x)}$ for all $x\in S$. We define a $(\omega_i,\omega^\eta_k)$-{\emph bilegendrian surface} to be an integral surface of $W$ such that
$$
\phi^*\omega_i=\phi^*\omega^\eta_k=0\ .
$$
In what follows, we will refer to such surfaces simply as {\emph bilegendrian surfaces} when no ambiguity arises. By Lemma \procref{InvariantPlanesAreBiLagrangian}, an integral surface of $W$ is a bilegendrian surface if and only if every one of its tangent planes is $J^\eta$-invariant. We will say that a bilegendrian surface is of {\emph elliptic type} (resp {\emph hyperbolic type}) whenever $J^\eta$ is positive (resp. negative). Note that bilegendrian surfaces of elliptic type are immersed $J^\eta$-holomorphic curves, justifying our terminology.
\par
Given a bilegendrian surface $S$, we define its {\emph cubic form} by
$$
C(\xi,\mu,\nu) := \omega_i(\nabla_\xi D\phi\cdot\mu,D\phi\cdot\nu)\ .
$$
Since $\omega_i$ vanishes over every tangent plane of $S$, its cubic form is a tensor. It will also be useful to define its {\emph complementary cubic form} by
$$
\hat{C}(\xi,\mu,\nu) := \omega_i(\nabla_\xi D\phi\cdot\mu,\alpha(D\phi\cdot\nu))\ ,
$$
even though this object is not a tensor. The cubic form plays the same role in lagrangian and legendrian geometries as the second fundamental form does in riemannian hypersurface geometry. Indeed,
$$
C(\xi,\mu,\nu) = \hat{g}^\eta(\nabla_\xi D\phi\cdot\mu,\hat{I}^\eta D\phi\cdot\nu)\ ,
$$
and since $\hat{I}^\eta$ sends every $\hat{J}^\eta$-invariant plane to its $\hat{g}^\eta$-orthogonal complement, $C$ indeed encodes the normal component of $\nabla_\xi D\phi\cdot\mu$ in $W$ at every $\hat{g}$-non-degenerate point of $S$.
\proclaim{Lemma \nextprocno}
\noindent The cubic form $C$ is symmetric.
\endproclaim
\proclabel{SymmetryOfCubicForm}
\proof Indeed, since $S$ is $\omega_i$-lagrangian,
$$\eqalign{
C(\xi,\mu,\nu) &= \omega_i(\nabla_\xi D\phi\cdot\mu,D\phi\cdot\nu)\cr
&=\omega_i(\nabla_\mu D\phi\cdot\xi + D\phi\cdot[\xi,\mu],D\phi\cdot\nu)\cr
&=\omega_i(\nabla_\mu D\phi\cdot\xi,D\phi\cdot\nu)\cr
&=C(\mu,\xi,\nu)\ .\cr}
$$
Likewise, since $\omega_i$ is covariant constant and vanishes over $TS$,
$$\eqalign{
C(\xi,\mu,\nu) &= \omega_i(\nabla_\xi D\phi\cdot\mu,D\phi\cdot\nu)\cr
&=D_\xi\omega_i(D\phi\cdot\mu,D\phi\cdot\nu) - \omega_i(D\phi\cdot\mu,\nabla_\xi D\phi\cdot\nu)\cr
&=\omega_i(\nabla_\xi D\phi\cdot\nu,D\phi\cdot\mu)\cr
&=C(\xi,\nu,\mu)\ .\cr}
$$
Since the transpositions $(12)$ and $(23)$ generate the group of permutations of $\{1,2,3\}$, it follows that $C$ is symmetric, as desired.\qed
\proclaim{Lemma \nextprocno}
\noindent $J^\eta$ is $C$-symmetric. That is, for all $\xi$, $\mu$ and $\nu$,
$$
C(J^\eta\xi,\mu,\nu) = C(\xi,J^\eta\mu,\nu) = C(\xi,\mu,J^\eta\nu)\ .\eqnum{\nexteqnno[CSymmetryOfJ]}
$$
\endproclaim
\proclabel{CSymmetryOfJ}
\proof Indeed, since $\nabla J^\eta=0$,
$$
C(\xi,J^\eta\mu,\nu)
=g^\eta(\nabla_\xi J^\eta\mu,I^\eta\nu)
=g^\eta(J^\eta\nabla_\xi\mu,I^\eta\nu)
=-g^\eta(\nabla_\xi\mu,J^\eta I^\eta\nu)
=g^\eta(\nabla_\xi\mu,I^\eta J^\eta\nu)
=C(\xi,\mu,J^\eta\nu)\ ,
$$
and the result follows by symmetry of $C$.\qed
\proclaim{Lemma \nextprocno}
\noindent Let $\xi$ and $\mu$ be tangent vector fields over $S$. If $\mu$ is a field of $J^\eta$-eigenvectors, then
$$
\hat{C}(\xi,\mu,\mu) = 0\ .\eqnum{\nexteqnno[VanishingValuesOfCCF]}
$$
\endproclaim
\proclabel{VanishingValuesOfCCF}
\proof Let $\epsilon\in\{\pm1\}$ denote the $J^\eta$-eigenvalue of $\mu$. Since $\hat{g}$ and $I$ are covariant constant, and since $I$ is $\hat{g}$-symmetric,
$$
\hat{g}(\nabla_\xi D\phi\cdot\mu,I\cdot D\phi\cdot\mu)
=\frac{1}{2}D_\xi\hat{g}(D\phi\cdot\mu,I\cdot D\phi\cdot\mu)
=\frac{\epsilon}{2}D_\xi\hat{g}(D\phi\cdot\mu,I\cdot J\cdot D\phi\cdot\mu)
=\frac{\epsilon}{2}D_\xi\hat{g}(D\phi\cdot\mu,K\cdot D\phi\cdot\mu)\ .
$$
Since $K$ is $\hat{g}$-antisymmetric, this vanishes, and the result follows.\qed
\medskip
Let $S$ be a bilegendrian surface. We say that a point $x\in S$ is {\emph exceptional} whenever $\hat{g}^\eta$ degenerates over the tangent space of $S$ at this point. We denote the set of exceptional points of $S$ by $\opExc(S)$. By Lemmas \procref{PrincipalDirectionsI} and \procref{PrincipalDirectionsII}, in a neighbourhood of every non-exceptional point, we may define a pair $(E_1,E_2)$ of (possibly complex) fields of principal vectors, which is unique up to reordering and change of sign. Furthermore, by Lemma \procref{PrincipalPlanes}, the principal vectors likewise define a pair $(L_1,L_2)$ of line-valued functions such that, at each point, $L_1(x)$ and $L_2(x)$ are orthogonal to one-another and also to the plane $\langle(\pi_1\circ\phi)(x),(\pi_2\circ\phi)(x)\rangle$. We will call these functions {\emph principal lines} of $S$. These will be of use in the sequel.
\par
Finally, when the principal directions are real, they integrate to yield two transverse foliations of $S\setminus\opExc(S)$ whose leaves we will call {\emph lines of curvature} of $S$.
\proclaim{Theorem \nextprocno, {\bf Existence of a flat metric}}
\noindent The metric $\phi^*\hat{g}^\eta$ is flat over $S\setminus\opExc(S)$, and the lines of curvature, when they exist, are geodesics of this metric.
\endproclaim
\proclabel{FlatMetric}
\proof Let $\hat{\nabla}^\eta$ denote the Levi-Civita covariant derivative of $\phi^*\hat{g}^\eta$ over $S\setminus\opExc(S)$. By definition of the principal vectors, we may suppose that $E_2=J^\eta E_1$. By $\hat{g}^\eta$-symmetry of $J^\eta$, for all $\xi$ tangent to $S$,
$$
\hat{g}^\eta (\nabla_\xi E_1,E_2) = \hat{g}^\eta (\nabla_\xi E_1,J^\eta\cdot E_1)
=\frac{1}{2}D_\xi \hat{g}^\eta (E_1,J^\eta\cdot E_1) - \hat{g}^\eta (E_1,(\nabla_\xi J^\eta)\cdot E_1)=0\ .
$$
On the other hand,
$$
\hat{g}^\eta (\nabla_\xi E_1,E_1) = \frac{1}{2}D_\xi \hat{g}^\eta (E_1,E_1) = 0\ .
$$
At every point where $\hat{g}^\eta$ is non-degenerate, this yields
$$
\hat{\nabla}^\eta_\xi E_1 = 0\ .
$$
In a similar manner, we show that
$$
\hat{\nabla}^\eta_\xi E_2 = 0\ .
$$
It follows that $(E_1,E_2)$ is a covariant constant frame of $\phi^*\hat{g}^\eta$, and this completes the proof.\qed
\newsubhead{Example I - Surfaces of constant positive curvature in $\Bbb{H}^3$}[SurfacesOfConstantPositiveCurvature]
Bilegendrian surfaces may be used to study CEC surfaces of non-zero curvature in space forms. We will not provide an exhaustive treatment of this subject, and will instead limit ourselves to two illustrative cases, namely positively curved CEC surfaces in hyperbolic $3$-space $\Bbb{H}^3$, and negatively curved CEC surfaces in euclidian $3$-space $\Bbb{R}^3$.
\par
Let $\bil$ be a non-degenerate, symmetric bilinear form of signature $(3,1)$. Recall that $\Bbb{H}^3$ identifies with a connected component of the submanifold
$$
X := \{ x\in\Bbb{R}^4\ |\ \bil(x,x)=-1 \}\ .
$$
Let $(S,e)$ be an immersed surface in $X$, let $\nu_e$ denote its unit normal vector field, let $\opI_e$, $\opII_e$ and $\opIII_e$ denote its three fundamental forms, and let $A_e$ denote its shape operator. Since the tangent bundle of $X$ naturally embeds into the trivial bundle $X\times\Bbb{R}^4$, we view $\nu_e$ as a function taking values in $\Bbb{R}^4$. For $k>0$, we define the $k$-{\emph Gauss lift} of $e$ by
$$
\hat{e}_k(x) := \bigg(e(x),\frac{1}{\sqrt{k}}\nu_e(x)\bigg)\ ,
$$
and we call $(S,\hat{e}_k)$ the $k$-{\emph Gauss lift} of $(S,e)$.
\proclaim{Lemma \nextprocno}
\noindent For all $k>0$, $\hat{e}_k$ takes values in $M$. Furthermore, for all $x\in S$, and for every tangent vector $\xi$ to $S$ at $x$,
$$
D\hat{e}_k(x)\cdot\xi = \bigg(De(x)\cdot\xi,\frac{1}{\sqrt{k}}De(x)\cdot A_e(x)\cdot\xi\bigg)\ .\eqnum{\nexteqnno[DerivativeOfGaussLift]}
$$
In particular, $(S,\hat{e}_k)$ is $W$-horizontal and $\omega_i$-lagrangian.
\endproclaim
\proclabel{DerivativeOfKGaussLift}
\proof Indeed, by hypothesis
$$
\bil(e(x),e(x)) = 0\ ,\ \bil(\nu_e(x),e(x)) =0\ ,\ \text{and}\
\bil(\nu_e(x),\nu_e(x)) = 1\ ,
$$
from which it follows that the $k$-Gauss lift takes values in $M$. Differentiating the second two relations yields
$$
\bil(D\nu_e(x)\cdot\xi,e(x)) = - \bil(De(x)\cdot\xi,\nu_e(x))\ \text{and}\
\bil(D\nu_e(x)\cdot\xi,\nu_e(x)) = 0\ .
$$
However, by definition
$$
\bil(De(x)\cdot\xi,\nu_e(x)) = 0\ ,
$$
so that
$$
\bil(D\nu_e(x)\cdot\xi,e(x)) = 0\ ,
$$
and $D\nu_e(x)$ therefore takes values in $\langle e(x),\nu_e(x)\rangle^\perp$. Since $e$ is an immersion in $\Bbb{H}^3$, this subspace is equal to the image of $De(x)$. However, for all other tangent vectors $\mu$ of $S$ at $x$,
$$
\bil(D\nu_e(x)\cdot\xi,De(x)\cdot\mu) = \opII_e(\xi,\mu) = \bil(De(x)\cdot A_e(x)\cdot\xi,De(x)\cdot\mu)\ ,
$$
from which it follows that
$$
D\nu_e(x)\cdot\xi = De(x)\cdot A_e(x)\cdot\xi\ ,
$$
and this completes the proof.\qed
\proclaim{Lemma \nextprocno}
\noindent For $k>0$, $(S,e)$ has constant extrinsic curvature equal to $k$ if and only if its $k$-Gauss lift is $J^{+1}$-invariant.
\endproclaim
\proclabel{KGaussLiftIsInvariant}
\proof Indeed, $(S,e)$ has constant extrinsic curvature equal to $k$ if and only if $\opDet(A_e)=k$. By \eqnref{DerivativeOfGaussLift}, this holds if and only if the $k$-Gauss lift is $\omega^{+1}_k$-lagrangian, and the result now follows by Theorem \procref{InvariantPlanesAreBiLagrangian}.\qed
\medskip
\noindent Theorem \procref{FlatMetric} now manifests itself in the following form.
\proclaim{Theorem \nextprocno, {\bf Existence of a translation structure}}
\noindent If $(S,e)$ has constant extrinsic curvature equal to $k$, then the metric
$$
h_e := \opI_e - \frac{1}{k}\opIII_e\eqnum{\nexteqnno[FlatMetricOfKSurface]}
$$
is flat away from the set of umbilic points of $(S,e)$. Furthermore, the lines of curvature are geodesics of this metric.
\endproclaim
\proclabel{FlatMetricOfKSurface}
\remark[FlatMetricOfKSurface] This metric has mixed signature away from the set of umbilic points of $e$.
\medskip
\proof Indeed, by \eqnref{DerivativeOfGaussLift},
$$
\hat{e}_k^* g^{-\eta} = \opI_e - \frac{1}{k}\opIII_e\ ,
$$
and the result now follows by Theorem \procref{FlatMetric}.\qed
\medskip
The metric constructed in Theorem \procref{FlatMetricOfKSurface} is interpreted geometrically as follows. Let $x$ be a non-umbilic point of $S$. Upon suitably reparametrizing a neighbourhood of $x$ by an open subset $\Omega$ of $\Bbb{R}^2$, we may suppose that
$$
h_e = dx^2 - dy^2\ ,
$$
and that the horizontal and vertical lines are lines of curvature of this immersion. It is then straightforward to show that the fundamental forms of $e$ are
$$
\opI_e = \pmatrix\opCosh^2(\theta)\hfill&0\hfill\cr0\hfill&\opSinh^2(\theta)\hfill\cr\endpmatrix\ ,\
\opII_e = \frac{\sqrt{k}}{2}\pmatrix\opSinh(2\theta)\hfill&0\hfill\cr0\hfill&\opSinh(2\theta)\hfill\cr\endpmatrix\ ,\ \text{and}\
\opIII_e = k\pmatrix\opSinh^2(\theta)\hfill&0\hfill\cr0\hfill&\opCosh^2(\theta)\hfill\cr\endpmatrix\ ,
$$
for some smooth function $\theta:\Omega\rightarrow]0,\infty[$.
\par
Recall now that the Hopf differential $\psi_e$ of $e$ is the trace-free component of the first fundamental form with respect to the second fundamental form. An elementary calculation then yields
$$
\psi_e = dzdz = (dx^2-dy^2) + i(dxdy + dydx)\ .
$$
In particular, Theorem \procref{FlatMetricOfKSurface} shows that the Hopf differential is holomorphic. This observation is in itself not very interesting, since holomorphicity of $\psi_e$ is readily proven using more direct techniques. It is more intriguing, however, to observe that $h_e$ is none other than the real part $\psi_e$. In other words, $h_e$ is the metric of the well-known flat $(\Bbb{Z}_2\ltimes\Bbb{R}^2)$-structure defined over the complement of the umbilic set of $S$ upon integrating the square root of its Hopf differential.
\newsubhead{Example II - Surfaces of constant negative curvature in $\Bbb{R}^3$}[SurfacesOfConstantNegativeCurvature]
Since $\Bbb{R}^3$ cannot be expressed as the level set of a non-degenerate bilinear form over $\Bbb{R}^4$, it is necessary to modify slightly our construction of $M$ and $W$. Although we could treat $\Bbb{R}^3$ as the limit of hyperbolic space as the length of the time direction tends to infinity, it is simpler to proceed more directly as follows. Let $\bil$ be a symmetric, non-degenerate bilinear form over $\Bbb{R}^3$. We define
$$
M := \{ (x,y)\in\Bbb{R}^3\times\Bbb{R}^3\ |\ b(y,y)\neq 0 \}\ ,
$$
and, for all $(x,y)\in M$, we define
$$
W_{(x,y)} := \{ (\xi,\mu)\in\Bbb{R}^3\times\Bbb{R}^3\ |\ b(\xi,x)=b(\mu,x)=0 \}\ .
$$
The Clifford structure is then constructed as before, with the pseudo-involutions $I^\eta$, $J^\eta$ and $K$ defined as in \eqnref{ExplicitFormulaOfIJK}. The reader may verify that the results of Sections \subheadref{TheCliffordBundle} and \subheadref{IntroducingBilegendrianSurfaces} continue to apply in this case.
\par
As before, let $(S,e)$ be an immersed surface in $\Bbb{R}^3$, let $\nu_e$ denote its unit normal vector field, let $\opI_e$, $\opII_e$ and $\opIII_e$ denote its three fundamental forms, and let $A_e$ denote its shape operator. We now view $\nu_e$ as a function taking values in $\Bbb{S}^2$ and, for $k>0$, we define the $k$-{\emph Gauss lift} of $e$ by
$$
\hat{e}_k(x) := \bigg(e(x),\frac{1}{\sqrt{k}}\nu_e(x)\bigg)\ ,
$$
and we call $(S,\hat{e}_k)$ the $k$-{\emph Gauss lift} of $(S,e)$. We readily verify that \eqnref{DerivativeOfGaussLift} continues to hold in this context, and we obtain the following result.
\proclaim{Lemma \nextprocno}
\noindent For $k>0$, $(S,e)$ has constant extrinsic curvature equal to $(-k)$ if and only if its $k$-Gauss lift is $J^{-1}$-invariant.
\endproclaim
\proclabel{KGaussLiftIsInvariantII}
\remark[KGaussLiftIsInvariantII] As in Lemma \procref{KGaussLiftIsInvariant}, the sign of $J$ is the same as that of the curvature. When the surfaces in question are of mixed signature, the situation is reversed, and the sign of $J$ is the opposite of that of the curvature. We will not address this case in the present paper.
\medskip
\proof Indeed, $(S,e)$ has constant extrinsic curvature equal to $(-k)$ if and only if $\opDet(A_e)=(-k)$. By \eqnref{DerivativeOfGaussLift}, this holds if and only if the $k$-Gauss lift is $\omega^{-1}_k$-lagrangian, and the result now follows by Theorem \procref{InvariantPlanesAreBiLagrangian}.\qed
\medskip
\noindent The manifestation of Theorem \procref{FlatMetric} in the present context is as follows.
\proclaim{Theorem \nextprocno, {\bf Existence of an asymptotic Chebyshev net}}
\noindent If $(S,e)$ has constant extrinsic curvature equal to $(-k)$, then the metric
$$
h_e := \opI_e + \frac{1}{k}\opIII_e\eqnum{\nexteqnno[FlatMetricOfKSurfaceB]}
$$
is flat. Furthermore, the lines of curvature are geodesics of this metric.
\endproclaim
\proclabel{FlatMetricOfKSurfaceB}
\remark[FlatMetricOfKSurfaceB] In this case, there are no umbilic points and $h_e$ is everywhere positive-definite.
\medskip
\proof Indeed, by \eqnref{DerivativeOfGaussLift},
$$
\hat{e}_k^* g^{-\eta} = \opI_e + \frac{1}{k}\opIII_e\ ,
$$
and the result now follows by Theorem \procref{FlatMetric}.\qed
\medskip
It remains only to discuss the geometric interpretation of the metric constructed in Theorem \procref{FlatMetricOfKSurfaceB}. We first suppose that $(S,e)$ is simply-connected. We will say that $(S,e)$ is {\emph quasicomplete} whenever it is complete with respect to the metric $\opI_e+\opIII_e$. Note that this property follows trivially from completeness. When this holds, $S$ may be globally parametrized by $\Bbb{R}^2$ such that
$$
h_e = dx^2 + dy^2\ .
$$
It is now straightforward to show that the fundamental forms of $e$ are
$$
\opI_e = \pmatrix\opCos^2(\theta)\hfill&0\hfill\cr0\hfill&\opSin^2(\theta)\hfill\cr\endpmatrix\ ,\
\opII_e = \frac{\sqrt{k}}{2}\pmatrix\opSin(2\theta)\hfill&0\hfill\cr0\hfill&-\opSin(2\theta)\hfill\cr\endpmatrix\ ,\ \text{and}\
\opIII_e = k\pmatrix\opSin^2(\theta)\hfill&0\hfill\cr0\hfill&\opCos^2(\theta)\hfill\cr\endpmatrix\ .\eqnum{\nexteqnno[FFSNegCurv]}
$$
for some smooth function $\theta:\Bbb{R}^2\rightarrow]0,\pi/2[$. In particular, the diagonals in $\Bbb{R}^2$ are the asymptotic lines of $e$, so that this parametrization is none other than the asymptotic Chebyshev net of $S$, used in the proof of Hilbert's theorem concerning the non-existence of complete, negatively curved CEC surfaces in $\Bbb{R}^3$ (see \cite{Spivak}).
\newsubhead{Example III - Hilbert's theorem}[HilbertsTheorem]
Our initial aim in writing this paper was to review the proof of Hilbert's theorem as part of a survey of recent developments in the theory of CEC surfaces in $3$-dimensional space forms. Since the proof of this result is both interesting and elementary once the asymptotic Chebyshev net has been constructed, for the sake of completeness, we discuss it here.
\par
We will address a more general case than that treated in Section \subheadref{SurfacesOfConstantNegativeCurvature}. For $c\in\Bbb{R}$, we denote by $X^3_c$ the $3$-dimensional riemannian space-form of constant sectional curvature equal to $c$. For $k>0$, let $(S,e)$ be a complete CEC surface in $X^3_c$ of curvature equal to $(-k)$. Reasoning as in Section \subheadref{SurfacesOfConstantNegativeCurvature}, we show that $S$ is globally parametrized by $\Bbb{R}^2$ in such a manner that the first, second and third fundamental forms are given by \eqnref{FFSNegCurv}, for some smooth function $\theta:\Bbb{R}^2\rightarrow]0,\pi/2[$.
\proclaim{Lemma \nextprocno}
\noindent The function $\theta$ satisfies
$$
(\theta_{xx}-\theta_{yy})dxdy = (k-c)\opdArea_e\ ,\eqnum{\nexteqnno[CurvatureEquationNegativeCase]}
$$
where $\opdArea_e$ here denotes the area form of the first fundamental form of $e$.
\endproclaim
\proclabel{CurvatureEquationNegativeCase}
\proof Indeed, using the above notation, we define the oriented, orthonormal frame $(u_x,u_y)$ of $\opI_e$ by
$$
u_x := \frac{1}{\opCos(\theta)}\partial_x\ \text{and}\ u_y := \frac{1}{\opSin(\theta)}\partial_y\ .
$$
The commutator of this frame is given by
$$
[u_x,u_y] = \frac{1}{\opCos(\theta)\opSin(\theta)}\big(-\opCot(\theta)(\partial_x\theta)\partial_y-\opTan(\theta)(\partial_y\theta)\partial_x\big)\ ,
$$
and its connection form is therefore
$$
\alpha = [u_x,u_y]^\flat = -(\partial_y\theta)dx - (\partial_x\theta)dy\ ,
$$
where $\flat$ here denotes Berger's musical isomorphism. It follows by Gauss' equation that
$$
(\theta_{xx}-\theta_{yy})dxdy = -d\alpha = (k-c)\opdA_e\ ,
$$
as desired.\qed
\proclaim{Corollary \nextprocno}
\noindent The function $\theta$ satisfies the hyperbolic sine-Gordon equation
$$
\theta_{xx}-\theta_{yy} = \frac{(k-c)}{2}\opSin(2\theta)\ .\eqnum{\nexteqnno[HyperbolicSineGordon]}
$$
\endproclaim
\proclabel{HyperbolicSineGordon}
\proof Indeed, with the above notation, the area form of $\opI_e$ is
$$
\opdArea_e = \frac{1}{2}\opSin(2\theta)dxdy\ ,
$$
and the result follows by \eqnref{CurvatureEquationNegativeCase}.\qed
\medskip
Hilbert's theorem now follows as a straightforward consequence of Hazzidaki's lemma, which we now recall. Let $[\theta]_0$ denote the total oscillation of $\theta$, that is
$$
[\theta]_0 := \msup_{x\neq y\in\Bbb{R}^2}\left|\theta(x)-\theta(y)\right|\ .
$$
\proclaim{Lemma \nextprocno, {\bf Hazzidaki}}
\noindent Let $\theta:\Bbb{R}^2\rightarrow]0,\pi/2[$ be a smooth function, consider the metric
$$
g_\theta := \opCos^2(\theta)dx^2 + \opSin^2(\theta)dy^2\ ,\eqnum{\nexteqnno[HazzidakiMetric]}
$$
and let $\kappa_\theta$ denote its curvature. If $\kappa_\theta$ does not change sign, then
$$
\int_{\Bbb{R}^2}\left|\kappa_\theta\right|\opdA_\theta \leq 4[\theta]_0\ .\eqnum{\nexteqnno[HazzidakiInequality]}
$$
\endproclaim
\proclabel{HazzidakiInequality}
\proof Indeed, by \eqnref{CurvatureEquationNegativeCase},
$$
\kappa_\theta\opdArea_\theta = (\theta_{xx}-\theta_{yy})dxdy\ .
$$
We define the new variables $u:=x+y$ and $v:=x-y$, and we verify that
$$
\int_{\Bbb{R}^2}\left|\kappa_\theta\right|\opdA_\theta
=\int_{\Bbb{R}^2}\left|\theta_{xx}-\theta_{yy}\right|dxdy
=\int_{\Bbb{R}^2}2\left|\theta_{uv}\right|dudv\ .
$$
Since $\theta_{uv}$ does not change sign,
$$\eqalign{
\int_{\Bbb{R}^2}2\left|\theta_{uv}\right|dudv &= \mlim_{R\rightarrow\infty}2\left|\int_{[-R,R]^2}\theta_{uv}dudv\right|\cr
&=\mlim_{R\rightarrow\infty}2\left|\theta(R,R)-\theta(R,-R)-\theta(-R,R)+\theta(-R,-R)\right|\vphantom{\bigg(}\cr
&\leq 4[\theta]_0\ ,\vphantom{\bigg(}\cr}
$$
and the result follows.\qed
\proclaim{Theorem \nextprocno, {\bf Hilbert}}
\noindent If $(-k)<\opMin(0,-c)$, then there exists no complete CEC surface in $X^3_c$ of curvature equal to $(-k)$.
\endproclaim
\proclabel{HilbertsTheorem}
\proof Suppose the contrary. Let $(S,e)$ be a complete CEC surface in $X^3_c$ of curvature equal to $(-k)$. We may suppose that $S$ is simply-connected, and since it is complete, in particular, it is quasicomplete, so that, by Hazzidaki's theorem, it has finite total curvature. However, since $(c-k)<0$, up to rescaling by a constant factor, $(S,\opI_e)$ is isometric to $\Bbb{H}^2$, which has infinite total curvature. This is absurd, and the result follows.\qed
\medskip
\noindent It is worth noting that when $c>0$ and $-c<-k<0$, Hilbert's theorem still holds, but with a slightly different proof.
\proclaim{Theorem \nextprocno, {\bf Hilbert}}
\noindent If $c>0$ and $-c<k<0$, then there exists no complete CEC surface in $X^3_c$ of curvature equal to $(-k)$.
\endproclaim
\proclabel{HilbertsTheoremB}
\proof Suppose the contrary, and let $(S,e)$ be such a surface. We may suppose that $S$ is simply connected. Since it is complete, it is quasicomplete, so that, by the preceding disussion, it is globally parametrized by $\Bbb{R}^2$. However, since $S$ is complete and of constant, positive curvature, it is diffeomorphic to the sphere. This is absurd, and the result follows.\qed
\medskip
We conclude this section with the following intriguing corollary of Hilbert's theorem.
\proclaim{Theorem \nextprocno}
\noindent For all $\epsilon>0$, there exists no smooth solution $\theta:\Bbb{R}^2\rightarrow[\epsilon,\pi-\epsilon]$ of the hyperbolic sine-Gordon equation
$$
\theta_{xx}-\theta_{tt}-\opSin(\theta) = 0\ .\eqnum{\nexteqnno[NoSolutionOfSinGordon]}
$$
\endproclaim
\proclabel{NoSolutionOfSinGordon}
\noindent It is worthwhile reflecting on the physical meaning of this result. The $1$-dimensional sine-Gordon equation
$$
\phi_{tt} + \opSin(\phi) = 0
$$
models the oscillations of a classical pendulum (see, for example, \cite{Goldstein}). The corresponding analogue of Theorem \procref{NoSolutionOfSinGordon} would thus affirm the intuitively obvious fact that a classical pendulum must eventually approach its axis. The $2$-dimensional sine-Gordon equation models the motion of an analogous system which we may view, either as the limit of an infinite row of such pendula joined to their neighbours by springs, or, more intuitively, as a very long vertical blind. Theorem \procref{NoSolutionOfSinGordon} thus states, once again, that such a system must also eventually approach its axis. It is quite striking that the proof of such a seemingly reasonable fact should rest on something as abstract as the infinite total curvature of hyperbolic space!
\newhead{Bilegendrian surfaces in $\opUS$}[BilegendrianSurfacesOverTheSphere]
\newsubhead{Overview}[OverviewBilegendrianSurfaces]
We conclude this paper by studying bilegendrian surfaces in $\opUS$, which generalize flat surfaces in $\Bbb{S}^3$. The first known examples of flat surfaces in $\Bbb{S}^3$ were probably the Clifford torus, introduced by Clifford himself in \cite{Clifford} (see also \cite{Volkert}), followed by those constructed by Bianchi in \cite{Bianchi}, presented in a modern manner using Hopf fibrations by Spivak in \cite{Spivak} (see also the work \cite{Sasaki} of Sasaki). Spivak also presented a simple construction of compact, embedded flat surfaces using Hopf projections, which was later studied in greater detail by Pinkall in \cite{Pinkall}.
\par
In \cite{Yau}, Yau posed the problem of classifying all compact flat surfaces in $\Bbb{S}^3$. This problem was solved independently by Kitigawa in \cite{KitigawaII} and by Weiner in \cite{Weiner} using entirely different techniques. Indeed, whilst Kitigawa builds on Spivak's work, Weiner uses the factorizations of the Gauss maps of flat surfaces developed by Enomoto in \cite{EnomotoI} and \cite{EnomotoII}.
\par
Our approach will contain elements of \cite{KitigawaII} and \cite{Weiner}, but will differ from both. In particular, the approaches of both these papers break down when the projection onto the base space degenerates. Instead, we will show in Sections \subheadref{BilegendrianSurfaces} and \subheadref{Factorization} how the $J$-invariant condition naturally interacts with the Lie group structure of $\Bbb{S}^3$ to yield in Theorem \procref{FactorizationOfBilegendrians} canonical factorizations of complete bilegendrian immersions, allowing us, in particular, to construct all such surfaces. Furthermore, in Theorem \procref{DoublyPeriodicCase}, we factorize the period lattices of compact bilegendrian immersions, and, together with Theorem \procref{HolonomyCondition}, this will establish a general formula for constructing all such immersions. We note, in particular, that both of these factorizations follow in a straightforward manner from a general factorization criterion for Lie group-valued functions over $\Bbb{R}^2$, and are thus, in a sense, more algebraic in nature than geometric.
\par
Complete bilegendrian immersions in $\opUS$ strike us as more natural objects than complete flat immersions in $\Bbb{S}^3$. Indeed, as we will see in Theorem \procref{ClassificationOfCompleteImmersions}, every complete bilegendrian immersion in $\opUS$ is uniquely defined, up to isometries of $\Bbb{S}^3$, by its angle function $\theta:\Bbb{R}^2\rightarrow\Bbb{R}$, which is a smooth solution of the classical wave equation
$$
\partial_x^2\theta - \partial_y^2\theta = 0\ .\eqnum{\nexteqnno[IntroWaveEquation]}
$$
However, wheras every such solution of \eqnref{IntroWaveEquation} is the angle function of some complete, bilegendrian immersion in $\opUS$, it is still a hard and unsolved problem to fully determine necessary and sufficient conditions for such a solution of \eqnref{IntroWaveEquation} to be the angle function of a {\emph complete} flat immersion in $\Bbb{S}^3$.
\par
In the doubly-periodic case, the distinction is even more striking, as the family of compact bilegendrian immersions in $\opUS$ is notably larger than that of compact flat immersions in $\Bbb{S}^3$. Indeed, wheras any pair of closed curves in $\Bbb{S}^2$ bounding areas which are rational multiples of $2\pi$ can arise as the factors of some doubly periodic bilegendrian immersion in $\opUS$, by Theorem \procref{AreaIsInteger}, only those curves which bound areas which are {\emph integer} multiples of $2\pi$ are candidates to arise as factors of compact flat immersions in $\Bbb{S}^3$. This has interesting consequences. For example, for any fixed pair $(q_1,q_2)$ of rational numbers, we expect a rigidity result analogous to the result \cite{EnoKitWei} of Enomoto--Kitigawa--Weiner to hold for compact bilegendrian immersions with factors that are respectively $q_1$- and $q_2$-quasiperiodic. We propose to address this problem in later work.
\newsubhead{Quaternions}[Quaternions]
Our study of complete bilegendrian immersed surfaces in $\opUS$ will rely heavily on the structure of the algebra $\Bbb{H}$ of quaternions, whose basic properties we now review. With the terminology of Section \subheadref{CliffordAlgebras}, $\Bbb{H}$ is nothing other than the Clifford algebra of a positive-definite, symmetric bilinear form over a $2$-dimensional real vector space. That is, it is the real algebra generated by the identity together with elements $\qi$, $\qj$ and $\qk$ satisfying
$$
\qi^2 = \qj^2 = \qk^2 = -1\ \text{and}\ \qi\cdot\qj\cdot\qk = -1\ .
$$
Let $\Cal{R}$ and $\Cal{I}$ denote respectively its real and imaginary subspaces. We define the {\emph inner product} $\bil$ over $\Bbb{H}$ by
$$
\bil(x,y) := \Cal{R}(x\cdot\overline{y})\ .
$$
In this case, $\bil$ is positive-definite, so that the set of unit-length quaternions identifies with the unit sphere $\Bbb{S}^3$, that is
$$
\Bbb{S}^3 = \{ x\in\Bbb{H}\ |\ \|x\|_\bil^2 = 1 \}\ .
$$
Since length is multiplicative, $\Bbb{S}^3$ also inherits from $\Bbb{H}$ the structure of a Lie group, with Lie algebra $\Cal{I}$. In particular, the inverse of any element $x$ of $\Bbb{S}^3$ is just its quaternionic conjugate $\overline{x}$.
\par
We view $\opUS$ as a submanifold of $\Bbb{H}\times\Bbb{H}$ view the identification
$$
\opUS := \{ (x,y)\in\Bbb{H}\times\Bbb{H}\ |\ \bil(x,x)=\bil(y,y)=1\ ,\ \bil(x,y)=0 \}\ .
$$
In particular, $\opUS$ is a codimension $2$ submanifold of the manifold $M$ studied in Section \subheadref{TheCliffordBundle}. The subbundle $W$ of $TM$ restricts to a codimension $1$ subbundle of the tangent bundle of $\opUS$, which we also denote by $W$.
\proclaim{Lemma \& Definition \nextprocno}
\noindent For all $(x,y)\in\opUS$, and for all $z\in\langle x,y\rangle^\perp$,
$$
-z\cdot\overline{x}\cdot y = y\cdot\overline{x}\cdot z \in \langle x,y\rangle^\perp\ .\eqnum{\nexteqnno[DefinitionOfAQuatCaseI]}
$$
We thus define
$$
A\cdot z := A_{(x,y)}\cdot z := -z\cdot\overline{x}\cdot y = y\cdot\overline{x}\cdot z\ .\eqnum{\nexteqnno[DefinitionOfAQuatCaseII]}
$$
This map is a rotation by $\pi/2$ radians about the axis $y$ in the positive direction with respect to the canonical orientation of $\Bbb{S}^3$.
\endproclaim
\proclabel{DefinitionOfAQuatCase}
\proof We first show that $z\cdot \overline{x}\cdot y=-y\cdot \overline{x}\cdot z$. Since $x$, $y$ and $z$ are pairwise orthogonal to one-another, the elements $\overline{x}\cdot y$ and $\overline{x}\cdot z$ are imaginary, and also orthogonal to one-another. Since orthogonal, imaginary quaternions anticommute,
$$
z\cdot\overline{x}\cdot y = x\cdot(\overline{x}\cdot z)\cdot(\overline{x}\cdot y) = -x\cdot(\overline{x}\cdot y)\cdot(\overline{x}\cdot z) = -y\cdot \overline{x}\cdot z\ ,
$$
as desired.
\par
We now show that this product lies in $\langle x,y\rangle^\perp$. Indeed, by orthogonality,
$$
\overline{x}\cdot y = -\overline{y}\cdot x\ ,
$$
so that
$$
\bil(z\cdot\overline{x}\cdot y,x) = \Cal{R}(z\cdot \overline{x}\cdot y\cdot\overline{x})
=-\Cal{R}(z\cdot\overline{y}\cdot x\cdot\overline{x}) = -\Cal{R}(z\cdot\overline{y})=-\bil(z,y) = 0\ .
$$
Likewise
$$
\bil(z\cdot\overline{x}\cdot y,y) = \Cal{R}(z\cdot\overline{x}\cdot y\cdot\overline{y}) = \Cal{R}(z\cdot\overline{x}) = \bil(z,x) = 0\ ,
$$
so that $z\cdot\overline{x}\cdot y$ is indeed an element of $\langle x,y\rangle^\perp$, as asserted.
\par
It remains only to verify the geometric properties of $A$. However,
$$
A^2\cdot z = z\cdot\overline{x}\cdot y\cdot \overline{x}\cdot y = -z\cdot\overline{x}\cdot x\cdot\overline{y}\cdot y = -z\ ,
$$
so that $A$ is a positive pseudo-involution, that is, a rotation by $\pi/2$ radians about the axis $y$. Finally, setting $(x,y)=(1,\qk)$ and $z=\qi$ yields $A\cdot z=-\qi\cdot\qk = \qj$, so that $A$ indeed turns in the positive direction about $y$ at this point. The result for general $(x,y)$ then holds by continuity, and this completes the proof.\qed
\medskip
We now set $\eta=-1$ and will suppress this superscript throughout the sequel. In terms of $A$, the endomorphism fields $I$, $J$ and $K$ constructed in Section \subheadref{TheCliffordBundle} are given by
$$
I = \pmatrix 0\hfill& \opId\hfill\cr \opId\hfill& 0\hfill\cr\endpmatrix\ ,\
J = \pmatrix 0\hfill& -A\hfill\cr A\hfill& 0\hfill\cr\endpmatrix\ ,\ \text{and}\
K = \pmatrix A\hfill& 0\hfill\cr 0\hfill& -A\hfill\cr\endpmatrix\ .\eqnum{\nexteqnno[DefinitionOfIJKForFlatSurfaces]}
$$
In particular,
$$
J^2 = \opId\ .
$$
\par
Finally, we identify $\Bbb{S}^2$ with the set of unit-length imaginary quaternions, that is
$$
\Bbb{S}^2 := \{ x\in\Cal{I}\ |\ \|x\|_\bil^2 = 1 \}\ ,
$$
and we view $\opT\Bbb{S}^2$ as a submanifold of $\Cal{I}\times\Cal{I}$ via the identification
$$
\opT\Bbb{S}^2 := \{ (\xi,\mu)\ |\ \|\xi\|=1\ ,\ \langle\xi,\mu\rangle=0\}\ .
$$
We recall the {\emph adjoint action} $\opad:\Bbb{S}^3\rightarrow\opSO(\Cal{I})$, given by
$$
\opad(x)\xi := x\cdot\xi\cdot x^{-1}\ .
$$
For any unit-length imaginary $x$, we define the {\emph Hopf projection} $\pi_x:\Bbb{S}^3\rightarrow\Bbb{S}^2$ by
$$
\pi_x(y) := \opad(y)x\ .
$$
For all such $x$, the preimage of $x$ under the Hopf projection is the $1$-parameter subgroup
$$
\Bbb{S}^1_x:=\Bbb{S}^3\minter\langle 1,x\rangle\ .
$$
In particular, $\pi_x$ projects to a diffeomorphism from $\Bbb{S}^3/\Bbb{S}^1_x$ to $\Bbb{S}^2$, which we also denote by $\pi_x$.
\newsubhead{Algebraic properties of bilegendrian surfaces}[BilegendrianSurfaces]
Let $(S,\phi)$ be a complete, simply-connected, bilegendrian immersed surface in $\opUS$, let $C$ and $\hat{C}$ denote respectively its cubic and complementary cubic forms, and let $X,Y:S\rightarrow\Bbb{S}^3$ denote its two components, that is
$$
X := \pi_1\circ\phi\ \text{and}\ Y := \pi_2\circ\phi\ .
$$
By Theorem \procref{FlatMetric}, we may suppose that $S=\Bbb{R}^2$, and that
$$
\phi^*\hat{g} = 2(dx_1^2 + dx_2^2)\ .\eqnum{\nexteqnno[FlatMetricOfBilegImm]}
$$
We may suppose furthermore that the diagonals $(\partial_1+\partial_2)/2$ and $(\partial_1-\partial_2)/2$ are principal vectors. Since
$$
J\cdot(\partial_1+\partial_2)/2 = (\partial_1-\partial_2)/2\ ,
$$
it follows that $\partial_1$ and $\partial_2$ are eigenvectors of $J$ with respective eigenvalues $(+1)$ and $(-1)$. In particular, since $J$ is $g$-antisymmetric, for each $i$,
$$
\phi^*g(\partial_i,\partial_i) = 0\ ,
$$
so that
$$
\bil(\partial_i X,\partial_i X) = \frac{1}{2}(g+\hat{g})(\partial_i\phi,\partial_i\phi) = 1\ .\eqnum{\nexteqnno[BLengthOfIthVector]}
$$
Finally, since $J$ is $C$-symmetric,
$$
C(\partial_1,\partial_2,\partial_2) = C(\partial_2,\partial_1,\partial_1) = 0\ .\eqnum{\nexteqnno[COverDiagB]}
$$
Likewise, by \eqnref{VanishingValuesOfCCF}, for every vector field $\xi$, and for each $i$,
$$
\hat{C}(\xi,\partial_i,\partial_i) = 0\ .\eqnum{\nexteqnno[COverDiagC]}
$$
\proclaim{Lemma \nextprocno}
\noindent The first derivatives of $X$ and $Y$ satisfy the following algebraic relations.
$$
\partial_1(Y\cdot X^{-1}) = \partial_2(X^{-1}\cdot Y) = 0\ .\eqnum{\nexteqnno[NormalIsLeftAndRightTransported]}
$$
\endproclaim
\proclabel{NormalIsLeftAndRightTransported}
\proof It suffices to prove the first identity, as the proof of the second is identical. By the product and chain rules,
$$
\partial_1(Y\cdot X^{-1}) = (\partial_1 Y - Y\cdot X^{-1}\cdot\partial_1 X)\cdot X^{-1}\ .
$$
However, since $\partial_1$ is an eigenvector of $J$ with eigenvalue $(+1)$, by \eqnref{DefinitionOfAQuatCaseII} and \eqnref{DefinitionOfIJKForFlatSurfaces},
$$
\partial_1 Y = A\cdot\partial_1 X = Y\cdot X^{-1}\cdot\partial_1 X\ ,
$$
so that
$$
\partial_1(Y\cdot X^{-1}) = 0\ ,
$$
as desired.\qed
\proclaim{Lemma \nextprocno}
\noindent The second derivatives of $X$ and $Y$ satisfy the following algebraic relations.
$$
\partial_2(X^{-1}\cdot(\partial_1 X)) = \partial_1((\partial_2 X)\cdot X^{-1}) = 0\ .\eqnum{\nexteqnno[ProductCriterionIsSatisfied]}
$$
\endproclaim
\proclabel{ProductCriterionIsSatisfied}
\proof It suffices to prove the first identity, as the proof of the second is identical. Since $X$ takes values in $\Bbb{S}^3$, $X^{-1}\cdot(\partial_1 X)$ takes values in $\Cal{I}$. Since $(X^{-1}\cdot(\partial_1 X),X^{-1}\cdot Y,X^{-1}\cdot Y\cdot X^{-1}\cdot(\partial_1 X))$ is an orthonormal basis of $\Cal{I}$, it suffices to show that the inner-products of $\partial_2(X^{-1}\cdot(\partial_1 X))$ with each of these elements vanish.
\par
By \eqnref{BLengthOfIthVector}, $\partial_1 X$ has constant unit length, and therefore so too does $X^{-1}\cdot(\partial_1 X)$. It follows that
$$
\bil(\partial_2(X^{-1}\cdot(\partial_1X)),X^{-1}\cdot(\partial_1X)) = \frac{1}{2}\partial_2\bil(X^{-1}\cdot(\partial_1X),X^{-1}\cdot(\partial_1X)) = 0\ ,
$$
as desired.
\par
By definition, $Y$ is orthogonal to $\partial_1 X$, so that $X^{-1}\cdot Y$ is orthogonal to $X^{-1}\cdot(\partial_1X)$. By \eqnref{NormalIsLeftAndRightTransported}, it follows that
$$
\bil(\partial_2(X^{-1}\cdot(\partial_1 X)),X^{-1}\cdot Y)=-\bil(X^{-1}\cdot(\partial_1 X),\partial_2(X^{-1}\cdot Y))=0\ ,
$$
as desired.
\par
Finally, by the chain and product rules,
$$\eqalign{
\bil(\partial_2(X^{-1}\cdot(\partial_1X)),X^{-1}\cdot Y\cdot X^{-1}\cdot(\partial_1X))
&=-\bil(X^{-1}\cdot(\partial_2X)\cdot X^{-1}\cdot(\partial_1X),X^{-1}\cdot Y\cdot X^{-1}\cdot(\partial_1 X))\cr
&\qquad\qquad+\bil(X^{-1}\cdot(\partial_2\partial_1X),X^{-1}\cdot Y\cdot X^{-1}\cdot(\partial_1 X))\ .\cr}
$$
Since $Y$ and $\partial_2X$ are orthogonal to one-another, and since multiplications on the left and on the right by unit quaternions yield isometries of $\Bbb{H}$, the first term on the right-hand side vanishes. Likewise, by \eqnref{DefinitionOfIJKForFlatSurfaces}, the second term is equal to
$$
\bil(D_{\partial_2}(\partial_1 X),A\cdot\partial_1 X)=\frac{1}{2}(g+\hat{g})(\nabla_{\partial_2}(\partial_1\phi),K\cdot\partial_1\phi)\ .
$$
Since $K=IJ$, and since $\partial_1$ is an eigenvector of $J$ with eigenvalue $(+1)$, this is equal to
$$
\frac{1}{2}(g+\hat{g})(\nabla_{\partial_2}(\partial_1\phi),I\cdot\partial_1\phi)
=\frac{1}{2}C(\partial_2,\partial_1,\partial_1) + \frac{1}{2}\hat{C}(\partial_2,\partial_1,\partial_1)\ ,
$$
which vanishes by \eqnref{COverDiagB} and \eqnref{COverDiagC}. This completes the proof.\qed
\newsubhead{Factorization}[Factorization]
We now use the algebraic identities of the preceding section to derive factorizations of complete bilegendrian immersions and, in the doubly-periodic case, of their period lattices. We first recall some further elementary definitions from the theory of Lie groups. We define the {\emph left} and {\emph right translations} $T_l,T_r:\Bbb{S}^3\times\Bbb{H}\rightarrow\Bbb{S}^3\times\Bbb{H}$ by
$$
T_l(x,\xi) := (x,x\cdot\xi)\ \text{and}\ T_r(x,\xi) := (x,\xi\cdot x)\ .
$$
Both $T_l$ and $T_r$ restrict to bundle isomorphisms from $\Bbb{S}^3\times\Cal{I}$ into $\opT\Bbb{S}^3$. For every imaginary quaternion $\xi$, we define its {\emph left} and {\emph right translates} respectively by
$$
\xi_l(x) := T_l(x,\xi)\ \text{and}\ \xi_r(x) := T_r(x,\xi)\ .
$$
Finally, when $\xi$ is non-zero, we say that a smooth curve $\gamma:I\rightarrow\Bbb{S}^3$ is {\sl left $\xi$-horizontal} (resp. {\sl right $\xi$-horizontal}) whenever $\dot{\gamma}$ is everywhere orthogonal to $\xi_l\circ\gamma$ (resp. $\xi_r\circ\gamma$).
\proclaim{Theorem \nextprocno, {\bf Factorization of immersions}}
\noindent Let $\phi:\Bbb{R}^2\rightarrow\opUS$ be a complete bilegendrian immersion. There exist a pair of orthogonal elements $a,b\in\Bbb{S}^3$, a right $(\overline{a}\cdot b)$-horizontal arc-length parametrized curve $\gamma_1:\Bbb{R}\rightarrow\Bbb{S}^3$, and a left $(\overline{b}\cdot a)$-horizontal arc-length parametrized curve $\gamma_2:\Bbb{R}\rightarrow\Bbb{S}^3$, such that $\gamma_1(0)=\gamma_2(0)=1$, and
$$
\phi(x_1,x_2) = (\gamma_2(x_2)\cdot a\cdot\gamma_1(x_1),\gamma_2(x_2)\cdot b\cdot\gamma_1(x_1))\ .\eqnum{\nexteqnno[FactorizationOfBilegendrians]}
$$
Furthermore, $a$, $b$, $\gamma_1$ and $\gamma_2$ are unique. Conversely, every function of this form is a complete, bilegendrian immersion.
\endproclaim
\proclabel{FactorizationOfBilegendrians}
In the case where the bilegendrian surface projects to a flat immersed surface in $\Bbb{S}^3$, two distinct factorizations are already known from the work \cite{KitigawaII} of Kitigawa and \cite{Weiner} of Weiner. It is worth reviewing how these factorizations relate to our own. We first address Kitigawa's construction. Here, the first component of \eqnref{FactorizationOfBilegendrians} identifies with the factorization described at the beginning of Section $4$ of \cite{KitigawaII}. Note, however, that whilst the factors are constructed in \cite{KitigawaII} by modifying horizontal lifts of curves in $\Bbb{S}^2$ (see Section $3$ of that paper), in the present case, the factors are obtained directly as horizontal lifts without modification upon correctly choosing the Hopf projections, as we will see presently.
\par
We now discuss Weiner's construction, which uses the following ideas developed by Enomoto in \cite{EnomotoI} and \cite{EnomotoII}. Let $\opGr^+(2,4)$ denote the grassmannian of oriented linear planes in $\Bbb{R}^4$. In the case where $X$ is an immersion, Enomoto defines its {\emph Gauss map} to be the function $G:\Bbb{R}^2\rightarrow\opGr^+(2,4)$ which maps every point $(x_1,x_2)$ of $\Bbb{R}^2$ to the tangent plane of $X$ at this point. In order to express Enomoto's Gauss map in terms of the curves $(\gamma_1,\gamma_2)$ of the factorization \eqnref{FactorizationOfBilegendrians}, it is necessary to recall how $\opGr^+(2,4)$ is described as a quotient of $\Bbb{S}^3\times\Bbb{S}^3$. We denote $P_0:=\langle 1,\qi\rangle$, and we furnish this plane with its canonical orientation. Let $\rho:\Bbb{S}^3\times\Bbb{S}^3\rightarrow\opSO(\Bbb{H})$ be the homomorphism given by
$$
\rho(x,y)\cdot z := x\cdot z\cdot y^{-1}\ .
$$
This homomorphism defines a transitive action of $\Bbb{S}^3\times\Bbb{S}^3$ on $\opGr^+(2,4)$, and its stabilizer of $P_0$ is $\Bbb{S}^1_\qi\times\Bbb{S}^1_\qi$. It follows that
$$
\opGr^+(2,4) = (\Bbb{S}^3/\Bbb{S}^1_\qi)\times(\Bbb{S}^3/\Bbb{S}^1_\qi) = \Bbb{S}^2\times\Bbb{S}^2\ .\eqnum{\nexteqnno[GrassmannianOfOrientedPlanes]}
$$
Let $\pi_\qi:\Bbb{S}^3\rightarrow\Bbb{S}^2$ denote the Hopf projection, as in Section \subheadref{Quaternions}. Observe now that the normal plane to $X$ at each point $(x_1,x_2)\in\Bbb{R}^2$ is spanned by $X(x_1,x_2)$ and $Y(x_1,x_2)$. It follows that, if $m,n\in\Bbb{S}^3$ are such that
$$
a = m\cdot\qj\cdot n^{-1}\ \text{and}\ b = m\cdot\qk\cdot n^{-1}\ ,
$$
then, up to reversal of orientation,
$$
G(x_1,x_2) = \langle\gamma_2(x_2)\cdot m\cdot 1\cdot n^{-1}\cdot\gamma_1(x_1),\gamma_2(x_2)\cdot m\cdot\qi\cdot n^{-1}\cdot\gamma_1(x_1)\rangle\ ,
$$
so that, with the identification \eqnref{GrassmannianOfOrientedPlanes},
$$
G(x_1,x_2) = (\pi_\qi(\gamma_2(x_2)\cdot m),\pi_\qi(\gamma_1(x_1)^{-1}\cdot n))\ .
$$
In this manner, we recover Weiner's and Enomoto's factorization from that of Theorem \procref{FactorizationOfBilegendrians}.
\par
Theorem \procref{FactorizationOfBilegendrians} is a consequence of the following general factorization result for smooth Lie group valued functions.
\proclaim{Lemma \nextprocno}
\noindent Let $G$ be a Lie group, and let $M:\Bbb{R}^2\rightarrow G$ be a smooth function. There exist $C\in G$, and smooth functions $A,B:\Bbb{R}\rightarrow G$ such that $A(0)=B(0)=\opId$, and
$$
M(x_1,x_2) = B(x_2)\cdot C\cdot A(x_1)\eqnum{\nexteqnno[FactorisationCriterionI]}
$$
if and only if
$$
\partial_2(M^{-1}\cdot(\partial_1 M)) = \partial_1((\partial_2 M)\cdot M^{-1}) = 0\ .\eqnum{\nexteqnno[FactorisationCriterionII]}
$$
Furthermore, the triplet $(A,B,C)$ is unique.
\endproclaim
\proclabel{FactorisationCriterion}
\proof Suppose first that $M$ is a product of the form \eqnref{FactorisationCriterionI}. Then
$$
\partial_2(M^{-1}\cdot (\partial_1 M)) = \partial_2(A(x_1)^{-1}\cdot\dot{A}(x_1)) = 0\ .
$$
The first identity of \eqnref{FactorisationCriterionII} follows, and the second follows in a similar manner.
\par
Suppose now that \eqnref{FactorisationCriterionII} holds, and denote
$$
\alpha(x_1) := M^{-1}\cdot(\partial_1 M)\ \text{and}\ \beta(x_2) := (\partial_2 M)\cdot M^{-1}\ .
$$
Define $A,B:\Bbb{R}\rightarrow G$ such that $A(0)=B(0)=\opId$, and
$$
\dot{A}=A\cdot\alpha\ \text{and}\ \dot{B}=\beta\cdot B\ .
$$
Since $B$ only depends on $x_2$,
$$
\partial_1(B^{-1}\cdot M\cdot A^{-1}) = B^{-1}\cdot((\partial_1 M)-M\cdot A^{-1}\cdot\dot{A})\cdot A^{-1}
=B^{-1}\cdot M\cdot(\alpha - A^{-1}\cdot\dot{A})\cdot A^{-1} = 0\ .
$$
Likewise
$$
\partial_2(B^{-1}\cdot M\cdot A^{-1}) = 0\ .
$$
The product $C:=B^{-1}\cdot M\cdot A^{-1}$ is therefore constant, and the result follows.\qed
\medskip
{\bf\noindent Proof of Theorem \procref{FactorizationOfBilegendrians}:\ }Indeed, let $X:=\pi_1\circ\phi$ and $Y:=\pi_2\circ\phi$ denote the two components of $\phi$. By Lemma \procref{ProductCriterionIsSatisfied}, $X$ satisfies the hypotheses of Lemma \procref{FactorisationCriterion}, and there therefore exists a unique element $a\in\Bbb{S}^3$ and unique smooth functions $\gamma_1,\gamma_2:\Bbb{R}\rightarrow\Bbb{S}^3$ such that $\gamma_1(0)=\gamma_2(0)=1$, and
$$
X(x_1,x_2) = \gamma_2(x_2)\cdot a\cdot\gamma_1(x_1)\ .
$$
Since $Y$ is normal to $X$, by Lemma \procref{NormalIsLeftAndRightTransported}, $\gamma_1$ and $\gamma_2$ are respectively right- and left-horizontal. Furthermore, up to isometries of $\Bbb{H}$, $\gamma_1$ and $\gamma_2$ are the respective restrictions of $X$ to the $x_1$- and $x_2$-axes, and it follows by \eqnref{BLengthOfIthVector} that these curves are arc-length parametrized.
\par
We now determine the form of $Y$. By symmetry, there exist a unique element $b\in\Bbb{S}^3$, and unique arc-length parametrized curves $\eta_1,\eta_2:\Bbb{R}\rightarrow\Bbb{S}^3$ such that $\eta_1(0)=\eta_2(0)=1$, and
$$
Y(x_1,x_2) = \eta_2(x_2)\cdot b\cdot\eta_1(x_1)\ .
$$
Applying \eqnref{NormalIsLeftAndRightTransported} with $x_2=0$ yields
$$
b\cdot\eta_1\cdot\gamma_1^{-1}\cdot\overline{a} = b\cdot\overline{a}\ ,
$$
so that
$$
\eta_1 = \gamma_1\ .
$$
Likewise
$$
\eta_2 = \gamma_2\ ,
$$
and \eqnref{FactorizationOfBilegendrians} follows.
\par
It remains to verify the geometric properties of $a$, $b$, $\gamma_1$ and $\gamma_2$. Upon evaluating \eqnref{FactorizationOfBilegendrians} at $(0,0)$, we see that $a$ is orthogonal to $b$. Since $\phi$ is an integral curve of $W$, applying \eqnref{NormalIsLeftAndRightTransported} with $x_2=0$ yields, for all $x_1$,
$$
a\cdot\dot{\gamma}_1(x_1) = \partial_1X(x_1,0)\perp Y(x_1,0) = Y(0,0)\cdot X^{-1}(0,0)\cdot X(x_1,0) = b\cdot\gamma_1(x_1)\ .
$$
It follows that
$$
\dot{\gamma}_1\cdot\gamma_1^{-1}\perp\overline{a}\cdot b\ ,
$$
so that $\gamma_1$ is right $(\overline{a}\cdot b)$-horizontal. In a similar manner, we show that $\gamma_2$ is left $(\overline{b}\cdot a)$-horizontal, as desired.
\par
Finally, we readily verify that any map of this form defines a bilegendrian immersion in $\opUS$, and this completes the proof.\qed
\medskip
We now show how the period lattice factorizes in the doubly-periodic case. We will say that a period lattice of a doubly-periodic immersion is {\emph maximal} whenever it is not contained in any larger period lattice.
\proclaim{Theorem \nextprocno, {\bf Factorization of period lattices}}
\noindent Let $\phi:\Bbb{R}^2\rightarrow\opUS$ be a complete bilegendrian immersion. If $\phi$ is doubly-periodic with maximal period lattice $\Lambda$, then there exist real numbers $p_1,p_2>0$, and rational numbers $q_1,q_2>0$ such that
$$
\Lambda = \left\{(m p_1,n p_2)\ |\ (m q_1 - n q_2)\in\Bbb{Z}\right\}\ ,\eqnum{\nexteqnno[DoublyPeriodicCaseI]}
$$
Furthermore, if $\gamma_1$ and $\gamma_2$ denote the factors of $\phi$, then, for all $x$,
$$
\gamma_1(x+p_1) = e^{-2\pi q_1 (b\cdot\overline{a})}\cdot\gamma_1(x)\ \text{and}\ \gamma_2(x+p_2) = \gamma_2(x)\cdot e^{2\pi q_2 (\overline{a}\cdot b)}\ .\eqnum{\nexteqnno[DoublyPeriodicCaseII]}
$$
In particular, each of these curves is also periodic, and their Hopf projections in $\Bbb{S}^2$ are periodic with respective periods $p_1$ and $p_2$.
\endproclaim
\proclabel{DoublyPeriodicCase}
\proof We use the factorization \eqnref{FactorizationOfBilegendrians}. Upon multiplying $\phi$ on the left by $\overline{a}$, and replacing $\gamma_2$ and $b$ respectively by $\gamma_2':=(\overline{a}\cdot\gamma_2\cdot a)$ and $b':=(\overline{a}\cdot b)$, we may suppose that $a=1$ and that $b$ is a unit imaginary quaternion. Let $\Bbb{S}^1_b$ denote the subgroup of $\Bbb{S}^3$ generated by $b$.
\par
We first show that, for all $(p_1,p_2)\in\Lambda$, there exists an element $\xi\in\Bbb{S}^1_b$ such that
$$
\gamma_1(p_1) = \gamma_2(p_2)^{-1} = \xi\ .
$$
Indeed, by double-periodicity of the first component,
$$
\gamma_2(p_2)\cdot\gamma_1(p_1) = X(p_1,p_2) = X(0,0) = 1\ ,
$$
so that
$$
\gamma_1(p_1) = \gamma_2(p_2)^{-1} =: \xi\ .
$$
Next, by double-periodicity of the second component,
$$
\gamma_2(p_2)\cdot b\cdot\gamma_1(p_1) = Y(p_1,p_2) = Y(0,0) = b\ ,
$$
so that
$$
b\cdot\xi = b\cdot\gamma_1(p_1) = \gamma_2(p_2)^{-1}\cdot b = \xi\cdot b\ .
$$
That is $\xi$ commutes with $b$. Since the centralizer of $b$ in $\Bbb{H}$ is $\langle 1,b\rangle$, it follows that $\xi\in\langle 1,b\rangle\cap\Bbb{S}^3=\Bbb{S}^1_b$, as asserted.
\par
We now claim that, for all $x$,
$$
\gamma_1(x+p_1) = \xi\cdot\gamma_1(x)\ \text{and}\ \gamma_2(x+p_2) = \gamma_2(x)\cdot\xi^{-1}\ .
$$
Indeed, it suffices to prove the first identity, as the proof of the second is identical. However,
$$
\gamma_1(x+p_1)^{-1}\cdot\dot{\gamma}_1(x+p_1) = X(x+p_1,0)^{-1}\cdot (\partial_1 X)(x + p_1,0)\ .
$$
By \eqnref{ProductCriterionIsSatisfied}, this is equal to
$$
X(x+p_1,p_2)^{-1}\cdot (\partial_1 X)(x+p_1,p_2) = X(x,0)^{-1}\cdot(\partial_1 X)(x,0) = \gamma_1(x)^{-1}\cdot\dot{\gamma}_1(x)\ .
$$
The curves $\gamma_1(\cdot+p_1)$ and $\xi\cdot\gamma_1$ therefore solve the same first-order ordinary differential equation with the same initial condition, so that, for all $x$,
$$
\gamma_1(x+p_1) = \xi\cdot\gamma_1(x)\ ,
$$
as asserted.
\par
Let $\Lambda_1$ and $\Lambda_2$ denote the respective images of $\Lambda$ under projection onto the first and second factors. By the preceding discussion, there exist homomorphisms $\xi_1:\Lambda_1\rightarrow\Bbb{S}^1_b$ and $\xi_2:\Lambda_2\rightarrow\Bbb{S}^1_b$ such that, for all $p_1\in\Lambda_1$, $p_2\in\Lambda_2$ and $x\in\Bbb{R}$,
$$
\gamma_1(x+p_1) = \xi_1(p_1)\cdot\gamma_1(x)\ \text{and}\ \gamma_2(x+p_2) = \gamma_2(x)\cdot\xi_2(p_2)^{-1}\ .
$$
We claim that $\Lambda_1$ and $\Lambda_2$ are discrete. Indeed, otherwise, without loss of generality, $\Lambda_1$ would be dense in $\Bbb{R}$, so that, by continuity, $\gamma_1$ would take values in $\Bbb{S}^1_b$. In particular, this would yield
$$
(\partial_1 X)(0,0) = \dot{\gamma}_1(0) = \pm b = \pm Y(0,0)\ ,
$$
which would be absurd, since $\phi$ is an integral curve of $W$. For each $i$, let $p_i$ be a positive generator of $\Lambda_i$, and let $q_i$ be such that
$$
e^{-2\pi q_i b} = \xi_i(p_i)\ .
$$
By definition, $\Lambda\subseteq\Lambda_1\times\Lambda_2$. Since $\Lambda$ is maximal, for $m,n\in\Bbb{Z}$, $mp_1+np_2\in\Lambda$ if and only if
$$
e^{2\pi n q_2 b}\cdot e^{-2\pi m q_1 b} = 1\ ,
$$
which in turn holds if and only if $m q_1 - n q_2\in\Bbb{Z}$. Finally, since $\Lambda$ is two-dimensional, there exists a non-singular integer-valued matrix $M$ such that
$$
M\cdot (q_1,q_2)^t \in\Bbb{Z}^2\ ,
$$
from which it follows that $q_1$ and $q_2$ are both rational. This completes the proof.\qed
\newsubhead{Quasiperiodic horizontal curves in $\Bbb{S}^3$}[UnitTorsionCurvesAndTheHopfFibration]
We now describe a general construction of doubly-periodic bilegendrian immersions in $\opUS$. Throughout this section, we fix a unit-length quaternion $\xi$, we denote by $\pi_\xi$ its Hopf projection, as defined in Section \subheadref{Quaternions}, and we denote
$$
\alpha_\xi := \xi_l^\flat\ \text{and}\ H_\xi := \opKer(\alpha_\xi)\ ,\eqnum{\nexteqnno[DefinitionOfAlphaAndW]}
$$
where $\flat$ here denotes Berger's musical isomorphism. We also restrict attention to left $\xi$-horizontal curves, as right $\xi$-horizontal curves are treated in an identical manner.
\par
Given a real number $q$, we say that a left $\xi$-horizontal curve $\gamma$ in $\Bbb{S}^3$ is $q$-{\emph quasiperiodic} whenever there exists $p\in\Bbb{R}$ such that, for all $x$,
$$
\gamma(x+p) = \gamma(x)\cdot e^{2\pi q\xi}\ .
$$
By Theorem \procref{DoublyPeriodicCase}, the classification of doubly-periodic bilegendrian immersions in $\opUS$ reduces to the classification of $q$-quasiperiodic curves in $\Bbb{S}^3$ for rational $q$, which we now address.
\proclaim{Lemma \nextprocno}
\noindent For all $\xi\in\Bbb{S}^3$, and for all real $q$, the image under $\pi_\xi$ of any $q$-quasiperiodic left $\xi$-horizontal curve is a periodic curve in $\Bbb{S}^2$.
\endproclaim
\proclabel{HopfProjectionOfQuasiperiodicCurve}
\proof Indeed, let $\gamma$ be a $q$-quasiperiodic left $\xi$-horizontal curve with period $p$. For all $x\in\Bbb{R}$,
$$
(\pi_\xi\circ\gamma)(x+p) = \gamma(x+p)\cdot\xi\cdot\gamma(x+p)^{-1}
=\gamma(x)\cdot e^{2\pi q\xi}\cdot\xi\cdot e^{-2\pi q\xi}\cdot\gamma(x)^{-1}
=\gamma(x)\cdot\xi\cdot\gamma(x)^{-1} = (\pi_\xi\circ\gamma)(x)\ ,
$$
and periodicity follows.\qed
\medskip
\noindent We now examine left $\xi$-horizontal lifts of periodic curves in $\Bbb{S}^2$. We first verify that they exist.
\proclaim{Lemma \nextprocno}
\noindent The Hopf projection $\pi_\xi$ defines a surjective riemannian submersion from $(\Bbb{S}^3,\bil)$ to $(\Bbb{S}^2,\bil/4)$. Furthermore, the kernel of $D\pi_\xi$ is the linear span of $\xi_l$, and its $\bil$-orthogonal complement is $H_\xi$.
\endproclaim
\proclabel{RiemannianSubmersion}
\proof Since multiplication on the left acts transitively by isometries on $\Bbb{S}^3$, and since $\opad$ acts transitively by isometries on $\Bbb{S}^2$, it suffices to prove the result at the identity. However, for all $\mu\in\Cal{I}$,
$$
D\pi_\xi(1)\cdot\mu = \frac{d}{dt}e^{t\mu}\cdot\xi\cdot e^{-t\mu}\bigg|_{t=0} = [\mu,\xi]\ .
$$
Suppose now that $\mu\in\langle\xi\rangle^\perp$. Since $\mu$ is imaginary, it anticommutes with $\xi$, so that
$$
D\pi_\xi(1)\cdot\mu = 2\xi\cdot\mu\ .
$$
Hence
$$
\|D\pi_\xi(1)\cdot\mu\|_\bil^2 = \|2\xi\cdot\mu\|_\bil^2 = 4\|\mu\|_\bil^2\ ,
$$
and the riemannian submersion property follows. Finally, observe that $\xi\in\opKer(D\pi_\xi)$. It follows by the rank nullity theorem that $\opKer(D\pi_\xi)$ coincides with the linear span of $\xi$, and this completes the proof.\qed
\medskip
\noindent By Lemma \procref{RiemannianSubmersion}, every smooth curve $\gamma$ in $\Bbb{S}^2$ lifts through the Hopf projection to a smooth, left $\xi$-horizontal curves in $\Bbb{S}^3$. We now establish a necessary and sufficient condition for any Hopf lift to be $q$-quasiperiodic. For this, we will require the concept of signed area bounded by closed curves in $\Bbb{S}^2$, which we now define. Note first that, since $\Bbb{S}^2$ is simply-connected, every piecewise smooth closed curve $\gamma:\Bbb{S}^1\rightarrow\Bbb{S}^2$ extends to a piecewise smooth function $\tilde{\gamma}:\overline{B}_1(0)\rightarrow\Bbb{S}^2$.
\proclaim{Lemma \& Definition \nextprocno}
\noindent If $\tilde{\gamma}_1,\tilde{\gamma}_2:\overline{B}_1(0)\rightarrow\Bbb{S}^2$ are two extensions of $\gamma$, then
$$
\int_{\overline{B}_1(0)}\tilde{\gamma}_1^*\opdArea = \int_{\overline{B}_1(0)}\tilde{\gamma}_2^*\opdArea\ \text{{\rm mod}}\ 4\pi\ ,\eqnum{\nexteqnno[DefinitionOfArea]}
$$
where $\opdArea$ here denotes the canonical area form of $\Bbb{S}^2$. We will call this integral the {\emph signed area} bounded by $\gamma$, and we denote it by $\opArea(\gamma)$.
\endproclaim
\proclabel{DefinitionOfArea}
\proof Viewing $\Bbb{S}^1$ as the equator of $\Bbb{S}^2$, we identify the lower hemisphere of $\Bbb{S}^2$ with $B_1(0)$, and the upper hemisphere with $B_1(0)$ with its orientation reversed. Let $\phi:\Bbb{S}^2\rightarrow\Bbb{S}^2$ be a piecewise smooth map equal to $\tilde{\gamma}_1$ in the lower hemisphere and equal to $\tilde{\gamma}_2$ in the upper hemisphere. Since $\phi^*\opdArea$ is an integral cohomology class over $\Bbb{S}^2$,
$$
\int_{\Bbb{S}^2}\phi^*\opdArea = 0\ \text{mod}\ 4\pi\ .
$$
Thus, since
$$
\int_{\Bbb{S}^2}\phi^*\opdArea = \int_{\overline{B}_1(0)}\tilde{\gamma}_1^*\opdArea - \int_{\overline{B}_1(0)}\tilde{\gamma}_2^*\opdArea\ ,
$$
the result follows.\qed
\medskip
\noindent It remains only to relate the holonomy of the $\xi$-horizontal lift to the signed area.
\proclaim{Lemma \nextprocno}
\noindent The contact form $\alpha_\xi$ satisfies
$$
\|d\alpha_\xi\wedge\alpha_\xi\|_\bil=2\ .\eqnum{\nexteqnno[ContactPropertyOfAlpha]}
$$
\endproclaim
\proclabel{ContactPropertyOfAlpha}
\proof Let $\mu$ and $\nu$ be unit-length imaginary quaternions. The exterior derivative of $\alpha_\xi$ is given by
$$
d\alpha_\xi(\mu_l,\nu_l) = D_{\mu_l}\alpha_\xi(\nu_l) - D_{\nu_l}\alpha_\xi(\mu_l) - \alpha_\xi([\mu_l,\nu_l])\ .
$$
Since the first two terms on the right-hand side vanish, and since $[\mu_l,\nu_l]=[\mu,\nu]_l$, it follows that
$$
d\alpha_\xi(\mu_l,\nu_l) = -\alpha_\xi([\mu,\nu]_l) = -\bil(\xi,[\mu,\nu])\ .
$$
In particular, if $\mu$, $\nu$ and $\xi$ are pairwise orthogonal, then
$$
[\mu,\nu]=\pm 2\xi\ ,
$$
from which it follows that
$$
\|d\alpha^\xi\wedge\alpha^\xi\|_\bil = 2\ ,
$$
as desired.\qed
\proclaim{Theorem \nextprocno, {\bf Quasiperiodicity of horizontal lifts}}
\noindent Let $\gamma$ be a periodic curve in $\Bbb{S}^2$ and with left $\xi$-horizontal lift $\hat{\gamma}$. If the signed area bounded by $\gamma$ is equal to $(-4\pi q)$ modulo $4\pi$, then $\hat{\gamma}$ is $q$-quasiperiodic.
\endproclaim
\proclabel{HolonomyCondition}
\noindent Combining Theorem \procref{HolonomyCondition} with Theorems \procref{FactorizationOfBilegendrians} and \procref{DoublyPeriodicCase}, we see that every doubly-periodic bilegendrian immersion in $\opUS$ may be constructed by the following ansatz. First choose two orthogonal unit quaternions $a$ and $b$. Next choose two arc-length parametrized periodic curves $\gamma_1$ and $\gamma_2$ in $\Bbb{S}^2$, bounding signed areas equal to rational multiples of $4\pi$, such that $\gamma_1(0)=\overline{a}\cdot b$ and $\gamma_2(0)=\overline{b}\cdot a$. Let $\hat{\gamma}_1$ denote the unique right $(\overline{a}\cdot b)$-horizontal lift of $\gamma_1$ such that $\gamma_1(0)=1$, and let $\hat{\gamma}_2$ denote the unique left $(\overline{b}\cdot a)$-horizontal lift of $\gamma_2$ such that $\gamma_2(0)=1$. The desired doubly-periodic bilegendrian immersion is now given by \eqnref{FactorizationOfBilegendrians}. Note, in particular, that since the set of area preserving deformations of any closed curve in $\Bbb{S}^2$ is infinite-dimensional, so too is the local deformation space of any doubly-periodic bilegendrian immersion in $\opUS$.
\medskip
{\bf\noindent Proof of Theorem \procref{HolonomyCondition}:\ }Without loss of generality, we may suppose that $\gamma(0)=\gamma(1)=\xi$. Let $q$ be such that $\hat{\gamma}(1)=\hat{\gamma}(0)\cdot e^{2\pi q\xi}$. By uniqueness of the horizontal lift, for all $x$,
$$
\hat{\gamma}(x+1) = \hat{\gamma}(x)\cdot e^{2\pi q\xi}\ .
$$
We now relate $q$ to the signed area bounded by $\gamma$. We define $\mu:I\rightarrow\Bbb{S}^3$ by
$$
\mu(t) := e^{2\pi tq\xi}\ .
$$
Since $\hat{\gamma}$ is horizontal, and since the restriction of $\alpha_\xi$ to $\Bbb{S}^1_\xi$ has constant unit length,
$$
2\pi q = \int_\mu\alpha_\xi = -\int_{\mu^{-1}}\alpha_\xi = -\int_{\mu^{-1}\cdot\hat{\gamma}}\alpha_\xi\ .
$$
We now view $\mu^{-1}\cdot\hat{\gamma}$ as a piecewise smooth function from $\Bbb{S}^1$ to $\Bbb{S}^3$. Since $\Bbb{S}^3$ is simply-connected, it extends to a piecewise smooth function $\phi:\overline{B}_1(0)\rightarrow\Bbb{S}^3$. By Stokes' theorem
$$
2\pi q = -\int_{\overline{B}_1(0)}\phi^*d\alpha_\xi\ .
$$
However, by Lemmas \procref{ContactPropertyOfAlpha} and \procref{RiemannianSubmersion},
$$
d\alpha_\xi = \frac{1}{2}\pi_\xi^*\opdArea\ ,
$$
so that
$$
4\pi q = -\int_{\overline{B}_1(0)}\phi^*\pi_\xi^*\opdArea = -\int_{\overline{B}_1(0)}(\pi_\xi\circ\phi)^*\opdArea\ .
$$
Since $\pi_\xi\circ\phi$ extends $\gamma$, it follows that
$$
4\pi q = -\opArea(\gamma)\ \text{mod}\ 4\pi\ ,
$$
as desired.\qed
\newsubhead{Curvature and the angle function}[Curvature]
We now return to the classical case of complete, flat immersions in $\Bbb{S}^3$. Indeed, although Theorems \procref{FactorizationOfBilegendrians}, \procref{DoublyPeriodicCase} and \procref{HolonomyCondition} together yield full classifications of complete and compact bilegendrian immersions in $\opUS$, it is clearly of interest to establish the conditions under which such immersions project to immersions in $\Bbb{S}^3$, that is, the conditions under which they arise as Gauss lifts of complete or compact flat surfaces in the sphere. Recall that these conditions were established by Kitigawa in \cite{KitigawaII} and Weiner in \cite{Weiner}. We now show how these results are recovered in the present framework.
\par
We first introduce the angle function. To this end, it will be useful to introduce the new coordinate system
$$
(u_1,u_2) := (x_1+x_2,x_1-x_2)\ .
$$
We denote the corresponding coordinate vector fields by $\hat{\partial}_1$ and $\hat{\partial}_2$, so that
$$
\hat{\partial}_1 = \frac{1}{2}(\partial_1+\partial_2)\ \text{and}\ \hat{\partial}_2 = \frac{1}{2}(\partial_1-\partial_2)\ .
$$
By construction, $\hat{\partial}_1$ and $\hat{\partial}_2$ are fields of principal vectors. Let $L_1$ and $L_2$ denote their corresponding principal lines and, for each $i$, let $e_i:\Bbb{R}^2\rightarrow\Bbb{S}^3$ be a global unit section of $L_i$ such that the quadruplet $(X,e_1,e_2,Y)$ is positively oriented. We define the {\sl angle function} $\theta:\Bbb{R}^2\rightarrow\Bbb{R}$ such that
$$
\hat{\partial}_1\phi = (\opCos(\theta)e_1,\opSin(\theta)e_1)\ .\eqnum{\nexteqnno[ProjectionOfFirstVector]}
$$
Since $\hat{\partial}_2\phi = J\cdot \hat{\partial}_1\phi$, it follows that
$$
\hat{\partial}_2\phi = (-\opSin(\theta)e_2,\opCos(\theta)e_2)\ .\eqnum{\nexteqnno[ProjectionOfSecondVector]}
$$
\proclaim{Lemma \nextprocno}
\noindent For every vector field $\xi$ over $\Bbb{R}^2$,
$$
C(\xi,\hat{\partial}_1,\hat{\partial}_1) = C(\xi,\hat{\partial}_2,\hat{\partial}_2) = -D_\xi\theta\ .\eqnum{\nexteqnno[ExplicitFormulaForC]}
$$
\endproclaim
\proclabel{ExplicitFormulaForC}
\noindent By symmetry, this suffices to determine every component of $C$. Furthermore, since $C$ identifies with the second fundamental form of $\phi$, the angle function $\theta$ thus suffices to determine the immersion $\phi$ up to postcomposition with isometries of $\Bbb{S}^3$. We will examine presently the necessary and sufficient conditions for $\theta$ to arise as the angle function of a bilegendrian immersion.
\proclaim{Corollary \nextprocno}
\noindent For every vector field $\xi$ over $\Bbb{R}^2$,
$$
C(\partial_1,\partial_1,\partial_1) = -4\partial_1\theta\ \text{and}\
C(\partial_2,\partial_2,\partial_2) = -4\partial_2\theta\ .\eqnum{\nexteqnno[COverDiagA]}
$$
\endproclaim
\proclabel{COverDiagA}
{\bf\noindent Proof of Lemma \procref{ExplicitFormulaForC}:\ } Since $e_1$ has constant unit length, for every vector field $\xi$ over $\Bbb{R}^2$,
$$
\bil(D_\xi e_1,e_1) = \frac{1}{2}D_\xi\bil(e_1,e_1) = 0.
$$
It follows that
$$\eqalign{
C(\xi,\hat{\partial}_1,\hat{\partial}_1)
&= g(D_\xi(\opCos(\theta)e_1,\opSin(\theta)e_1),I\cdot(\opCos(\theta)e_1,\opSin(\theta)e_1))\cr
&= \bil(D_\xi(\opCos(\theta)e_1),\opSin(\theta)e_1) - \bil(D_\xi(\opSin(\theta)e_1),\opCos(\theta)e_1)\cr
&= -\opSin^2(\theta)D_\xi\theta - \opCos^2(\theta)D_\xi\theta\cr
&= -D_\xi\theta\ ,\cr}
$$
as desired. Next, since $J$ is $C$-symmetric,
$$
C(\xi,\hat{\partial}_2,\hat{\partial}_2) = C(\xi,J\cdot\hat{\partial}_1,J\cdot\hat{\partial}_1) = C(\xi,\hat{\partial}_1,J^2\cdot\hat{\partial}_1) = C(\xi,\hat{\partial}_1,\hat{\partial}_1) = -D_\xi\theta\ ,
$$
and this completes the proof.\qed
\medskip
We now relate the above formulae to the geometry of asymptotic curves. To this end, we recall some elementary definitions of classical curve theory. We define a {\emph framed curve} in $\Bbb{S}^3$ to be a smooth function $\gamma:I\rightarrow\opSO(\bil)$ such that, for all $t$, the matrix $\gamma(t)^{-1}\cdot\dot{\gamma}(t)$ is tridiagonal. Recall that a matrix is said to be {\emph tridiagonal} whenever its only non-vanishing elements are those lying just above, on, and just below the diagonal. Given a framed curve $\gamma$, abusing notation, we denote
$$
\gamma=:(\gamma,T,N,B)\ ,
$$
and we call $T$, $N$ and $B$ respectively its {\emph tangent}, {\emph normal} and {\emph binormal}. We will say that the framed curve is {\emph arc-length parametrized} whenever its first component has this property. When this holds, we will suppose, furthermore, that $T=\dot{\gamma}$. Given an arc-length parametrized framed curve, we define its {\emph curvature} $\kappa$ and {\emph torsion} $\tau$ by
$$
\kappa := \bil(\dot{T},N)\ \text{and}\ \tau := \bil(\dot{N},B)\ .
$$
For each $i$, we now define
$$
\gamma_i(t) := X(t\partial_i)\ ,\ T_i(t) := (\partial_iX)(t\partial_i)\ ,\ N_i(t) := A\cdot T_i(t)\ ,\ \text{and}\ B_i(t) := Y(t\partial_i)\ .
$$
\proclaim{Lemma \nextprocno}
\noindent $\gamma_i:=(\gamma_i,T_i,N_i,B_i)$ is a framed curve in $\Bbb{S}^3$.
\endproclaim
\proclabel{FramedCurve}
\proof By definition of $A$, this quadruplet is positively oriented. We now show that $\gamma_i^{-1}\cdot\gamma_i$ is tridiagonal. Since this matrix is antisymmetric, and since $\dot{\gamma}_i=T_i$, it suffices to show that $\bil(\dot{T}_i,B_i)$ vanishes. However, since $Y$ is normal to $X$,
$$
\bil(\dot{T}_i,B_i) = \bil(D_{\partial_i}(\partial_iX),Y) = -\bil(\partial_iX,\partial_iY)\ .
$$
An elementary calculation involving \eqnref{ProjectionOfFirstVector} and \eqnref{ProjectionOfSecondVector} shows that this vanishes, and it follows that $\gamma_i^{-1}\cdot\gamma_i$ is tridiagonal, as desired.\qed
\proclaim{Lemma \nextprocno}
\noindent For each $i$, the curvature of $\gamma_i$ is given by
$$
\kappa_i = -2\epsilon_i \theta_i'\ .\eqnum{\nexteqnno[CurvatureOfAsymptoticCurves]}
$$
where $\epsilon_i$ here denotes the $J$-eigenvalue of $\partial_i$.
\endproclaim
\proclabel{CurvatureOfAsymptoticCurves}
\proof Indeed,
$$
\kappa_i = \bil(\dot{T}_i,N_i)
=\bil(D_{\partial_i}\partial_iX,A\cdot \partial_iX)
=\frac{1}{2}(g+\hat{g})(\nabla_{\partial_i}\partial_i\phi,K\cdot \partial_i\phi)\ .
$$
Since $K=IJ$, this yields
$$
\kappa_i = \frac{\epsilon_i}{2}(g+\hat{g})(\nabla_{\partial_i}\partial_i\phi,I\cdot \partial_i\phi)
= \frac{\epsilon_i}{2}C(\partial_i,\partial_i,\partial_i) + \frac{\epsilon_i}{2}\hat{C}(\partial_i,\partial_i,\partial_i)\ ,
$$
so that, by \eqnref{COverDiagC} and \eqnref{COverDiagA},
$$
\kappa_i = -2\epsilon_i\partial_i\theta\ ,
$$
as desired.\qed
\proclaim{Lemma \nextprocno}
\noindent For each $i$, the torsion of $\gamma_i$ is given by
$$
\tau_i = -\epsilon_i\ .\eqnum{\nexteqnno[TorsionOfAsymptoticCurves]}
$$
where $\epsilon_i$ here denotes the $J$-eigenvalue of $\partial_i$.
\endproclaim
\proclabel{TorsionOfAsymptoticCurves}
\remark[TorsionOfAsymptoticCurves] In fact, this follows immediately from Lemma \procref{NormalIsLeftAndRightTransported}. Indeed, it is an interesting exercise to show that a smooth framed curve in $\Bbb{S}^3$ is left (resp. right) horizontal if and only if it has unit positive (resp. negative) torsion (see, for example, Lemma $44$ of \cite{Spivak}).
\medskip
\proof Indeed, since $Y$ is normal to $X$,
$$
\tau_i = \bil(\dot{N}_i(t),B_i(t)) = \bil(D_{\partial_i}(A\cdot \partial_iX),Y)
= -\bil(A\cdot DX\cdot\partial_i,DY\cdot\partial_i)\ ,
$$
and the result now follows by an elementary calculation using \eqnref{DefinitionOfIJKForFlatSurfaces}, \eqnref{ProjectionOfFirstVector} and \eqnref{ProjectionOfSecondVector}.\qed
\proclaim{Lemma \nextprocno}
\noindent The angle function $\theta$ satisfies the wave equation
$$
\partial_1\partial_2\theta = 0\ .\eqnum{\nexteqnno[WaveOfThetaVanishes]}
$$
In particular, there exist smooth functions $\theta_1,\theta_2:\Bbb{R}\rightarrow\Bbb{R}$ such that
$$
\theta(x_1,x_2) = \theta_1(x_1) + \theta_2(x_2)\ .\eqnum{\nexteqnno[WaveOfThetaVanishesII]}
$$
\endproclaim
\proclabel{WaveOfThetaVanishes}
\remark[WaveOfThetaVanishes] In fact, this also follows immediately from \eqnref{FactorizationOfBilegendrians} and \eqnref{COverDiagA}. Indeed, since left and right multiplication by unit quaternions yield isometries of $\Bbb{H}$, the curvatures of horizontal (resp. vertical) lines are independent of the vertical (resp. horizontal) coordinate. However, the following, much longer, proof is of interest because it is also valid for more general bilegendrian surfaces.
\medskip
\proof It will be easier to prove the equivalent identity
$$
\hat{\partial}_1^2\theta - \hat{\partial}_2^2\theta = 0\ .
$$
By \eqnref{ExplicitFormulaForC} and symmetry of $C$,
$$
\hat{\partial}_1^2\theta = - D_{\hat{\partial}_1}C(\hat{\partial}_1,\hat{\partial}_1,\hat{\partial}_2) = -D_{\hat{\partial}_1}\omega_i(\nabla_{\hat{\partial}_2}(\hat{\partial}_1\phi),(\hat{\partial}_2\phi))\ .
$$
Since $\omega_i$ is covariant constant, this yields
$$
\hat{\partial}_1^2\theta = -\omega_i(\nabla_{\hat{\partial}_1}\nabla_{\hat{\partial}_2}(\hat{\partial}_1\phi),(\hat{\partial}_2\phi)) - \omega_i(\nabla_{\hat{\partial}_2}(\hat{\partial}_1\phi),\nabla_{\hat{\partial}_1}(\hat{\partial}_2\phi))\ .
$$
However,
$$
\omega_i(\nabla_{\hat{\partial}_2}(\hat{\partial}_1\phi),\nabla_{\hat{\partial}_1}(\hat{\partial}_2\phi)) = \hat{g}(\nabla_{\hat{\partial}_2}(\hat{\partial}_1\phi),\hat{I}\cdot\nabla_{\hat{\partial}_1}(\hat{\partial}_2\phi))\ .
$$
Since $\phi^*\hat{g}$ is the euclidean metric, the coordinate vector fields are $\phi^*\hat{g}$-covariant constant. The first argument of $\hat{g}$ on the right-hand side of the above identity is therefore $\hat{g}$-normal to $TS$, whilst the second argument is tangential. Their inner product therefore vanishes, and we conclude that
$$
\hat{\partial}_1^2\theta = -\omega_i(\nabla_{\hat{\partial}_1}\nabla_{\hat{\partial}_2}(\hat{\partial}_1\phi),(\hat{\partial}_2\phi))\ .
$$
Likewise,
$$
\hat{\partial}_2^2\theta = -\hat{\partial}_2 C(\hat{\partial}_1,\hat{\partial}_1,\hat{\partial}_2) = -\omega_i(\nabla_{\hat{\partial}_2}\nabla_{\hat{\partial}_1}(\hat{\partial}_1\phi),(\hat{\partial}_2\phi))\ ,
$$
and since $[\hat{\partial}_1,\hat{\partial}_2]=0$,
$$
\hat{\partial}_1^2\theta - \hat{\partial}_2^2\theta = -\omega_i(R_{\hat{\partial}_1\hat{\partial}_2}\hat{\partial}_1,\hat{\partial}_2)\ .
$$
By \eqnref{CurvatureFormula}, this vanishes, and the result follows.\qed
\medskip
We now conclude the classification of complete bilegendrian immersions in $\opUS$. Let $\hat{\Cal{B}}$ denote the space of smooth bilegendrian immersions $\phi:\Bbb{R}^2\rightarrow\opUS$ satisfying \eqnref{FlatMetricOfBilegImm}, furnished with the $C^\infty_\oploc$ topology. Let $\Cal{B}$ denote the quotient of this space under the action of postcomposition by elements of $\opSO(\Bbb{H})$. Let $C^\infty(\Bbb{R})$ denote the space of smooth, real-valued functions over $\Bbb{R}$, also furnished with the $C^\infty_\oploc$ topology, and let $\Bbb{S}^1$ denote the unit circle in $\Bbb{C}$. We construct the function $\Theta:\hat{\Cal{B}}\rightarrow\Bbb{S}^1\times C^\infty(\Bbb{R})\times C^\infty(\Bbb{R})$ as follows. Let $\phi:\Bbb{R}^2\rightarrow\opUS$ be a smooth, bilegendrian immersion satisfying \eqnref{FlatMetricOfBilegImm}, and consider the factorization \eqnref{FactorizationOfBilegendrians}. Upon post-composing with an element of $\opSO(\Bbb{H})$, we may suppose that $a=1$, $b=\qk$, and $\dot{\gamma}_1(0)=2\qi$, from which it follows that
$$
\dot{\gamma}_2(0) = 2e^{\theta_0\qk}\qi\ ,
$$
for some $\theta_0\in[0,2\pi[$. We now define $\Theta(\phi)$ by
$$
\Theta(\phi) := (e^{i\theta_0},\partial_1\theta,\partial_2\theta)\ ,
$$
where $\theta$ here denotes the angle function of $\phi$. We now define the set-valued function $\opI$ over $C^\infty(\Bbb{R})$ by
$$
\opI(f) := \opIm(F)\ ,
$$
where
$$
F(x) := \int_0^x f(y)dy\ .
$$
\proclaim{Theorem \nextprocno, {\bf Classification of complete immersions}}
\noindent $\Theta$ descends to a homeomorphism from $\Cal{B}$ into $\Bbb{S}^1\times C^\infty(\Bbb{R})\times C^\infty(\Bbb{R})$. Furthermore, $\phi\in\hat{\Cal{B}}$ projects to a flat immersed surface in $\Bbb{S}^3$ if and only if
$$
(\theta_0/2 + \opI(\theta_x) + \opI(\theta_y))\minter(\pi/2)\Bbb{Z} = \emptyset\ .\eqnum{\nexteqnno[ClassificationOfCompleteImmersions]}
$$
\endproclaim
\proclabel{ClassificationOfCompleteImmersions}
\proof We first show that $\Theta$ is invertible. Note first that every element of $\opSO(\Bbb{H})$ are uniquely determined by an orthonormal triple in $\Bbb{H}$, corresponding to the image of $(1,\qk,\qi)$. Given a framed curve $\gamma$, let $\kappa(\gamma)$ denote its curvature, and define $\hat{\Theta}:\hat{\Cal{B}}\rightarrow\opSO(\Bbb{H})\times\Bbb{S}^1\times C^\infty(\Bbb{R})\times C^\infty(\Bbb{R})$ by
$$
\hat{\Theta}(\phi) := ((a,b,\dot{\gamma}_1(0)/2),e^{i\theta_0},-\kappa(\gamma_1)/2,\kappa(\gamma_2)/2)\ ,
$$
where $(a,b,\gamma_1,\gamma_2)$ are as in the factorization \eqnref{FactorizationOfBilegendrians}, and $\theta_0$ denotes the angle between $\dot{\gamma}_1(0)$ and $\dot{\gamma}_2(0)$, measured in the positive direction. Since left and right horizontal curves in $\Bbb{S}^2$ are uniquely defined by their initial values, their initial directions, and their curvatures, $\hat{\Cal{B}}$ is a homeomorphism, and therefore so too is $\Cal{B}$.
\par
It remains only to establish when an element of $\hat{\Cal{B}}$ projects to an immersion in $\Bbb{S}^3$. However, by \eqnref{ProjectionOfFirstVector} and \eqnref{ProjectionOfSecondVector}, this holds whenever $\opCos(\theta)\opSin(\theta)$ never vanishes, that is, when
$$
\opIm(\theta)\cap(\pi/2)\Bbb{Z}=\emptyset\ .
$$
Thus, since
$$
\opIm(\theta) = \theta_0/2 + \opI(f) + \opI(g)\ ,
$$
the result follows.\qed
\medskip
\noindent The classification of doubly-periodic immersed bilegendrian surfaces in $\opUS$ does not lend itself to such a simple formulation. We will thus limit ourselves to the following result, which recovers Theorem $3.9$ of \cite{KitigawaII}.
\proclaim{Theorem \nextprocno}
\noindent Let $\phi:\Bbb{R}^2\rightarrow\opUS$ be a smooth, doubly-periodic bilegendrian immersion, and let $(a,b,\gamma_1,\gamma_2)$ denote its factorization as in \eqnref{FactorizationOfBilegendrians}. If $(\pi_\xi\circ\phi)$ is also an immersion, then each of $\gamma_1$ and $\gamma_2$ is either periodic, or $(1/2)$-quasiperiodic.
\endproclaim
\proclabel{AreaIsInteger}
\proof It suffices to prove the result for $\gamma_1$, as the proof for $\gamma_2$ is identical. Let $q\in\Bbb{Q}$ and $p\in\Bbb{R}$ be such that $\gamma_1$ is $q$-quasiperiodic with period $p$. Let $\theta$ denote the angle function of $\phi$, and let $\kappa_1$ denote the geodesic curvature of $(\pi_\xi\circ\gamma_1)$. By \eqnref{CurvatureOfAsymptoticCurves},
$$
\kappa_1 = -2(\partial_1\theta)\ .
$$
In particular,
$$
\int_0^p\kappa_1(x)dx = -2\int_0^p(\partial_1\theta)(x)dx = -2(\theta(p,0)-\theta(0,0))\ .
$$
This vanishes, for otherwise, by quasiperiodicity, $\theta$ would be unbounded, contradicting \eqnref{ClassificationOfCompleteImmersions}. By the Gauss--Bonnet Theorem, the curve $(\pi_\xi\circ\gamma_1)$ therefore bounds a signed area equal to an integer multiple of $2\pi$ modulo $4\pi$, and the result now follows by Theorem \procref{HolonomyCondition}.\qed
\inappendicestrue
\bigskip
\global\headno=0
\newhead{Calculating the curvature}[CalculatingTheCurvature]
\noindent Let $W$ denote the bundle defined in Section \subheadref{TheCliffordBundle}. We now determine certain components of its curvature. The calculation is rather technical and uninformative, and for this reason we have chosen to relegate it to this appendix.
\proclaim{Lemma \nextprocno}
\noindent At every point $(x,y)$ of $M$, and for every vector $\xi\in W_{(x,y)}$ such that $\pi_1(\xi)$ and $\pi_2(\xi)$ are colinear,
$$
\omega_i(R_{X,J^\eta X}X,J^{\eta} X) = \frac{-\epsilon\|(x,y)\|_g^2}{2\|x\|_{\bil}^2\|y\|_{\bil}^2}\hat{g}(X,I^\eta X)\hat{g}(X,X)\ .\eqnum{\nexteqnno[CurvatureFormula]}
$$
\endproclaim
\proclabel{CurvatureFormula}
\proof Indeed, let $\lambda_x$ and $\lambda_y$ denote respectively the $\bil$-lengths of $x$ and $y$. The normal bundle of $W$ is generated by the vector fields
$$
N_1 := \frac{1}{\lambda_x}(x,0),\ N_2 := \frac{1}{\lambda_y}(y,0),\ N_3 := \frac{1}{\lambda_x}(0,x),\ \text{and}\ N_4 := \frac{1}{\lambda_y}(0,y)\ .
$$
For $X$ and $Y$ vector fields taking values in $W$, the respective shape operators of these normals are
$$\eqalign{
B_1(X,Y) &= -\frac{1}{\lambda_x}\bil(\pi_1(X),\pi_1(Y))\ ,\cr
B_2(X,Y) &= -\frac{1}{\lambda_y}\bil(\pi_2(X),\pi_1(Y))\ ,\cr
B_3(X,Y) &= -\frac{\eta}{\lambda_x}\bil(\pi_1(X),\pi_2(Y))\ ,\ \text{and}\cr
B_4(X,Y) &= -\frac{\eta}{\lambda_y}\bil(\pi_2(X),\pi_2(Y))\ .\cr}
$$
We now suppose that $\pi_1(X)$ and $\pi_2(X)$ are colinear. Furthermore, we recall that
$$\eqalign{
\pi_1\circ J^\eta &= \eta A\circ\pi_2\ ,\cr
\pi_1\circ K &= A\pi_1\ ,\cr
\pi_2\circ J^\eta &= A\circ\pi_1\ ,\ \text{and}\cr
\pi_2\circ K &= -A\pi_2\ .\cr}
$$
Recall also that
$$
A^2 = -\epsilon\opId\ ,
$$
where $\epsilon$ here denotes the sign of $\bil$. For each $i$, let $\sigma_i$ denote the component of curvature arising from $B_i$. We obtain
$$\eqalign{
\sigma_1(X,J^\eta X,X,KX)
&=-\frac{1}{\|x\|_\bil^2}\bil(\pi_1(X),\pi_1(X))\bil(\pi_1(J^\eta X),\pi_1(KX))\cr
&=-\frac{\eta}{\|x\|_\bil^2}\bil(\pi_1(X),\pi_1(X))\bil(A\pi_2(X),A\pi_1(X))\cr
&=-\frac{\epsilon\eta}{\|x\|_\bil^2}\bil(\pi_1(X),\pi_1(X))\bil(\pi_1(X),\pi_2(X))\ .\cr}
$$
In a similar manner, we obtain
$$\eqalign{
\sigma_2(X,J^\eta,X,KX) &= -\frac{\epsilon}{\|y\|_\bil^2}\bil(\pi_1(X),\pi_1(X))\bil(\pi_2(X),\pi_1(X))\ ,\cr
\sigma_3(X,J^\eta,X,KX) &= \frac{\epsilon}{\|x\|_\bil^2}\bil(\pi_1(X),\pi_2(X))\bil(\pi_2(X),\pi_2(X))\ ,\ \text{and}\cr
\sigma_4(X,J^\eta,X,KX) &= \frac{\epsilon\eta}{\|y\|_\bil^2}\bil(\pi_2(X),\pi_2(X))\bil(\pi_1(X),\pi_1(X))\ .\cr}
$$
Combining these identities yields
$$\eqalign{
g(R_{X,J^\eta X}X,KX)
&= \bil(\pi_1(X),\pi_2(X))\bigg(-\frac{\epsilon\eta}{\|x\|_\bil^2}\hat{g}(X,X)-\frac{\epsilon}{\|y\|_\bil^2}\hat{g}(X,X)\bigg)\cr
&= \frac{-\epsilon\|(x,y)\|_g^2}{2\|x\|_\bil^2\|y\|_\bil^2}\hat{g}(X,I^\eta X)\hat{g}(X,X)\ ,\cr}
$$
as desired.\qed
\newhead{Bibliography}[Bibliography]
{\leftskip = 5ex \parindent = -5ex
\leavevmode\hbox to 4ex{\hfil \cite{BarBegZegh}}\hskip 1ex{Barbot T., B\'eguin F., Zeghib A., Prescribing Gauss curvature of surfaces in $3$-dimensional spacetimes: application to the Minkowski problem in the Minkowski space, {\sl Ann. Inst. Fourier}, {\bf 61}, no. 2, (2011), 511–-591}
\smallskip\leavevmode\hbox to 4ex{\hfil \cite{Bianchi}}\hskip 1ex{Bianchi L., Sulle superficie a curvatura nulla in geometria ellittica, {\sl Ann. Mat. Pura Appl.}, {\bf 24}, 93--129, (1896)}
\smallskip\leavevmode\hbox to 4ex{\hfil \cite{BonMonSchI}}\hskip 1ex{Bonsante F., Mondello G., Schlenker J. M., A cyclic extension of the earthquake flow I, {\sl Geom. Topol.}, {\bf 17}, no. 1, (2013), 157--234}
\smallskip\leavevmode\hbox to 4ex{\hfil \cite{BonMonSchII}}\hskip 1ex{Bonsante F., Mondello G., Schlenker J. M., A cyclic extension of the earthquake flow II, {\sl Ann. Sci. Ec. Norm. Sup\'er.}, {\bf 48}, no. 4, (2015), 811–-859}
\smallskip\leavevmode\hbox to 4ex{\hfil \cite{Clifford}}\hskip 1ex{Clifford W. K., Preliminary sketch of biquaternons, {\sl Proc. London Math. Soc.}, {\bf IV}, 381--395, (1873)}
\smallskip\leavevmode\hbox to 4ex{\hfil \cite{CliffordI}}\hskip 1ex{Clifford W. K., Applications of Grassmann's extensive algebra, {\sl Amer. J. Math.}, {\bf 1}, 350--358, (1878)}
\smallskip\leavevmode\hbox to 4ex{\hfil \cite{CliffordII}}\hskip 1ex{Clifford W. K., On the classification of geometric algebras, in {\sl Mathematical papers of William Kingdon Clifford}, Tucker R. (ed), MacMillan, London, (1882)}
\smallskip\leavevmode\hbox to 4ex{\hfil \cite{CruForGad}}\hskip 1ex{Cruceanu V., Fortuny P., Gadea P. M., A survey on paracomplex geometry, {\sl Rocky Mt. J. Math.}, {\bf 26}, no. 1, 83--115, (1996)}
\smallskip\leavevmode\hbox to 4ex{\hfil \cite{DajNom}}\hskip 1ex{Dajczer M., Nomizu K., On flat surfaces in $\Bbb{S}^3$ and $\Bbb{H}^3$, in {\sl Manifolds and Lie Groups}, Hano J., Morimoto A., Murakami S., Okamoto K., Ozeki H. (eds.), Progress in Mathematics, {\bf 14}, Birkh\"auser-Verlag, (1981)}
\smallskip\leavevmode\hbox to 4ex{\hfil \cite{EnomotoI}}\hskip 1ex{Enomoto K., The Gauss image of flat surfaces in $\Bbb{R}^4$, {\sl Kodai Math. J.}, {\bf 9}, 19--32, (1986)}
\smallskip\leavevmode\hbox to 4ex{\hfil \cite{EnomotoII}}\hskip 1ex{Enomoto K., Global properties of the Gauss image of flat surfaces in $\Bbb{R}^4$, {\sl Kodai Math. J.}, {\bf 10}, 272--284, (1987)}
\smallskip\leavevmode\hbox to 4ex{\hfil \cite{EnoKitWei}}\hskip 1ex{Enomoto K., Kitigawa Y., Weiner J. L., A rigidity theorem for the Clifford tori in $\Bbb{S}^3$, {\sl Proc. Amer. Math. Soc.}, {\bf 124}, 265--268, (1996)}
\smallskip\leavevmode\hbox to 4ex{\hfil \cite{FillSmi}}\hskip 1ex{Fillastre F., Smith G., Group actions and scattering problems in Teichm\"uller theory, in {\sl Handbook of group actions IV}, Advanced Lectures in Mathematics, {\bf 40}, (2018), 359--417}
\smallskip\leavevmode\hbox to 4ex{\hfil \cite{GalMarMil}}\hskip 1ex{G\'alvez J. A., M\'artinez A., Mil\'an F., Flat surfaces in the hyperbolic $3$-space, {\sl Math. Ann.}, {\bf 316}, no. 3, 419--435, (2000)}
\smallskip\leavevmode\hbox to 4ex{\hfil \cite{GalMirI}}\hskip 1ex{G\'alvez J. A., Mira P., Embedded isolated singularities of flat surfaces in hyperbolic $3$-space, {\sl Cal. Var. PDE.}, {\bf 24}, 239--260, (2005)}
\smallskip\leavevmode\hbox to 4ex{\hfil \cite{GalMirII}}\hskip 1ex{G\'alvez J. A., Mira P., Isometric immersions of $\Bbb{R}^2$ into $\Bbb{R}^4$ and perturbation of Hopf tori, {\sl Math. Zeit.}, {bf 266}, 207--227, (2010)}
\smallskip\leavevmode\hbox to 4ex{\hfil \cite{Garling}}\hskip 1ex{Garling D. J. H., {\sl Clifford algebras: an introduction}, Cambridge University Press, (2011)}
\smallskip\leavevmode\hbox to 4ex{\hfil \cite{Goldstein}}\hskip 1ex{Goldstein H. P., {\sl Classical mechanics}, Pearson, (2001)}
\smallskip\leavevmode\hbox to 4ex{\hfil \cite{HarveyLawson}}\hskip 1ex{Harvey F. R., Lawson H. B., Pseudoconvexity for the special lagrangian potential equation, {\sl Calc.Var. PDEs.}, {\bf 60}, (2021)}
\smallskip\leavevmode\hbox to 4ex{\hfil \cite{HitchinII}}\hskip 1ex{Hitchin N. J., The moduli space of special lagrangian submanifolds, {\sl, Annali della Scuola Normale Superiore de Pisa, Classe de Scienze 4${}^{\text{e}}$ s\'erie}, {\bf 25}, no. 3--4, 503--515, (1997)}
\smallskip\leavevmode\hbox to 4ex{\hfil \cite{HitchinI}}\hskip 1ex{Hitchin N. J., The moduli space of complex lagrangian submanifolds, {\sl Asian J. Math.}, {\bf 3}, no. 1, 77--92, (1999)}
\smallskip\leavevmode\hbox to 4ex{\hfil \cite{KitigawaII}}\hskip 1ex{Kitigawa Y., Periodicity of the asymptotic curves in flat tori in $\Bbb{S}^3$, {\sl J. Math. Soc. Japan}, {\bf 40}, 457--476, (1988)}
\smallskip\leavevmode\hbox to 4ex{\hfil \cite{KitigawaI}}\hskip 1ex{Kitagawa Y., Isometric deformations of flat tori in the $3$-sphere with nonconstant mean curvature, {\sl Tohoku Math. J.}, {\bf 52}, no. 2, 283--298, (2000)}
\smallskip\leavevmode\hbox to 4ex{\hfil \cite{KokUmeYam}}\hskip 1ex{Kokubu M., Umehara M., Yamada K., Flat fronts in hyperbolic $3$-space, {\sl Pac. J. Math.}, {\bf 216}, 149--175, (2004)}
\smallskip\leavevmode\hbox to 4ex{\hfil \cite{KokRossSajUmeYam}}\hskip 1ex{Kokubu M., Rossman W., Saji K., Umehara M., Yamada K., Singularities of flat fronts in hyperbolic $3$-space, {\sl Pac. J. Math.}, {\bf 221}, 303--351, (2005)}
\smallskip\leavevmode\hbox to 4ex{\hfil \cite{LabIII}}\hskip 1ex{Labourie F., M\'etriques prescrites sur le bord des vari\'et\'es hyperboliques de dimension $3$, {\sl J. Diff. Geom.}, {\bf 35}, no. 3, (1992), 609--626}
\smallskip\leavevmode\hbox to 4ex{\hfil \cite{LabI}}\hskip 1ex{Labourie F., Exemples de courbes pseudo-holomorphes en g\'eom\'etrie riemannienne, in {\sl Holomorphic curves in symplectic geometry}, Progress in Mathematics, {\bf 117}, Birkh\"auser, (1994)}
\smallskip\leavevmode\hbox to 4ex{\hfil \cite{LabII}}\hskip 1ex{Labourie F., Probl\`emes de Monge-Amp\`ere, courbes pseudo-holomorphes et laminations, {\sl Geom. Funct. Anal.}, {\bf 7}, (1997), 496--534}
\smallskip\leavevmode\hbox to 4ex{\hfil \cite{LabIV}}\hskip 1ex{Labourie F., Un lemme de Morse pour les surfaces convexes, {\sl Inv. Math.}, {\bf 141}, No 2, (2000), 239--297}
\smallskip\leavevmode\hbox to 4ex{\hfil \cite{Mira}}\hskip 1ex{Le\'on-Guzm\'an M. A., Mira P., Pastor J. A., The space of Lorentzian flat tori in anti de-Sitter $3$-space, {\sl Trans. AMS.}, {\bf 363}, no. 12, 6549--6573, (2011)}
\smallskip\leavevmode\hbox to 4ex{\hfil \cite{Lounesto}}\hskip 1ex{Lounesto P., {\sl Clifford algebras and spinors}, London Mathematical Society Lecture Note Series, {\bf 286}, (2001)}
\smallskip\leavevmode\hbox to 4ex{\hfil \cite{LychRub}}\hskip 1ex{Lychagin V. V., Rubtsov V. N., Local classification of Monge-Amp\`ere differential equations (Russian), {\sl Dokl. Akad. Nauk SSSR}, {\bf 272}, no. 1, (1983), 34--38}
\smallskip\leavevmode\hbox to 4ex{\hfil \cite{LychRubChekI}}\hskip 1ex{Lychagin V. V., Rubtsov V. N., Chekalov I. V., A classification of Monge-Amp\`ere equations, {\sl Ann. Sci. ENS.}, {\bf 26}, no. 3, 281--308, (1993)}
\smallskip\leavevmode\hbox to 4ex{\hfil \cite{McDuffSalamon}}\hskip 1ex{McDuff D., Salamon D., {\sl J-holomorphic curves and symplectic topology}, Colloquium Publications, {\bf 52}, Amer. Math. Soc., (2012)}
\smallskip\leavevmode\hbox to 4ex{\hfil \cite{Milnor}}\hskip 1ex{Milnor T. K., Harmonic maps and classical surface theory in minkowski $3$-space, {\sl Trans. AMS.}, {\bf 280}, no. 1, 161--185, (1983)}
\smallskip\leavevmode\hbox to 4ex{\hfil \cite{Pinkall}}\hskip 1ex{Pinkall U., Hopf tori in $\Bbb{S}^3$, {\sl Inv. Math.}, {\bf 81}, 379--386, (1985)}
\smallskip\leavevmode\hbox to 4ex{\hfil \cite{Rubtsov}}\hskip 1ex{Rubtsov V. N., Geometry of Monge-Amp\`ere structures, in {\sl  Nonlinear PDEs, their geometry, and applications}, Kycia R. A. et al. (eds.), Springer-Verlag, (2019)}
\smallskip\leavevmode\hbox to 4ex{\hfil \cite{SajUmeYam}}\hskip 1ex{Saji K., Umehara M., Yamada K., The geometry of fronts, {\sl Ann. Math.}, {\bf 169}, 491--529, (2009)}
\smallskip\leavevmode\hbox to 4ex{\hfil \cite{Sasaki}}\hskip 1ex{Sasaki S., On complete surfaces with Gaussian curvature zero in the $3$-sphere, {\sl Colloq. Math.}, {\bf 26}, 165--174, (1972)}
\smallskip\leavevmode\hbox to 4ex{\hfil \cite{SasakiNomizu}}\hskip 1ex{Sasaki S., Nomizu K., {\sl Affine differential geometry: geometry of affine immersions}, Cambridge Tracts in Mathematics, {\bf 111}, Cambridge University Press, (1994)}
\smallskip\leavevmode\hbox to 4ex{\hfil \cite{Sch}}\hskip 1ex{Schlenker J. M., Hyperbolic manifolds with convex boundary, {\sl Inv. Math.}, {\bf 163}, (2006), 109--169}
\smallskip\leavevmode\hbox to 4ex{\hfil \cite{SmiSLC}}\hskip 1ex{Smith G., Special Lagrangian curvature, Math. Annalen, 355, no. 1, 57--95, (2013)}
\smallskip\leavevmode\hbox to 4ex{\hfil \cite{SmiEW}}\hskip 1ex{Smith G., On an Enneper-Weierstrass-type representation of constant Gaussian curvature surfaces in 3-dimensional hyperbolic space, in {\sl Minimal surfaces: integrable systems and visualisation}, Hoffmann T., Kilian M., Leschke K., Martin G. (eds.), Springer Proceedings in Mathematics and Statistics, {\bf 349}, Springer-Verlag, (2021)}
\smallskip\leavevmode\hbox to 4ex{\hfil \cite{SmiMS}}\hskip 1ex{Smith G., M\"obius structures, hyperbolic ends and $k$-surfaces in hyperbolic space, in {\sl In the tradition of Thurston, Vol. II}, Ohshika K., Papadopoulos A. (eds.), Springer Verlag, (2022)}
\smallskip\leavevmode\hbox to 4ex{\hfil \cite{SmiAS}}\hskip 1ex{Smith G., On the asymptotic geometry of finite-type $k$-surfaces in three-dimensional hyperbolic space, to appear in {\sl J. Eur. Math. Soc.}}
\smallskip\leavevmode\hbox to 4ex{\hfil \cite{SmiQMAK}}\hskip 1ex{Smith G., Quaternions, Monge-Amp\`ere structures, and $k$-surfaces, arXiv:2210.02664}
\smallskip\leavevmode\hbox to 4ex{\hfil \cite{SmiQCKS}}\hskip 1ex{Smith G., On quasicomplete $k$-surfaces in $3$-dimensional space-forms, arXiv:2211.14868}
\smallskip\leavevmode\hbox to 4ex{\hfil \cite{Spivak}}\hskip 1ex{Spivak M., {\sl A comprehensive introduction to differential geometry, Vol. IV}, Publish or Perish, Berkeley, (1977)}
\smallskip\leavevmode\hbox to 4ex{\hfil \cite{Tamburelli}}\hskip 1ex{Tamburelli A., Prescribing metrics on the boundary of anti-de Sitter $3$-manifolds, {\sl Int. Math. Res. Not.}, no. 5, (2018), 1281--1313}
\smallskip\leavevmode\hbox to 4ex{\hfil \cite{Volkert}}\hskip 1ex{Volkert K., Space forms: a history, {\sl Bulletin of the Manifold Atlas}, (2013)}
\smallskip\leavevmode\hbox to 4ex{\hfil \cite{Weiner}}\hskip 1ex{Weiner J. L., Flat tori in $\Bbb{S}^3$ and their Gauss maps, {\sl Proc. London Math. Soc.}, {\bf 62}, 54--76, (1991)}
\smallskip\leavevmode\hbox to 4ex{\hfil \cite{Yau}}\hskip 1ex{Yau S. T., Submanifolds with constant mean curvature II, {\sl Amer. J. Math.}, {\bf 97}, 76--100, (1975)}
\par}
%
%
%
%
\enddocument